\documentclass[a4paper,twoside,11pt]{article}

\usepackage[left=2cm, right=2cm, bottom=3cm, top=3cm]{geometry}

\title{A Resolution of the Poisson Problem for Elastic Plates}
\date{\today}
\author{Francesca Da Lio\thanks{ETH Department Mathematik,R\"amistrasse 101, CH-8092 Z\"urich, Switzerland}
\and Francesco Palmurella*,
\and Tristan Rivi\`ere*}

\usepackage{amsmath,amsthm,amssymb,amsxtra} 
\usepackage{bbm} 
\usepackage{todonotes} 
\usepackage{enumerate} 
\usepackage{hyperref} 
\usepackage{mathrsfs}
\usepackage{framed} 
\usepackage{graphicx}

\numberwithin{equation}{section}

\newcommand{\R}{\mathbb{R}}
\newcommand{\C}{\mathbb{C}}
\newcommand{\N}{\mathbb{N}}
\newcommand{\Z}{\mathbb{Z}}

\newcommand{\dd}{\mathrm{d}} 

\newcommand{\e}{\mathrm{e}}  
\newcommand{\hau}{\mathcal{H}} 
\newcommand{\leb}{\mathcal{L}} 
\DeclareMathOperator{\divop}{div}

\DeclareMathOperator{\dist}{dist}
\DeclareMathOperator{\supp}{supp}
\DeclareMathOperator{\sign}{sign}

\DeclareMathOperator{\pv}{p.v.}
\DeclareMathOperator{\spanop}{span}

\DeclareMathOperator{\II}{I\!I}
\DeclareMathOperator{\Area}{Area}

\DeclareMathOperator{\Gr}{Gr}

\DeclareMathOperator{\Id}{Id}

\def\Xint#1{\mathchoice
        {\XXint\displaystyle\textstyle{#1}}%
        {\XXint\textstyle\scriptstyle{#1}}%
        {\XXint\scriptstyle\scriptscriptstyle{#1}}%
        {\XXint\scriptscriptstyle\scriptscriptstyle{#1}}%
        \!\int}
\def\XXint#1#2#3{{\setbox0=\hbox{$#1{#2#3}{\int}$}
\vcenter{\hbox{$#2#3$}}\kern-.5\wd0}}
\def\dashint{\Xint-}

\theoremstyle{plain}
\newtheorem{definition}{Definition}[section]
\newtheorem{theorem}[definition]{Theorem}
\newtheorem*{theorem*}{Theorem} 
\newtheorem{lemma}[definition]{Lemma}

\newtheorem{remark}[definition]{Remark}

\theoremstyle{remark}
\newtheorem{example}[definition]{Example}
\newenvironment{bfproof}[1][\proofname]
{\begin{proof}[\normalfont\bfseries #1]}{\end{proof}} 

\begin{document}
\maketitle

\begin{abstract}
The Poisson problem
consists in finding an immersed surface $\Sigma\subset\R^m$
minimising Germain's elastic energy
(known as Willmore energy in geometry)
with prescribed boundary,
boundary Gauss map and area.
This problem represents a non-linear model for the equilibrium state of thin, clamped elastic plates
originating from the work of S. Germain and S.D. Poisson in the early XIX century.
We present a solution to this problem
in the case of boundary data of class $C^{1,1}$ and when the boundary curve is simple and closed.
The minimum is realised by an immersed disk, possibly with a finite number
of branch points in its interior, which is of class $C^{1,\alpha}$ up to the boundary
for some $0<\alpha<1$, and whose Gauss map extends
to a map of class $C^{0,\alpha}$ up to the boundary.
\end{abstract}

\medskip


\setcounter{tocdepth}{1}
\tableofcontents

\section{Introduction}
\begin{figure}[h]
\label{fig:Poisson-excerpts}
\centering
\includegraphics[scale=0.5]{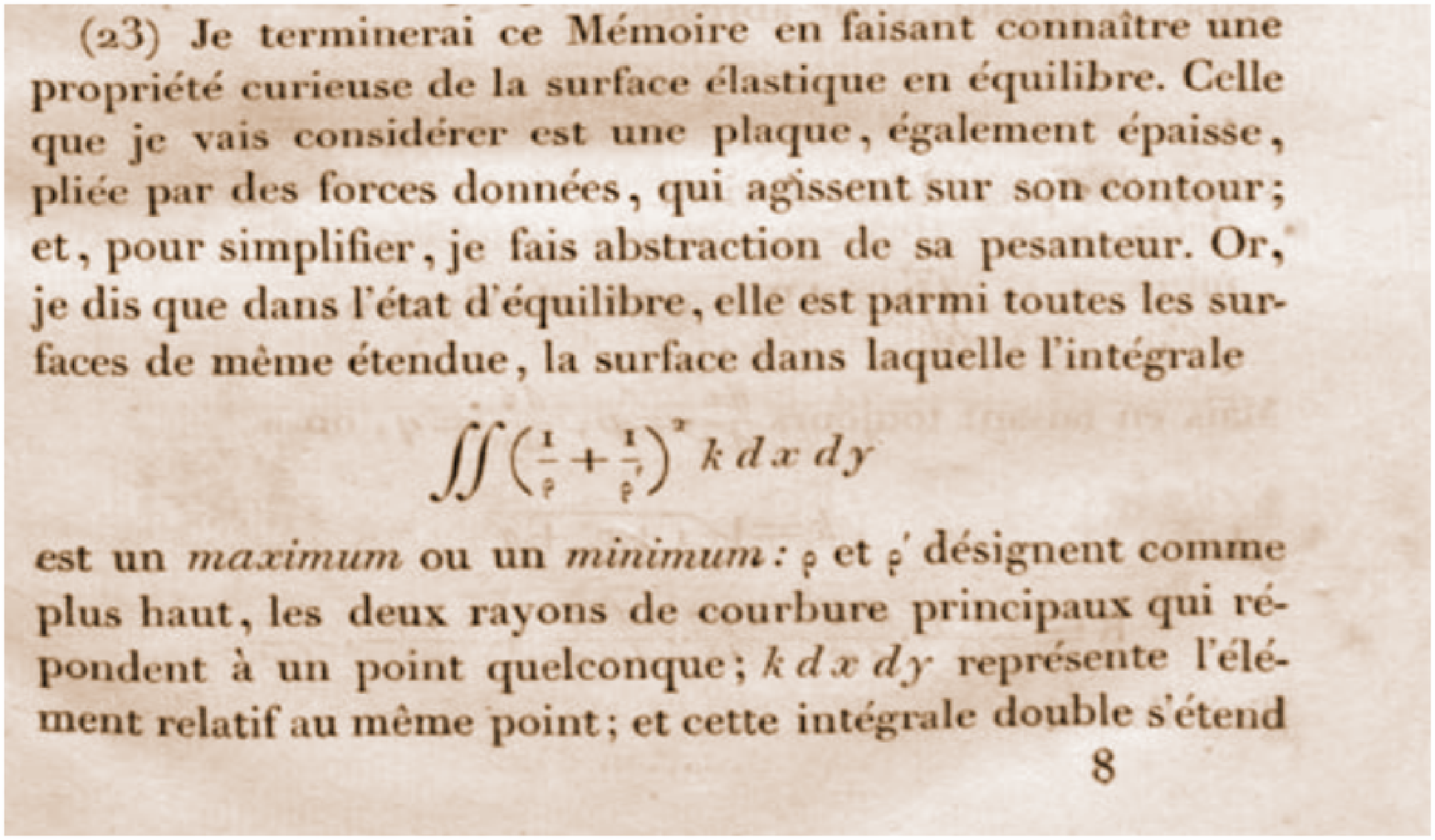}
\caption{An excerpt from Poisson's 1814 (pub. 1816) memoir where he considers
the Lagrangian nowadays known as Willmore functional.
}
\end{figure}

Given a curve $\Gamma$ in $\R^m$, $m\geq 3$, consisting possibly of several connected components,
a unit normal ($m-2$)-vector field $\vec{n}_0$ along $\Gamma$,
and a given positive number $A>0$, 
the \textit{Poisson problem for $\Gamma$, $\vec{n}_0$ and $A$} consists in finding 
an oriented, immersed surface $\Sigma\subset\R^m$ 
which is a minimiser of the energy
\begin{align}
\label{eq:Willmore-energy-Sigma}
W(\Sigma)=\int_\Sigma |\vec{H}_\Sigma|^2\,\dd vol_\Sigma,
\end{align}
among all surfaces having boundary $\Gamma$,
Gauss map along $\Gamma$ equal to $\vec{n}_0$
and area $A$. 
Here, $\vec{H}_\Sigma$ and $\dd vol_\Sigma$ denote, respectively,
the mean curvature vector of $\Sigma$
and the area element of $\Sigma$ induced by the immersion
of $\Sigma$ into $\R^m$. 
The name is after S.D. Poisson's memoir, \cite{PoissonMemoireElastique},
who considered equilibrium states of thin, clamped elastic plates,
and found that critical points of the energy
\eqref{eq:Willmore-energy-Sigma}
(see Fig. \ref{fig:Poisson-excerpts}).

The hypothesis that the elastic energy density be proportional to the mean curvature
is due to S. Germain's seminal research on elastic plates
(later collected in the monograph \cite{GermainReCherchesElastiques}),
who in turn built on earlier one-dimensional models concerning the vibration of elastic beams  investigated
by Euler and Jacques Bernoulli.
Variants of such hypothesis are still in use
nowadays for the study of equilibria of elastic plates,
albeit often in a ``linearised'' version
where \eqref{eq:Willmore-energy-Sigma} is replaced by the simpler biharmonic energy
of a graph.
We refer  for instance to \cite{MR884707, MR1486043, MR2667016} for more on
this subject,
and to \cite{MR619691,MR919065,SzaboHohereMechanik} for a historical perspective on the
development of the subject.
The energy \eqref{eq:Willmore-energy-Sigma} is nowadays known
as \emph{Willmore energy} after T.J. Willmore \cite{MR0202066} and R. Bryant 
\cite{MR772125} who re-introduced it in a modern and geometric perspective.
Finally, we remark that the Poisson problem may be seen 
as a generalisation of the Plateau's problem
(see for instance \cite{MR992402,MR2780140,MR2566897,MR2760441}),
since, for certain special choices of the field $\vec{n}_0$ and of $A$,
there will be minimal surfaces $\Sigma$ (i.e. satisfying $\vec{H}_\Sigma=0$)
bounding $\Gamma$ which are then absolute minimisers for \eqref{eq:Willmore-energy-Sigma}.
\bigskip

When the topology of the surface $\Sigma$ is fixed, the Poisson problem can be reformulated
in a more tractable way that we now describe.
If $g_\Sigma$ and $\vec{\II}_\Sigma$ denote the induced metric and the second fundamental
form of $\Sigma$ and, for every $x\in\Sigma$, $|\vec{\II}_\Sigma(x)|_{g_\Sigma}$
denotes the norm of $\vec{\II}_\Sigma(x)$ on 
$(T_x\Sigma, g_\Sigma(x))$\footnote{
Geometrically, we have $|\vec{\II}_\Phi(x)|^2 = k_1(x)^2+k_2(x)^2$,
where $k_1(x)$ and $k_2(x)$ are the principal curvatures of $\Sigma$ at $x$.},
the \emph{total curvature energy} of $\Sigma$ is defined as
\begin{align}
\label{eq:total-curvature-Sigma}
E(\Sigma)=\int_\Sigma
|\vec{\II}_{\Sigma}|^2_{g_\Sigma}\,\dd vol_\Sigma.
\end{align}
Since there holds
\begin{align}
\label{eq:II-norm-H-K}
|\vec{\II}_{\Sigma}(x)|^2_{g_\Sigma} = 4|\vec{H}_\Sigma(x)|^2-2K_\Sigma(x),
\end{align}
where $K_\Sigma$ denotes the Gauss curvature of $\Sigma$,
we have that
\begin{align*}
4W(\Sigma) - E(\Sigma) = 2\int_{\Sigma}K_\Sigma\,\dd vol_\Sigma.
\end{align*}
By virtue of the Gauss-Bonnet theorem (see for instance \cite[Chapter 6]{MR2964051})
there holds
\begin{align}\label{identityWE}
\int_\Sigma K_\Sigma\,\dd vol_\Sigma
=2\pi\chi(\Sigma) - \int_{\partial\Sigma}k_g,
\end{align}
where $\chi(\Sigma)$ denotes the Euler-Poincar\'e characteristic of $\Sigma$
and $\int_{\partial\Sigma} k_g$ is the integral (or sum of integrals when $\partial\Sigma$ consists of 
multiple connected components) of the geodesic curvature of $\partial\Sigma$
as a positively oriented curve in $\Sigma$.
Hence, if the topology of $\Sigma$ is fixed, the difference between 
the functionals $4W$ and $E$ is a null Lagrangian,
consequently the Willmore and the total curvature energy are equivalent from
the variational point of view for the treatment of the Poisson problem.
The advantage of working with the total curvature energy is that it has a better
coercivity property and that it controls the number of branch points\footnote{
We recall that, roughly speaking, a branch point is a point
so that the immersion $\Sigma\hookrightarrow\R^m$ degenerates.
A classical example is the following: in $\R^4\simeq\C^2$ consider the map
$
\Phi(z)= (z^2,z^3)
$,
for $z\in B_1(0)$. 
Then $\Sigma=\Phi(B_1(0))$ is a surface which is branched at the origin.
Branch points play an important role in the study of minimal and
Willmore surfaces. We refer for instance to the monographs \cite{MR2780140,RivierePCMI}.}
with
their multiplicities (see for instance \cite{MR3276119} or \cite{RivierePCMI}),
so the variational problem will be well-posed with $E$.
\par
We observe that there could be differences between minimizers of the two Lagrangians $W(\Sigma) $ and $ E(\Sigma)$ in the case of the presence of interior branch points, since the identity \eqref{identityWE} does not hold anymore.\par

 We present in this paper a solution to the Poisson problem,
in the   case where $\Gamma$ is a connected, simple, closed curve
and $\Sigma$ is (the image of) a parametrised, possibly branched, immersed disk 
$\Phi:\overline{B_1(0)}\to\R^m$,
(here we denote $B_1(0)=\{z\in\R^2:|z|<1\}$).
Let us define the class of ``admissible'' data $(\Gamma,\vec{n}_0,A)$
for which we can solve the problem.
\begin{definition}
\label{def:admissible-triple}
A triple $(\Gamma, \vec{n}_0,A)$ 
curve $\Gamma\subset\R^m$, a unit-normal $(m-2)$-vector field
$\vec{n}_0$ and a real number $A>0$ is called \emph{admissible for the Poisson problem}
if $\Gamma$ and $\vec{n}_0$ are of class $C^{1,1}$,
$\Gamma$ is simple and closed, and if there is at least one map
$\Phi\in W^{1,2}(B_1(0),\R^m)$ so that
\begin{enumerate}[(i)]
\item  $\Phi$ is conformal, i.e. it satisfies
\begin{align*}
|\partial_1\Phi|
=|\partial_2\Phi|
\quad\text{and}\quad\langle\partial_1\Phi,
\partial_2\Phi\rangle=0
\quad\text{a.e. in } B_1(0),
\end{align*}
and it is a $C^1$-immersion \underline{up to the boundary} away from a finite number of \underline{interior} branch points,
\item the Gauss map $\vec{n}_\Phi$ of $\Phi$   is of class $W^{1,2}$, i.e.
\begin{align*}
\int_{B_1(0)}|\nabla\vec{n}_\Phi|^2\,\dd x<+\infty,
\end{align*}

\item if $\gamma:[0,\mathcal{H}^1(\Gamma)]/\sim\to\R^m$
is a chosen arc-length parametrisation of $\Gamma$, 
there exist a homeomorphism $\sigma_\Phi:S^1\to[0,\mathcal{H}^1(\Gamma)]/\sim$
so that, for every $x\in\partial B_1(0)=S^1$
there holds
\begin{align*}
\Phi(x)=\gamma(\sigma_\Phi(x))\quad\text{and}\quad
\vec{n}_\Phi(x) 
=\vec{n}_0(\gamma(\sigma_\Phi(x))),
\end{align*}
\item there holds 
\begin{align*}
\Area(\Phi)=\frac{1}{2}\int_{B_1(0)}|\nabla\Phi|^2\,\dd x=A.
\end{align*}
 \end{enumerate}
\end{definition}
Note that, for a given triple $(\Gamma, \vec{n}_0, A)$ it is not so obvious
to determine directly whether is it admissible of not.
However, an elementary application of the $h-principle$
(see \cite{MR1909245,MR864505}) allows us, for any given $\Gamma$ and
$\vec{n}_0$ as in definition \ref{def:admissible-triple},
to prove the existence of some $A_0>0$ such that, for every $A\geq A_0$
the triple $(\Gamma,\vec{n}_0,A)$ is admissible.
Such a construction is presented in appendix \ref{app:contruction-competitors}
(lemma \ref{lemma:existence-competitors}).
In particular when $m=3$, if one requires the map $\Phi$ as in the definition
\ref{def:admissible-triple} not to have any branch points, $(\Gamma,\vec{n}_0)$
need to satisfy a topological constraint, namely, if $\mathbf{t}$ denotes the
tangent vector of $\Gamma$, the map
\begin{align*}
x\mapsto
(\mathbf{t}\times\vec{n}_0,\mathbf{t},\vec{n}_0)(x),\quad x\in S^1,
\end{align*}
has to define a non-nullhomotopic loop in the space of 
special orthogonal matrices $SO(3)$.
\bigskip

The first main result of this paper is the following.
\begin{theorem}
\label{thm:existence-solution-Poisson-problem}
Let $(\Gamma,\vec{n}_0,A)$ be an admissible triple for the Poisson problem.
Then, there exists a map $\Phi_\infty:B_1(0)\to\R^m$ verifying  this data which
is a $C^1(B_1(0)\setminus\{a_{1},\ldots,a_{\ell}\})$ conformal immersion, where $\{a_{1},\ldots,a_{\ell}\}$ is the finite set (possibly empty) of  branched points in $B_1(0)$,
minimising the total curvature energy $E$ in this  class.
\end{theorem}

The second main result is the following.
\begin{theorem}
\label{thm:regularity-solution-Poisson-problem}
Let $(\Gamma,\vec{n}_0,A)$ be an admissible triple for the Poisson problem.
Every  minimising map   $\Phi_\infty$ as in theorem \ref{thm:existence-solution-Poisson-problem}
satisfies the Euler-Lagrange equation
\begin{align}
\label{eq:Euler-Lagrange-Poisson-problem}
\divop\left(
\nabla\vec{H}_{\Phi_\infty} -3\pi_{\vec{n}_{\Phi_\infty}}(\nabla\vec{H}_{\Phi_\infty})
+\star(\nabla^\perp\vec{n}_{\Phi_\infty}\wedge\vec{H}_{\Phi_\infty})
+ c\nabla \Phi_{\infty}
\right)
=0
\quad\text{in }\mathcal{D}'(B_1(0),\R^m),
\end{align}
where $\pi_{n_{\Phi_\infty}}$ is the orthogonal projection
onto the normal bundle of $\Phi_\infty$, 
$\nabla^\perp\vec n_{\Phi_\infty}
= (-\partial_2\vec{n}_{\Phi_\infty},\partial_1\vec{n}_{\Phi_\infty})$
and $c\in\R$.
Moreover, there exists $0<\alpha<1$ so that $\Phi_\infty$ is of class
$C^{1,\alpha}$ up to the boundary and $\vec{n}_{\Phi_\infty}$
extends to a map of class $C^{0,\alpha}$ up to the boundary.
\end{theorem}
We recall that any map $\Phi_\infty$ as in theorem \ref{thm:regularity-solution-Poisson-problem}
is smooth in $B_1(0)$ away from its branch points
(this follows from an adaptation of the results in \cite{MR2430975},
similarly as done in \cite{MR2989995, MR3276119})
and, since the equation \eqref{eq:Euler-Lagrange-Poisson-problem} is satisfied
through the branch points, from the study of the singularities of Willmore
surfaces (see \cite{MR2430975, MR3096502} and \cite{MR2119722,MR2318282})
it follows that, for every $0<\beta<1$,
$\Phi_\infty$ is of class $C^{2,\beta}$ at the branch points
and its Gauss map $\vec{n}_{\Phi_\infty}$ extends to a map
of class $C^{1,\beta}$  at the branch points.
\begin{remark}
\label{rmk:expected-improvements}
Notice that the assumptions $\Gamma,\vec{n_0}\in W^{2,\infty}$ are far from being optimal. We expect that the conclusions of theorems \ref{thm:existence-solution-Poisson-problem}
\ref{thm:regularity-solution-Poisson-problem} to hold under weaker assumptions on the admissible boundary data.
 \end{remark}

Let us put our results in a broader context.
Nitsche \cite{MR1218374} discussed various boundary conditions
for the Willmore and related type of functionals, and proved existence
and uniqueness results for a class of such problems,
also considering a volume constraint, when the surfaces are graphs in $\R^3$ and the boundary data
are sufficiently small in $C^{4,\alpha}$-norm.
Recently Deckelnick-Grunau-R\"ogers \cite{MR3665920} also consider
the minimisation over graphs in $\R^3$ of the Willmore functional
(also plus a constant times integral of the Gauss curvature) subject to
various boundary conditions and deduced compactness results in the
$L^1$-topology, and from this, also a lower-semicontinuity for a suitably
defined relaxation of the Willmore functional.
A considerable series of results
(we refer to \cite{
MR3606563,MR3665920,MR3116012,MR3010281,MR2770424,MR2729304,MR2568881,MR2480063}
and the references therein) is available when considering boundary value problems
for the Willmore functional under the hypothesis that the surfaces in consideration
are surfaces of revolution around an axis in $\R^3$ (hence the boundary consist of two circles).
Sch\"atzle \cite{MR2592972},
by working on the sphere $S^m\subset\R^{m+1}$, has proved the existence,
for arbitrary smooth boundary data $\Gamma$ and $\vec{n}_0$ and without area constraint,
of a branched immersion, smooth away from the finitely many branch points,
satisfying the strong form of the Willmore equation in $S^m$, namely
\begin{align}\label{ELeq}
\Delta^{S^m}_{g}\vec{H}_{\Phi,S^m}
+Q(\mathring{\II}_{\Phi,S^m})\vec{H}_{\Phi,S^m}=0,
\end{align}
away from the branch points
(we refer to the paper in question or to
the monograph by Kuwert and the same author in \cite{MR2906051}
for more detail on this equation).
However this equation has no meaning at the branch points and
it is not proved that, starting from arbitrary data $\Gamma,\vec{n}_0$
in $\R^m$, projecting them stereographically into $S^m$,
and then considering the projected-back Willmore surface obtained in $S^m$
(recall that the Willmore equation is conformally invariant), one gets
a surface which  is ``Poisson-minimal'', i.e. that it solves the Poisson problem
or also whether it passes through $\infty$ or not.\par

\begin{figure}
\centering
\includegraphics[scale=0.5]{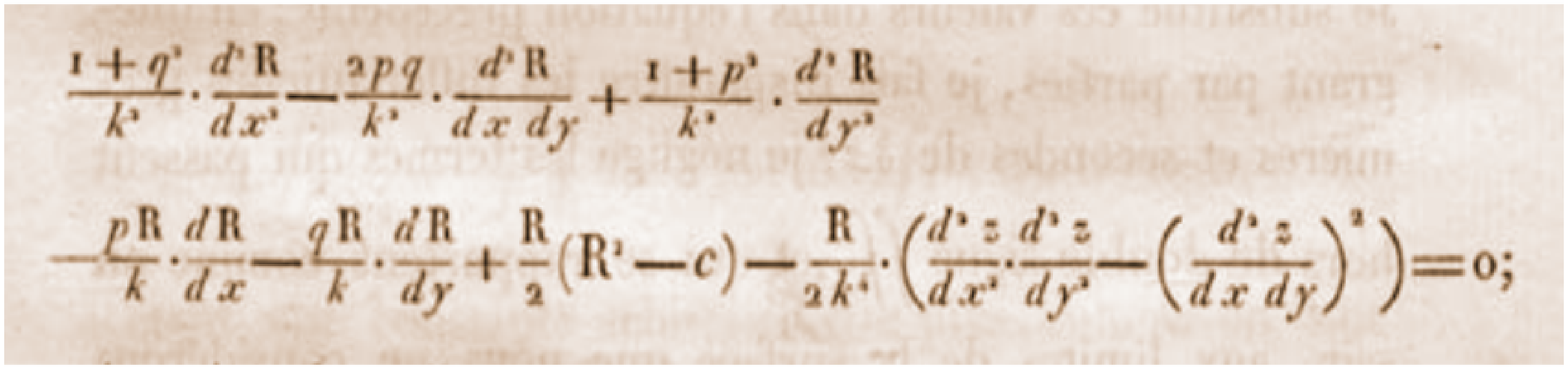}
\caption{Euler-Lagrange Equation for Willmore Energy derived by Poisson.}
\label{fig:Poisson-excerptsbis}
\end{figure}
Concerning the equation \eqref{ELeq}, we mention that it already appeared   in Poisson's  memoirs for special cases of graphs, (see Fig. \ref{fig:Poisson-excerptsbis}). This is actually 
striking since at the time Lagrange's theory was not so developed.\par
Alexakis-Mazzeo \cite{MR3415767} consider smooth, properly embedded and complete Willmore
surfaces in the hyperbolic space $\mathbb{H}^3$ and relate the regularity
of their asymptotic boundary with the smallness of a suitable version of the Willmore energy.
In a recent paper Alessandroni-Kuwert
\cite{MR3463464}, considering a free-boundary problem for the Willmore functional,  have proved the existence (and non-uniqueness) of smooth Willmore disk-type surfaces in $\R^3$
with prescribed but small value of the area whose boundary lays on
the boundary of a smooth, bounded domain.
\bigskip

At this point we would like to mention some  interesting questions and open problems related to the 
the Poisson problem, some of which are the aim of future investigation.
Beside our expectations regarding the  boundary regularity
(see remark \ref{rmk:expected-improvements}) one may for instance consider:
\begin{itemize}
\item existence and properties of non-minimising i.e. saddle-type Willmore surfaces
having prescribed boundary data and area
(a partial answer is in \cite{MR2592972} described above);
\item the solution to the Poisson problem in the case of a boundary curve consisting
of multiple connected components
(in the spirit of its Plateau-counterpart as done in the classical work by Meeks-Yau
\cite{MR670745}) or among manifold with a non-trivial topology
(an important tool for this sould be the work by Bauer-Kuwert \cite{MR1941840}),
or both;
\item an investigation on a version of the Poisson problem where
no area bound is prescribed;
\item free-boundary versions of the Poisson problem in the case of surfaces with
arbitrary prescribed area and boundary laying in some submanifold of
$\R^m$ such as the sphere (as done in the forementioned work \cite{MR3463464} but
with no restriction on the prescribed value of the area).
\end{itemize}
\bigskip

Let us discuss some central issues involved in our proofs.
Our approach to the Poisson problem
is in the framework of the so-called weak immersions with
$L^2$-bounded second fundamental form, developed by the third author
in \cite{MR2430975,MR3276154} (see also the monographs \cite{RiviereConformal,RivierePCMI}).
As it has been discovered in \cite{MR2430975}, a central role
in the analysis of Willmore surfaces  is played by peculiar estimates that are valid
for PDEs involving Jacobian  non-linearities.
In a nutshell, one considers the following Dirichlet problem:
\begin{align}
\label{eq:Jacobian-Dirichlet}
\left\{
\begin{aligned}
-\Delta u &=  \partial_1a\,\partial_2 b - \partial_2 a\,\partial_1 b
&& \text{in}\, B_1(0),\\
u&=0 &&\text{on}\, \partial B_1(0),
\end{aligned}\right.
\end{align}
where $a,b\in W^{1,2}(B_1(0))$.
It is well-known since Wente's work (see\cite{MR0336563, MR592104} and \cite{MR733715})
that, for every $a,b$,   the solution  $u$ to
\eqref{eq:Jacobian-Dirichlet} is $ C^0(\overline{B_1(0)})\cap W^{1,2}(B_1(0))$ and satisfies  
\begin{align}
\label{eq:Wente-estimateintr}
\|u\|_{L^\infty(B_1(0))}+\|\nabla u\|_{L^2(B_1(0))}
\leq C\|\nabla a\|_{L^2(B_1(0))}\|\nabla b\|_{L^2(B_1(0))}.
\end{align}
for some constant $C>0$ independent of $u,a$ and $b$ 
(see appendix \ref{app:integrability-compensation}).
The estimate \eqref{eq:Wente-estimateintr} is a priori not trivial since, being
on a first glance $\Delta u$ only in $L^1(B_1(0))$, the standard Calder\'on-Zygmund
estimates do not apply.
Such sharp estimates allow for a successful treatment of many issues regarding both
the compactness and the regularity of Willmore surfaces (we refer to \cite{RivierePCMI}).
In the case of the Poisson problem, however, it turns out that 
the central role is played by a corresponding Neumann problem, namely
\begin{align}
\label{eq:Jacobian-Neumann}
\left\{
\begin{aligned}
-\Delta u &=  \partial_1a\,\partial_2 b - \partial_2 a\,\partial_1 b
&& \text{in}\, B_1(0),\\
\partial_\nu u&=\partial_{\tau} a\, b &&\text{on}\, \partial B_1(0),
\end{aligned}\right.
\end{align}
for $a,b\in W^{1,2}(B_1(0))$,
where $\partial_{\tau}$ denotes the tangential derivative along $\partial B_1(0)$.
The first and second author has recently shown in \cite{MR3764916}
(see also \cite{HirschWente}) that, although any solution of \eqref{eq:Jacobian-Neumann}
is naturally in $W^{1,2}(B_1(0))$ (we refer again to the discussion in appendix \ref{app:integrability-compensation}),
there are explicit counter-examples to the $L^\infty$-estimate
if no further assumptions are made on $a$ and $b$
(one could even find counter-examples to the $W^{1,2}$-estimate if
a homogeneous boundary condition is chosen, see \cite{MR3764916}).
We will prove that
through suitable a-priori estimates on the general Neumann problem
similar to those appearing in \cite{MR1132783, MR3399138} and 
a careful analysis on the boundary behaviour of the solution
of \eqref{eq:Jacobian-Neumann}, a sufficient condition to achieve such crucial estimate
will be on the trace of $a,b$ to be in $W^{1,p}(\partial B_1(0))$, for 
some $p>1$.
\par
The regularity {\em up to the boundary} of a minimising  map of \eqref{eq:total-curvature-Sigma} in theorem \ref{thm:regularity-solution-Poisson-problem} follows from a boundary $\varepsilon$-regularity Lemma which is in turn obtained by suitable comparison arguments.
It is expected such an $\varepsilon$-regularity Lemma up to the boundary to hold for general critical points of the  Energy \eqref{eq:total-curvature-Sigma} (an interior $\varepsilon$-regularity for 
such critical points has been already  proved by the third author in \cite{MR2430975}).\par
 We conclude by introducing  some basic notation and geometric notions that will be used throughout the paper.

\subsection{Basic notation}
Beside the common notation,
in this paper we denote:
\begin{align*}
\R^2_+&=\{(x^1,x^2):x^2>0\}
&&\text{upper half-space},\\
B_r^+(0)&=B_r(0)\cap \R^2_+\quad\text{with } r>0
&&\text{upper half-ball of radius }r,\\
\vec{\varepsilon}_1&,\ldots,\vec{\varepsilon}_m
&&\text{canonical basis of }\R^m.
\end{align*}
We will write
$
\partial B_r^+(0)= rI + rS,
$
where
\begin{align*}
rI &= \{(x^1,0):-r<x^1<r\}\simeq(-r,r)
&&\text{base diameter},\\
rS &= \{(x^1,\sqrt{1-(x^1)^2}):-r<x^1<r\}
&&\text{upper semi-circle}.\\
\end{align*}
When $\Omega$ is a domain of $\R^2$, and
$f=(f^1,\ldots,f^n):\Omega\subseteq{\R^2}\to\R^n$,
is a given function, we adopt the following notation:
$$\nabla f = 
\begin{pmatrix}
\partial_1 f\\
\partial_2 f
\end{pmatrix},\qquad~~~\nabla^\perp f = 
\begin{pmatrix}
-\partial_2 f\\
\partial_1 f
\end{pmatrix},$$
 and if $g:\Omega\subseteq{\R^2}\to\R^n$ is another function, we set
\begin{align*}
\langle \nabla^\perp f,g\rangle 
&=\begin{pmatrix}
-\langle\partial_2 f,g\rangle\\
\langle\partial_1 f,g\rangle
\end{pmatrix},\\
\langle\nabla^\perp f,\nabla g\rangle
&=\divop(\langle\nabla^\perp f,g\rangle)
=\sum_{i=1}^n\langle\nabla^\perp f^i,\nabla g^i\rangle
=\sum_{i=1}^n\partial_1 f^i\partial_2 g^i - \partial_2 f^i\partial_1 g^i.
\end{align*}

\subsection{Geometric notions}
We list here the geometric objects that will appear in the paper.
For more information, we refer the reader 
to textbooks of differential geometry of curves and surfaces 
such as \cite{MR2964051}, \cite{MR2572292} and \cite{MR615912}
and to the monograph \cite{RivierePCMI}.
For a given immersion $\Phi:\overline{B_1(0)}\to(\R^m,\langle\;,\;\rangle)$
into the Euclidean $m$-dimensional space, we denote:
\begin{enumerate}[-]
\item the induced metric on $B_1(0)$:
 $g_\Phi=(g_\Phi)_{ij} = (\langle\partial_i\Phi,\partial_j\Phi\rangle)_{ij}$,
 with inverse $(g_\Phi^{ij})_{ij}$,
\item the induced area element: $\dd vol_\Phi=\sqrt{\det g_\Phi}\,\dd x$,
\item the Gauss map:
$\vec{n}_\Phi:B_1(0)\to\Gr_{m-2}(\R^m)$  defined as
\begin{align}
\label{eq:Gauss-map}
\vec{n}_\Phi(x)=
\star\frac{\partial_1\Phi(x)\wedge\partial_2\Phi(x)}
{|\partial_1\Phi(x)\wedge\partial_2\Phi(x)|},
\end{align}
where $\star$ is the Hodge operator in $\R^m$
(here $\Gr_{m-2}(\R^m)$ denotes the $(m-2)$-dimensional Grassmannian of $\R^m$
and  $\star$ is the Hodge operator in $\R^m$ see e.g. \cite[Chapter 2]{MR615912})\footnote{
when $m=3$, we have the canonical identification
$\star(V\wedge W)=V\times W$, where $\times$ is the vector product,
so in this case $\vec{n}_\Phi$ can be considered as a $S^2$-valued map.}.
It can be represented as 
$\vec{n}_\Phi = \vec{n}_{\Phi,1}\wedge\ldots\wedge
\vec{n}_{\Phi,m-2}$,
where the $\vec{n}_{\Phi,i}$ are normal vector fields for $\Phi$ so
that the $m$-ple
$
(\partial_1\Phi,\partial_2\Phi,\vec{n}_{\Phi,1},\ldots,\vec{n}_{\Phi,m-2})
$
defines a positively oriented basis of $\R^m$,
\item the Hessian: $\nabla^2\Phi$,
and its norm:
$
|\nabla^2\Phi|_{g_\Phi}^2=\sum_{i,j,k,l}
g_{\Phi}^{ik}g_{\Phi}^{jl}\langle\partial^2_{ij}\Phi,\partial^2_{kl}\Phi\rangle,
$
\item the second fundamental form: $\vec{\II}_\Phi = \pi_{\vec{n}_\Phi}(\nabla^2\Phi)$,
where $\pi_{\vec{n}_\Phi}$ denotes the orthogonal projection onto the normal bundle of $\Phi$,
and its norm:
$
|\vec{\II}_\Phi|_{g_\Phi}^2=\sum_{i,j,k,l}
g_{\Phi}^{ik}g_{\Phi}^{jl}\langle\vec\II_{ij},\vec\II_{kl}\rangle,
$
\item the Gauss curvature: $K_\Phi$ and the mean curvature vector: $\vec{H}_\Phi$,
\item  the area functional:
$\Area(\Phi) = \int_{B_1(0)}\,\dd vol_\Phi$,
\item the total curvature functional:
$E(\Phi)=\int_{B_1(0)}|\vec{\II}_{\Phi}|^2_{g_\Phi}\,\dd vol_\Phi$.
\end{enumerate}
When $\Phi$ is conformal with conformal factor $\e^{2\lambda_\Phi}=|\nabla\Phi|/\sqrt{2}$,
some of the above mentioned objects can be written in an easier way in that
\begin{align*}
&(g_\Phi)_{ij} = \e^{2\lambda_\Phi}\delta_{ij},\\
&(g_\Phi)^{ij} = \e^{-2\lambda_\Phi}\delta^{ij},\\
&\dd vol_\Phi = \e^{2\lambda_\Phi}\,\dd x,\\
&|\nabla^2\Phi|^2_{g_\Phi} = \e^{-4\lambda_\Phi}|\nabla^2\Phi|^2,\\
&|\vec\II_\Phi|^2_{g_\Phi} = \e^{-4\lambda_\Phi}|\vec\II_\Phi|^2.
\end{align*}
A direct computation reveals that the general identity
$
|\nabla^2\Phi|_{g_\Phi}
= |\pi_{T_\Phi}(\nabla^2\Phi)|_{g_\Phi}
+|\vec{\II}_\Phi|^2_{g_\Phi},
$
where $\pi_{T_\Phi}$ denotes the orthogonal projection onto the tangent bundle of $\Phi$,
can be rewritten, when $\Phi$ is conformal, as
\begin{align}
\label{eq:pointwise-Hessian-identity}
|\nabla^2\Phi|^2_{g_\Phi} 
&= 4\e^{-2\lambda_\Phi}|\nabla\lambda_\Phi|^2 +|\vec{\II}_\Phi|^2_{g_\Phi}.
\end{align}
If $\Gamma\subset\R^m$ is a simple, closed curve
with a chosen arc-length parametrization $\gamma:[0,\hau^1(\Gamma)]/\sim\to\Gamma$
 and $\vec{n}_0=\vec{n}_0(\gamma(\cdot))$ is some unit-normal $(m-2)$-vector field
along $\Gamma$,
the geodesic curvature $k_g$ of $\Gamma$ (with respect to $\vec{n}_0)$
is defined as follows:
if $\mathbf{t}=\dot{\gamma}$
denotes the unit-tangent vector of $\Gamma$ we set
(see \cite[Chapter 5]{MR2964051})
\begin{align}
\label{eq:geodesic-curvature}
k_g=\left\langle
\dot{\mathbf{t}}
,\star(\vec{n}_0\wedge\mathbf{t})\right\rangle=
\left\langle
\ddot{\gamma}
,\star(\vec{n}_0\wedge\dot{\gamma})\right\rangle.
\end{align}

\begin{definition}[Weak notions of immersions]
\label{def:weak-immersions}
~
\begin{enumerate}[(i)]
\item 
$\mathcal{F}_{B_1(0)}$ is the set of
 Lipschitz maps $\Phi:B_1(0)\to\R^m$ 
so that there exists a bi-Lipschitz homeomorphism
$\psi:B_1(0)\to B_1(0)$ so that $\Psi:=\Phi\circ \psi$ is conformal and
there exists a finite (possibly empty) set $\{a_1,\ldots,a_N\}\subset \overline{B_1(0)}$
 so that, for every compact set $K\subset \overline{B_1(0)}\setminus \{a_1,\ldots,a_n\}$ there
is a constant $C_K>0$ so that there holds
\begin{align}
 \sqrt{\det{g_\Psi}}=|\nabla\Psi|^2/2\geq C_K
\quad\text{a.e. in } K,
\end{align}
and if the Gauss map of $\Psi$, defined a.e. as in \eqref{eq:Gauss-map}
belongs to $W^{1,2}(B_1(0),\Gr_{m-2}(\R^m))$.
\item
If $\Gamma\subset\R^m$ is a simple, closed curve and
$\vec{n}_0$ is a unit-normal $(m-2)$-vector field 
along $\Gamma$,
$\mathcal{F}_{B_1(0)}(\Gamma,\vec{n}_0)$
is the set of Lipschitz maps $\Phi:B_1(0)\to\R^m$ so that
$\Phi$ is an element of $\mathcal{F}_{B_1(0)}$
and for some homeomorphism 
$\sigma_\Phi:S^1\to[0,\hau^1(\Gamma)]/\sim$ there holds
\begin{align}
\Phi(x)=\gamma(\sigma_\Phi(x))\quad\text{and}\quad
\vec{n}_\Phi(x) 
=\vec{n}_0(\sigma_\Phi(x)),
\quad \text{for }x\text{ in }\partial B_1(0)=S^1,
\end{align}
where $\gamma:[0,\hau^1(\Gamma)]/\sim\to\Gamma$
is a fixed arc-length parametrization of $\Gamma$.
\item For $A>0$, $\mathcal{F}_{B_1(0)}(\Gamma,\vec{n}_0,A)$
is the set of Lipschitz maps $\Phi:B_1(0)\to\R^m$
so that $\Phi\in\mathcal{F}_{B_1(0)}(\Gamma,\vec{n}_0)$
and 
$
\Area(\Phi)=A.
$
\end{enumerate}
\end{definition}
For $\Phi$ as in definition \ref{def:weak-immersions},
we may write:
\begin{align*}
E(\Phi)
=\int_{B_1(0)}|\vec{\II}_\Phi|^2_{g_\Phi}\,\dd vol_\Phi
=\int_{B_1(0)}|\nabla\vec{n}_{\Phi\circ\psi}|^2_{g_\Phi}\,\dd x
<+\infty,
\end{align*}
where $\psi:B_1(0)\to B_1(0)$ is the homeomorphism giving the conformal
parametrization of $\Phi$.
\smallskip

If $\Phi:B_1(0)\to\R^m$ is a possibly branched conformal immersion,
 the logarithm of its conformal factor $\lambda_\Phi=\log(|\nabla\Phi|/\sqrt{2})$
will be a weak solution of the so-called Liouville equation:
\begin{align}
\label{eq:Liouville}
\left\{
\begin{aligned}
-\Delta \lambda_\Phi &=  K_\Phi\e^{2\lambda_\Phi}-2\pi\sum_{i=1}^{N}n_{i}\delta_{a_{i}} 
&& \text{in}\, B_1(0),\\
\partial_{\nu}\lambda_{\Phi}&=k_g(\sigma_\Phi)\e^{\lambda_\Phi}-1 &&\text{on}\, \partial B_1(0),
\end{aligned}\right.
\end{align}
where
$a_1,\ldots, a_n$ are the branch points of $\Phi$ and $n_i\in \N$ are their respective multiplicities.

The relative abundance of maps satisfying
definition \ref{def:weak-immersions} is 
implied by the following result which makes use
 M\"uller-\v{S}ver\'ak theory of weak isothermic 
charts \cite{MR1366547} and H\'elein's moving frame technique
(for more deails we refer to \cite{MR1913803} and \cite{RivierePCMI}):
\begin{lemma}
\label{lemma:conformal-coordinates}
Let $\Phi\in W^{1,\infty}(B_1(0),\R^m)$   satisfy the following conditions:
\begin{enumerate}
\item[$i)$]
there exists some $C>0$ such that
\begin{align*}
|\nabla\Phi|^2\leq C\sqrt{\det g_\Phi}
\quad\text{a.e. in } B_1(0),
\end{align*}
\item[$ii)$] there exists a finite (possibly empty) set $\{a_1,\ldots, a_N\}\subset B_1(0)$
such that  for every compact set $K\subset \overline{B_1(0)}\setminus \{a_1,\ldots,a_n\}$  there holds
\begin{align*}
\|\nabla\Phi\|_{L^{\infty}(K)}\geq C_K,
\end{align*}
for some constant $C_K>0$,
\item[$iii)$] 
its  Gauss map, defined a.e. as in \eqref{eq:Gauss-map},
satisfies 
\begin{align*}
\int_{B_1(0)}|\nabla\vec{n}_\Phi|^2_{g_\Phi}\,\dd vol_\Phi
=\int_{B_1(0)}\sum_{\substack{i,j=1,2\\l=1,\ldots,m-2}}
g_{\Phi}^{ij}\langle\partial_i \vec{n}_{\Phi,l},\partial_j\vec{n}_{\Phi,l}\rangle
\dd vol_\Phi<+\infty.
\end{align*}
\end{enumerate}
Then there is a bi-Lipschitz diffeomorphism $\psi:B_1(0)\to B_1(0)$
so that $\Phi\circ\psi$ is conformal and $\Phi\circ\psi\in W^{2,2}(B_1(0),\R^m)$.
\end{lemma}

Some further information about moving frames is 
in appendix \ref{app:frames}.

\par
\bigskip
The paper is organised as follows:
\begin{enumerate}[-]
 \item
in section \ref{sec:exponential-integrability} we prove a priori estimate
on boundary exponential intergrability for the Neumann problem that will 
be needed in the sequel;
\item in section \ref{sec:exitence-minimiser}
we prove   the  first part of theorem \ref{thm:existence-solution-Poisson-problem}, i.e.  the existence of a  Lipschitz minimising immersion,
 
\item in section \ref{sec:regularity} we prove an $\varepsilon$-regularity result up to the boundary for minimisers  by constructing suitable
competitors that permits us to   conclude the proof of
theorems  \ref{thm:existence-solution-Poisson-problem}-\ref{thm:regularity-solution-Poisson-problem}.
\end{enumerate}

\section{An Estimate for the Neumann Problem}
\label{sec:exponential-integrability}
For given $f\in L^1(B_1(0))$, $g\in L^1(\partial B_1(0))$
we recall that a function $u\in W^{1,1}(B_1(0))$
is said to weakly solve the Neumann problem for the Poisson equation in $B_1(0)$:
\begin{align}
\label{eq:Neumann}
\left\{
\begin{aligned}
-\Delta u &=f &&\text{in }B_1(0),\\
\partial_\nu u & =g &&\text{on }\partial B_1(0),
\end{aligned}
\right.
\end{align}
if, for every $\psi\in C^{\infty}(\overline{B_1(0)})$,
there holds
\begin{align*}
\int_{B_1(0)} f\,\psi\,\dd x + \int_{\partial B_1(0)}g\,\psi\,\dd\hau^1
=\int_{B_1(0)}\langle\nabla u,\nabla\psi\rangle\,\dd x.
\end{align*}
From this expression, it is immediate to see that a necessary condition for 
the existence of a weak solution is that $\int_{B_1(0)} f = -\int_{\partial B_1(0)}g$
(see e.g. \cite{MR1282720} for more on weak formulations of Neumann problems).
Such condition is also sufficient and we have the following representation formula:
\begin{align*}
 u(x)-\dashint_{\partial B_1(0)} u\,\dd\hau^1
=\int_{B_1(0)}\mathcal{G}(x,y)\,f(y)\,\dd y
+\int_{\partial B_1(0)}\mathcal{G}(x,y) g(y)\,\dd\hau^1(y)
\quad x\in B_1(0),
\end{align*}
where $\mathcal{G}$ is the
Green function for the Neumann problem (with zero average on $\partial B_1(0)$),
that is the function:
\begin{align}
\label{eq:Green-ft-Neu}
\mathcal{G}(x,y)=
-\frac{1}{2\pi}
(\log|x-y|+\log|\tilde{x}-y||x|)
+\frac{|y|^2}{4\pi}-\frac{1}{4\pi},
\end{align}
which satisfies, for every $x\in B_1(0)$,
\begin{align}
\label{eq:Green-for-Neumann}
\left\{
\begin{aligned}
-\Delta_y\mathcal{G}(x,\cdot) &=\delta_x-\frac{1}{|D|}, &&\text{in }B_1(0),\\
\partial_{\nu_y}\mathcal{G}(x,\cdot) & =0 &&\text{on }\partial B_1(0),\\
\int_{\partial B_1(0)}\mathcal{G}(x,y)\,\dd\hau^1(y)&=0.
\end{aligned}
\right.
\end{align}
Two such solutions to \eqref{eq:Neumann} differ by a constant.
\bigskip

In what follows, we denote $\frac{1}{2\pi}\int_{\partial B_1(0)}\phi$ by $\overline{\phi}$.
The results we now present are in the spirit of Brezis-Merle \cite[Theorem 1]{MR1132783}
and Da Lio-Martinazzi-Rivi\`ere \cite[Theorem 3.2]{MR3399138}.

\begin{theorem}
\label{thm:global-exp-int}
Let $f\in L^1(B_1(0))$, $g\in L^1(\partial B_1(0))$ satisfy 
$\int_{B_1(0)} f = -\int_{\partial B_1(0)}g$ and let $ u\in W^{1,1}(B_1(0))$ be a weak solution
to the problem
\eqref{eq:Neumann}.
Then, for every $\varepsilon>0$ verifying
$
\|g\|_{L^1(\partial B_1(0))}+2\|f\|_{L^1(B_1(0))}<\pi -\varepsilon
$
 we have 
$$\|\e^{ u-\overline{ u}}\|_{L^p(\partial B_1(0))}\leq C_\varepsilon,$$
for some $C_\varepsilon>0$ depending on $\varepsilon$ and
\begin{align*}
p=\frac{\pi-\varepsilon}{\|g\|_{L^1(\partial B_1(0))}+2\|f\|_{L^1(B_1(0))}}.
\end{align*}
\end{theorem}

\begin{bfproof}[Proof of Theorem \ref{thm:global-exp-int}]
Set $k=\frac{1}{2\pi}\int_{B_1(0)} f = -\frac{1}{2\pi}\int_{\partial B_1(0)} g$.
We write the solution of \eqref{eq:Neumann} as:
\begin{align}
\label{eq:u-bar-u}
u(x) - \overline{ u}
= u_1(x)+ u_2(x) +\frac{k}{2}(1-|x|^2),\quad x\in B_1(0),
\end{align}
where:
\begin{align*}
\left\{
\begin{aligned}
-\Delta u_1 &=f-2k, &&\text{in }B_1(0),\\
\partial_\nu u_1 & =0 &&\text{on }\partial B_1(0),\\
\overline{ u_1}&=0,
\end{aligned}
\right.
\quad\text{and}\quad
\left\{
\begin{aligned}
-\Delta u_2 &=0, &&\text{in }B_1(0),\\
\partial_\nu u_2 & =g+k &&\text{on }\partial B_1(0),\\
\overline{ u_2}&=0.
\end{aligned}
\right.
\end{align*}
\emph{Step 1: study of $ u_1$.}
We set $F=f-2k=f-\dashint_{B_1(0)} f$;
note that there holds
$
\int_{B_1(0)} F=0$
and
$\|F\|_{L^1(B_1(0))}\leq 2\|f\|_{L^1(B_1(0))}$.
The Green function for the Neumann problem
can be written as
\begin{align*}
\mathcal{G}(x,y)=
\frac{1}{2\pi}
\left(
\log\left(\frac{2}{|x-y|}\right)+\log\left(\frac{2}{|\tilde{x}-y||x|}\right)
\right)
+\frac{|y|^2}{4\pi}-\frac{1}{4\pi}-\frac{1}{\pi}\log 2,
\end{align*}
where,
since $|x-y|\leq 2$ and $|\tilde{x}-y||x|\leq 2$,
the term in brackets is non negative.
Then we may write $ u_1$ as:
\begin{align*}
 u_1(x)
=\frac{1}{2\pi}\int_{B_1(0)}
\left(
\log\left(\frac{2}{|x-y|}\right)+\log\left(\frac{2}{|\tilde{x}-y||x|}\right)
\right)
F(y)\,\dd y
+\int_{B_1(0)}\frac{|y|^2}{4\pi}F(y)\,\dd y,
\end{align*}
and in particular, for $x\in\partial B_1(0)$, we have the formula:
\begin{align*}
 u_1(x)=
\frac{1}{\pi}\int_{B_1(0)}\log\left(\frac{2}{|x-y|}\right)F(y)\,\dd y
+\int_{B_1(0)}\frac{|y|^2}{4\pi}F(y)\,\dd y,
\quad x\in\partial B_1(0).
\end{align*}
For $\gamma>0$ then there holds:
\begin{align*}
\frac{\gamma| u_1(x)|}{\|F\|_{L^1(B_1(0))}}
\leq\frac{\gamma}{\pi}\int_{B_1(0)}
\log\left(\frac{2}{|x-y|}\right)
\frac{|F(y)|}{\|F\|_{L^1(B_1(0))}}\,\dd y + \frac{\gamma}{4\pi},
\end{align*}
so by Jensen's inequality (see e.g. \cite[Appendix B]{MR2597943}):
\begin{align*}
\exp\left(\frac{\gamma| u_1(x)|}{\|F\|_{L^1(B_1(0))}}\right)
\leq\int_{B_1(0)}
\left(\frac{2}{|x-y|}\right)^{\frac{\gamma}{\pi}}
\frac{|F(y)|}{\|F\|_{L^1(B_1(0))}}\,\dd y 
\cdot \e^{\frac{\gamma}{4\pi}}.
\end{align*}
Integrating on $\partial B_1(0)$ and using Tonelli's theorem
yields:
\begin{align*}
\int_{\partial B_1(0)}
\exp\left(\frac{\gamma| u_1(x)|}{\|F\|_{L^1(B_1(0))}}\right)\,\dd \hau^1(x)
\leq \int_{B_1(0)}\left\{
\int_{\partial B_1(0)}\left(\frac{2}{|x-y|}\right)^{\frac{\gamma}{\pi}}
\,\dd \hau^1(x)\right\}\frac{|F(y)|}{\|F\|_{L^1(B_1(0))}}\,\dd y 
\cdot \e^{\frac{\gamma}{4\pi}}.
\end{align*}
Then one sees that,
for $\gamma<\pi$, the integral:
\begin{align*}
\int_{\partial B_1(0)}\left(\frac{2}{|x-y|}\right)^{\frac{\gamma}{\pi}}
\,\dd \hau^1(x)
\end{align*} is convergent and its value uniformly bounded in $y\in B_1(0)$.
We conclude that:
\begin{align}
\label{eq:exp-u-1}
\int_{\partial B_1(0)}
\exp\left(\frac{\gamma|u_1(x)|}{\|F\|_{L^1(D)}}\right)\,\dd \hau^1(x)
\leq C_\gamma ,\quad\text{for }\gamma<\pi,
\end{align}
for some constant $C_\gamma>0$ depending on $\gamma.$\par
\emph{Step 2: study of $u_2$.}
Note  that for $x=\e^{i\phi}$, $y=\e^{i\theta}$, $\mathcal{G}(x,y)$
takes the form:
\begin{align}
\label{eq:definition-G-on-boundary}
\mathcal{G}(\e^{i\phi},\e^{i\theta})
=-\frac{1}{\pi}\log|x-y|
=-\frac{1}{2\pi}\log(2(1-\cos(\phi-\theta)))
=:G(\phi-\theta),
\end{align}
hence the boundary value of $u_2$ can be written as:
\begin{align*}
u_2(\phi)=u_2(\e^{i\phi})
=G*(g+k) = G*g,\quad\text{on }\partial B_1(0),
\end{align*}
where the last equality follows since $G$ has zero average.
Using again this property, we may write:
\begin{align*}
 u_2(\e^{i\phi})
&=-\frac{1}{2\pi}\int_{\partial B_1(0)}
\log(2(1-\cos(\phi-\theta)))g(\theta)\,\dd\theta\\
&=\frac{1}{2\pi}
\int_{\partial B_1(0)}
\log\left(\frac{1-\cos(\phi-\theta)}{2}\right)
g(\theta)\,\dd\theta,
\end{align*}
where now the argument in the logarithm is always bigger than $1$.
As in step 1, for $\gamma>0$ and Jensen's inequality
one deduces:
\begin{align*}
\exp\left(
\frac{\gamma u_1(\e^{i\theta})}{\|g\|_{L^1(\partial B_1(0))}}\right)
\leq \int_{\partial B_1(0)}
\left(\frac{2}{1-\cos(\phi-\theta)}\right)^{\frac{\gamma}{2\pi}}
\frac{|g(\theta)|}{\|g\|_{L^1(\partial B_1(0))}}\,\dd\theta,
\end{align*}
hence with Tonelli's theorem:
\begin{align*}
\int_{\partial B_1(0)}\exp\left(
\frac{\gamma u_1(\e^{i\theta})}{\|g\|_{L^1(\partial B_1(0))}}\right)\,\dd\theta
\leq \int_{\partial B_1(0)}
\left\{
\int_{\partial_{B_1(0)}}\left(\frac{2}{1-\cos(\phi-\theta)}\right)^{\frac{\gamma}{2\pi}}
\,\dd\phi\right\}
\frac{|g(\theta)|}{\|g\|_{L^1(\partial B_1(0))}}\,\dd\theta.
\end{align*}
Provided $\gamma<\pi$ the integral:
\begin{align*}
\int_{\partial B_1(0)}\left(\frac{2}{1-\cos(\phi-\theta)}\right)^{\frac{\gamma}{2\pi}}
\,\dd\phi
=\int_{\partial B_1(0)}\left(\frac{2}{1-\cos(\phi)}\right)^{\frac{\gamma}{2\pi}}
\,\dd\phi
\end{align*}
is convergent, hence
we conclude that:
\begin{align}
\label{eq:exp-u-2}
\int_{\partial B_1(0)}
\exp\left(\frac{\gamma| u_2(x)|}{\|g\|_{L^1(\partial B_1(0))}}\right)\,\dd x
\leq C_\gamma,\quad\text{for }\gamma<\pi,
\end{align}
and for some $C_\gamma>0$.\par
\emph{Step 3.}
Finally from \eqref{eq:u-bar-u},
we may write
\begin{align*}
\e^{u-\overline{u}}
=\e^{u_1}\e^{u_2}
\quad\text{on }\partial_{B_1(0)}.
\end{align*}
In particular, if $\frac{1}{p}=\frac{1}{p_1}+\frac{1}{p_2}$, by H\"older's inequality
there holds 
\begin{align*}
\|\e^{u-\overline{u}}\|_{L^p(\partial(B_1(0)))}\leq\|\e^{u_1}\|_{L^{p_1}(\partial B_1(0))}
\|\e^{u_2}\|_{L^{p_2}(\partial B_1(0))}.
\end{align*}
Choosing
\begin{align*}
p_1=\frac{\pi-\varepsilon}{2\|f\|_{L^1(B_1(0))}}
\quad\text{and}\quad
p_2=\frac{\pi-\varepsilon}{\|g\|_{L^1(\partial B_1(0))}},
\end{align*}
we reach  the conclusion by using estimates 
\eqref{eq:exp-u-1} and \eqref{eq:exp-u-2} with $\gamma=\pi-\varepsilon$.
This proves theorem \ref{thm:global-exp-int}.
\end{bfproof}

\begin{remark}
\label{rmk:invar-transl}
In step 2, $G$ defined in \eqref{eq:definition-G-on-boundary} has zero average, the computation is invariant
by translations of $g$.
Consequently, the assumption on $f$ and $g$ on theorem \ref{thm:global-exp-int}
may be replaced by
\begin{align*}
\|g-\alpha\|_{L^1(\partial B_1(0))}+2\|f\|_{L^1(B_1(0))}<\pi-\varepsilon
\quad\text{for some }\alpha\in\R.
\end{align*}
\end{remark}

The following is a localised version of
theorem \ref{thm:global-exp-int}.
\begin{lemma}
\label{lemma:local-exp-int}
Let $f\in L^1(B_1(0))$, $g\in L^1(\partial B_1(0))$,
$a_1,\ldots a_N$ be points in $B_1(0)$
and $\alpha_1,\ldots,\alpha_N$ be real numbers 
satisfying $\int_{B_1(0)} f +\sum_i \alpha_i = \int_{\partial B_1(0)}g$.
Let $u\in W^{1,1}(B_1(0))$ be a weak solution to the problem
\begin{align}
\label{eq:Neumann-dirac}
\left\{
\begin{aligned}
-\Delta u&=f+\sum_{i=1}^N\alpha_i\delta_{a_i} &&\text{in }B_1(0),\\
\partial_\nu u &=g &&\text{on }\partial B_1(0).
\end{aligned}
\right.
\end{align}
Assume that, for a given $x_0\in\partial B_1(0)$ and $0<r<1$,
$ B_r(x_0)\cap\{a_1,\ldots a_N\}=\emptyset$   and 
   $
\|g\|_{L^1(\partial B_1(0)\cap B_r(x_0))}
+2\|f\|_{L^1(B_1(0)\cap B_r(x_0))}<\pi-\varepsilon
$, for some  $0<\varepsilon<\pi$.
Then 
\begin{align}
\label{eq:exp-u}
\|\e^{u-\overline{u}}\|_{L^p(\partial B_1(0)\cap B_{r/2}(x_0))}
\leq C_\varepsilon\,
C_1
\left(\frac{1}{r}\right)^{C_2},
\end{align}

where    \begin{align*}
p=\frac{\pi-\varepsilon}{\|g\|_{L^1(\partial B_1(0)\cap B_r(x_0))}
+2\|f\|_{L^1(B_1(0)\cap B_r(x_0))}},
\end{align*}
  
 $C_\varepsilon>0$ is a constant depending on $\varepsilon$ and $C_1,C_2$ depend on $\|f\|_{L^1(B_1(0))},$ $\|g\|_{L^1(\partial B_1(0))})$ and the $\alpha_i's$.  \end{lemma}

\begin{bfproof}[Proof of Lemma \ref{lemma:local-exp-int}]
Let $\chi:B_1(0)\to\R$ be a function in $C^\infty(\overline{B_1(0)})$ so that
$\chi=1$ in $B_1(0)\cap B_{3r/4} (x_0)$ and with support in $B_1(0)\cap B_r (x_0)$.
We write the solution of \eqref{eq:Neumann-dirac} as:
$
 u-\overline{ u} =u_1+u_2,
$
where:
\begin{align*}
\left\{
\begin{aligned}
-\Delta u_1 &=f\chi &&\text{in }B_1(0),\\
\partial_\nu u_1 & =g\chi -c &&\text{on }\partial B_1(0),\\
\overline{u}_1 &=0,
\end{aligned}
\right.
\quad\text{and}\quad
\left\{
\begin{aligned}
-\Delta u_2 &=f(1-\chi) + \sum_{i-1}^N\alpha_i\delta_{a_i} &&\text{in }B_1(0),\\
\partial_\nu u_2 & =g(1-\chi) +c &&\text{on }\partial B_1(0),\\
\overline{u}_2 &=0,
\end{aligned}
\right.
\end{align*}
 with $c=\frac{1}{\pi}\int_{B_1(0)} f\chi + \frac{1}{2\pi}\int_{\partial B_1(0)}g\chi$.
Applying theorem 
\ref{thm:global-exp-int} together with remark \ref{rmk:invar-transl},
we deduce the existence of $C_\varepsilon>0$ such that:
\begin{align}
\label{eq:estimate-alpha1}
\|\e^{u_1}\|_{L^p(\partial B_1(0))}\leq C_\varepsilon.
\end{align}
To estimate $u_2$ we use the representation formula
 
\begin{align*}
u_2(x)
=\int_{B_1(0)} \mathcal{G}(x,y)f(y)(1-\chi(y))\,\dd y 
+\int_{\partial B_1(0)} \mathcal{G}(x,y)g(y)(1-\chi(y))\,\dd\hau^1(y)
+\sum_{i=1}^{N}\alpha_i\mathcal{G}(x,a_i).
\end{align*}
Notice that for $x\in\partial B_1(0)$ we have
$\e^{\alpha_i\mathcal{G}(x,a_i)}
=|x-a_i|^{-\alpha_i/\pi}\e^{\alpha_i(|a_i|^2-1)/4\pi}$.
Since none of the $a_i$'s is in $B_1(0)\cap B_r(x_0)$ we have
  $r/2\leq |x-a_i|\leq 2$ for $x\in \partial B_1(0)\cap B_{r/2}(x_0)$ and 
$|x-a_i|^{-\alpha_i/\pi}\lesssim r^{-|\alpha_i|/\pi}$. Observe also that
and $1-\chi$ vanishes in  $B_1(0)\cap B_{3r/4}(x_0)$, therefore for $x\in \partial B_1(0)\cap B_{r/2}(x_0)$ we have the estimate
 \begin{align}
\label{eq:exp-u2}
\e^{|u_2(x)|}\le C_1 \left(\frac{1}{r}\right)^{C_2}
\end{align}
where $C_1$ and $C_2$ depend on $ \|f\|_{L^1(B_1(0))},$ $\|g\|_{L^1(\partial B_1(0))})$ and 
$\sum_{i=1}^N|\alpha_i|$.
Hence
joining estimates \eqref{eq:estimate-alpha1} and \eqref{eq:exp-u2},
we then deduce the validity of \eqref{eq:exp-u}.
This proves the lemma.
\end{bfproof}

\section{Existence of a Minimiser for the Total Curvature Energy}
\label{sec:exitence-minimiser}
This section is devoted to prove the main part of theorem \ref{thm:existence-solution-Poisson-problem}, namely the existence    for an admissible
triple $(\Gamma,\vec{n}_0,A)$ of  a minimiser for the total curvature energy
\eqref{eq:total-curvature-Sigma} in the class $\mathcal{F}_{B_1(0)}(\Gamma,\vec{n}_0,A)$ introduced in definition
\ref{def:weak-immersions}.
In the section \ref{sec:regularity} it is shown that  such a minimiser is actually $C^1$ and not only Lipschitz .\smallskip

We need several preliminary lemmas.
Along this section, $\Gamma$ and $\vec{n}_0$ will be fixed as in the statement of theorem
\ref{thm:existence-solution-Poisson-problem},
$\gamma:[0,\hau^1(\Gamma)]/\sim\to\Gamma$ will denote
a fixed arc-length parametrization of $\Gamma$
and $k_g$ its geodesic curvature defined in \eqref{eq:geodesic-curvature}.
When dealing with a sequence of maps $(\Phi)_k\subset \mathcal{F}_{B_1(0)}(\Gamma,\vec{n}_0,A)$,
we denote with a subscript $k$ every quantity pertaining to the immersion $\Phi_k$
(e.g. the Gauss map $\vec{n}_{\Phi_k}$ will be simply denoted by $\vec{n}_k$).

\begin{lemma}
\label{lemma:uniform-control-quantities}
Let $(\Phi_k)_k$ be a sequence in of conformal maps in $\mathcal{F}_{B_1(0)}(\Gamma,\vec{n}_0)$
with 
$
E=\sup_{k}E(\Phi_k)<+\infty.
$
Then:
\begin{enumerate}[(i)]
\item the $L^1$-norm of the Gauss curvature $\|K_k\e^{2\lambda_k}\|_{L^1(B_1(0))}$
is uniformly bounded,
\item the number of branch points of $\Phi_k$ and their multiplicity
is uniformly bounded on $k$,
\item The $L^{(2,\infty)}$-norm of the gradient of   the conformal factor 
$\|\nabla\lambda_k\|_{L^{(2,\infty)}(B_1(0))}$ is uniformly bounded on $k$,
\item for any fixed $E_0>0$, the set of points $x\in\overline{B_1(0)}$ such that the curvature energy concentrates
above the level $E_0$, i.e.:
\begin{align*}
\liminf_{k\to+\infty}\left(\inf\{r>0:\|\nabla\vec{n}_k\|^2_{L^2(B_1(0)\cap B_r(x))}\geq E_0\}\right)=0
\end{align*}
is finite and its cardinality is uniformly bounded on $k$.
\end{enumerate}
All such bounds depend only on $E$
and on $\|k_g\|_{L^1(\Gamma)}$.
\end{lemma}
\begin{bfproof}[Proof of Lemma \ref{lemma:uniform-control-quantities}]
\emph{Proof of (i)}. The bound follows from the pointwise a.e. relation
(see appendix \ref{app:frames} for more  details)
\begin{align}
\label{eq:pointwise-estimate-Gauss-n}
K_k\e^{2\lambda_k}\leq\frac{|\nabla\vec{n}_k|^2}{2}.
\end{align}

\emph{Proof of (ii)}. Since $\lambda_k$ is a weak solution to the Liouville's equation
\begin{align}
\label{eq:liouville-k}
\left\{
\begin{aligned}
-\Delta \lambda_{k} &=  K_k\e^{2\lambda_k}-2\pi\sum_{i_k=1}^{N_k}n_{i_k}\delta_{a_{i_k}} && \text{in}\, B_1(0),\\
\partial_{\nu}\lambda_{k}&=k_g(\sigma_k)\e^{\lambda_k}-1 &&\text{on}\, \partial B_1(0),
\end{aligned}\right.
\end{align}
there must hold
\begin{align*}
-2\pi\sum_{i_k=1}^{N_k}n_{i_k}
+\int_{B_1(0)}K_k\e^{2\lambda_k}\,\dd x
=\int_{\partial B_1(0)}(k_g(\sigma_k)\e^{\lambda_k}-1)\,\dd\hau^1,
\end{align*}
hence
\begin{align*}
\sum_{i_k=1}^{N_k}|n_{i_k}|
\leq C\big(\|K_k\e^{2\lambda_k}\|_{L^1(B_1(0))} + \|k_g\|_{L^1(\Gamma)}+1\big),
\end{align*}
and the result then follows from (i).

\emph{Proof of (iii)}. Using Green's representation formula
(the green function $\mathcal{G}$ is given in \eqref{eq:Green-for-Neumann}),
we have
\begin{align*}
\nabla\lambda_k(x)
&=\int_{B_1(0)}\nabla_x\mathcal{G}(x,y)K_k(y)\e^{2\lambda_k(y)}\,\dd y
-2\pi\sum_{i_k=1}^{N_k}n_{i_k}\nabla_x\mathcal{G}(x,a_{i_k})\\
&\phantom{{}={}}
+\int_{\partial B_1(0)}\nabla_x\mathcal{G}(x,y)k_g(\sigma_k(y))\e^{\lambda_k(y)}\,\dd\hau^1(y),
\end{align*}
and since there holds
\begin{align*}
\sup_{y\in B_1(0)}\|\nabla_x\mathcal{G}(\cdot,y)\|_{L^{(2,\infty)}(B_1(0))}<+\infty,
\end{align*}
we can estimate
\begin{align}\label{eq:estgradlambda}
\|\nabla\lambda_k\|_{L^{(2,\infty)}(B_1(0))}
\leq C\big(
\|K_k\e^{2\lambda_k}\|_{L^1(B_1(0))}
+\sum_{i_k=1}^{N_k}|n_{i_k}|
+\|k_g\|_{L^1(\Gamma)} + 1
\big),
\end{align}
hence deduce that the right-hand-side of \eqref{eq:estgradlambda} does not depend on $k$ thanks to (i) and (ii).

\emph{Proof of (iv).}
The proof of (iv) follows from standard concentration-compactness arguments and we omit it.
 \end{bfproof}

Combining lemmas \ref{lemma:Helein} and \ref{lemma:uniform-control-quantities}, 
the following estimate
is obtained for points in $B_1(0)$ (see \cite[Theorem 4.5]{RivierePCMI}).
\begin{lemma}
\label{lemma:lemma:interior-estimate-conf-fact}
There exists an $\varepsilon_0>0$ with the following property.
Let $(\Phi_k)_k$ be a sequence of conformal maps in $\mathcal{F}_{B_1(0)}(\Gamma,\vec{n}_0)$
with 
$
E=\sup_{k}E(\Phi_k)<+\infty
$
and let $x_0\in B_1(0)$ and $0<r<1$ be so that 
$\overline{B_r(x_0)}\subset B_1(0)$.
If, for every $k\in\N$,
$\overline{B_r(x_0)}$ contains no branch points of $\Phi_k$
and there holds
\begin{align}
\|\nabla\vec{n}_k\|_{L^2(B_r(x_0))}^2\leq\varepsilon
\end{align}
for some $0<\varepsilon<\varepsilon_0$,
then $\lambda_k\in C^0(\overline{B_r(x_0)})\cap W^{1,2}(B_r(x_0))$ and there exist a constant 
$C=C(\Gamma,\vec{n}_0,E,\varepsilon_0,r)>0$ 
and a sequence $(c_k)_{k}\subset\R$
 so that
\begin{align}
\label{eq:estimate-lambda-interior}
\sup_k\bigg(
\|\nabla\lambda_k\|_{L^2(B_r(x_0))}+
\|\lambda_k-c_k\|_{L^\infty(B_{r/2}(x_0))}
\bigg)\leq C.
\end{align}
\end{lemma}

We have the following analogue result for boundary points.

\begin{lemma}
\label{lemma:boundary-estimate-conf-fact}
There exists an $\varepsilon_0>0$ with the following property.
Let $(\Phi_k)_k$ be a sequence of conformal maps in $\mathcal{F}_{B_1(0)}(\Gamma,\vec{n}_0)$
with 
$
E=\sup_{k}E(\Phi_k)<+\infty,
$
and let $x_0\in\partial B_1(0)$ and $0<r<1$.
If, for every $k\in\N$,  $B_1(0)\cap\overline{B_r(x_0)}$ contains no branch points of $\Phi_k$
and, having denoted $\vec{e}_k = \e^{-\lambda_k}(\partial_1\Phi_k,\partial_2\Phi_k)$ 
the ortho-normal frame
associated with $\Phi_k$, there holds
\begin{align}
\label{eq:smallness-2nd-fund}
\|\nabla\vec{n}_k\|_{L^2(B_1(0)\cap B_r(x_0))}^2
+\|k_g(\sigma_k)\e^{\lambda_k}\|_{L^1(\partial B_1(0)\cap B_r(x_0))}\leq\varepsilon
\quad\text{and}\quad
[\vec{e}_k]_{W^{1/2,2}(\partial B_1(0)\cap B_r(x_0))}^2\leq \varepsilon
\end{align}
for some $0<\varepsilon<\varepsilon_0$,
then
$\lambda_k\in C^0(\overline{B_1(0)\cap B_r(x_0)})\cap W^{1,2}(B_1(0)\cap B_r(x_0))$
and there exist constant a constant 
$C=C(\Gamma,\vec{n}_0,E,\varepsilon_0,r)>0$
independent of $\varepsilon$
 and a sequence $(c_k)_{k}\subset\R$  
so that
\begin{align}
\label{eq:estimate-lambda-boundary}
\sup_k\bigg(
\|\nabla\lambda_k\|_{L^2(B_1(0)\cap B_{r/4}(x_0))}+
\|\lambda_k-c_k\|_{L^\infty(B_1(0)\cap B_{r/4}(x_0))}\bigg)
\leq C.
\end{align}
\end{lemma}

\begin{bfproof}[Proof of Lemma \ref{lemma:boundary-estimate-conf-fact}]
On the one hand, 
from the pointwise a.e. relation
\eqref{eq:pointwise-estimate-Gauss-n}
and  \eqref{eq:smallness-2nd-fund}
we have
\begin{align*}
\|k_g(\sigma_k)\e^{\lambda_k}\|_{L^1(\partial B_1(0)\cap B_{r}(x_0))}
+2\|K_k\e^{2\lambda_k}\|_{L^1(B_1(0)\cap B_r(x_0))} \leq\varepsilon
\end{align*}
consequently since $\lambda_k$ is a weak solution to Liouville's equation
\eqref{eq:liouville-k},
by choosing an $\varepsilon_0$ small enough, 
by lemma \ref{lemma:local-exp-int}
 and lemma \ref{lemma:uniform-control-quantities}
we may find a $p=p(\varepsilon_0)>1$ so that,
uniformly on $k$, there holds
\begin{align}
\label{eq:a-priori-estimate-exp}
\|\e^{\lambda_k-\overline{\lambda}_k}\|_{L^p(\partial B_1(0)\cap B_{r/2}(x_0))}
\leq C.
\end{align}
On the other hand, possibly after reducing $\varepsilon_0$
we can invoke lemma \ref{lemma:local-frame-prescribed-boundary}
and deduce the existence of Coulomb ortho-normal frames
$\vec{g}_k=(\vec{g}_{k,1},\vec{g}_{k,2})\in W^{1,2}(B_1(0)\cap B_r(x_0))$
lifting $\vec{n}_k$ in $B_1(0)\cap B_r(x_0)$,
coinciding with $\vec{e}_k$ on $\partial B_1(0)\cap B_r(x_0)$
and so that uniformly in $k$ there holds
\begin{align}
\label{eq:bound-energy-controlled-frame}
\|\nabla\vec{g}_k\|^2_{L^2(B_1(0)\cap B_r(x_0))}\leq C.
\end{align}
In particular, we may write:
\begin{align*}
K_k\e^{2\lambda_k} &= \langle \nabla^\perp \vec{g}_{k,1},\nabla\vec{g}_{k,2}\rangle
\text{ in }B_1(0)\cap B_{r}(x_0),\\
k_g(\sigma_k)\e^{\lambda_k}-1 &= \langle\partial_\tau \vec{g}_{k,1},\vec{g}_{k,2}\rangle
\text{ on }\partial B_1(0)\cap B_{r}(x_0).
\end{align*}
From lemma \ref{lemma:estimate-localised-Wente-Neumann},
we deduce that $\lambda_k\in C^0(\overline{B_1(0)\cap B_{r/4}(x_0)})\cap W^{1,2}(B_1(0)\cap B_{r/4}(x_0))$
and that for some constant $c_k\in\R$ there holds
\begin{align}
\label{eq:estimate-lambda-k-long}
&\|\lambda_k-c_k\|_{L^{\infty}(B_1(0)\cap B_{r/4}(x_0))}
{+
\|\nabla\lambda_k\|_{L^2(B_1(0)\cap B_{r/4}(x_0))}}\\
&\notag
\phantom{{}={}}\leq C
\bigg(\|K_k\e^{2\lambda}\|_{L^1(B_1(0))} 
+ \|k_g(\sigma_k)\e^{\lambda_k}\|_{L^1(\partial B_1(0))} + \sum_{i_k=1}^{N_k} |\alpha_{i_k}|\\
\notag&\phantom{{}\leq{}}+
\|\nabla\vec{g}_{k,1}\|_{L^2(B_1(0)\cap B_{r/2}(x_0))}
\|\vec{g}_{k,2}\|_{W^{1,2}(B_1(0)\cap B_{r/2}(x_0))} \\
\notag&\phantom{{}\leq{}}+
\|\partial_\tau\vec{g}_{k,1}|_{\partial B_1(0)}\|_{L^p(\partial B_1(0)\cap B_{r/2}(x_0))}
\|\vec{g}_{k,2}|_{\partial B_1(0)}\|_{W^{1,p}(\partial B_1(0)\cap B_{r/2}(x_0))}
\bigg).
\end{align}
The first line on the right hand side of \eqref{eq:estimate-lambda-k-long} can be estimated
uniformly on $k$ by means of lemma \ref{lemma:uniform-control-quantities}-(i)-(ii).
The second line can be estimated
uniformly on $k$  with \eqref{eq:bound-energy-controlled-frame}.
Finally the third line is estimated
uniformly on $k$ with \eqref{eq:a-priori-estimate-exp}
since, for $i=1,2$, we have
\begin{align}\label{estgki}
\|\vec{g}_{k,i}|_{\partial B_1(0)}\|_{W^{1,p}(\partial B_1(0)\cap B_{r/2}(x_0))}
&\leq C\big(
\|(|k_g(\sigma_k)|+|\dot{\vec{n}}_0(\sigma_k)|)\e^{\lambda_k}
\|_{L^p(\partial B_1(0)\cap B_{r/2}(x_0))}+1\big)\nonumber\\
&\leq C\big(
\|\e^{\lambda_k-\overline{\lambda}_k}\|_{L^p(\partial B_1(0)\cap B_{r/2}(x_0))}
+1\big).
\end{align}
In \eqref{estgki} we use the fact that $k_g(\sigma_k),\dot{\vec{n}}_0(\sigma_k)\in L^{\infty}(\partial{B_1(0)}$ by the current assumptions on the boundary data and  that by Jensen's inequality it holds uniformly on $k$:
\begin{align}
\label{eq:Jensen-lambda}
\e^{\overline{\lambda}_k}
=\exp\dashint_{\partial{B_1(0)}}\lambda_k\leq
\dashint_{\partial{B_1(0)}}\exp(\lambda_k)=\frac{\hau^1(\Gamma)}{2\pi}.
\end{align}
This concludes the proof of the lemma.
\end{bfproof}

\begin{definition}
\label{def:three-point-cond}
Let $P_1,P_2,P_3$ be three distinct, fixed, consecutive points
in  $P_1,P_2,P_3$ on $\Gamma$
that is, $\gamma(s_j)=P_j$ for some $0\leq s_1<s_2<s_3<\mathcal{H}^1(\Gamma)$.
We denote by 
$\mathcal{F}^*_{B_1(0)}(\Gamma,\vec{n}_0,A)$ is the set of maps $\Phi\in\mathcal{F}_{B_1(0)}(\Gamma,\vec{n}_0,A)$
so that
\begin{align}
\label{eq:three-pt-cond}
\Phi\left(\e^{\frac{2\pi i }{3}j}\right)=P_j
\quad\text{for }j=1,2,3.
\end{align}
\end{definition}

\begin{remark}
\label{rmk:moebius-equivalent}
We note that:
\begin{enumerate}[(i)]
\item if $\sigma_\Phi$ defines the boundary parametrization of $\Phi$,
that is $\Phi|_{\partial B_1(0)}=\gamma\circ\sigma_\Phi$,
condition \eqref{eq:three-pt-cond} is equivalent to
$\sigma_\Phi\left(\frac{2\pi}{3}j\right)= s_j$ for $j=1,2,3$.

\item For every $\Phi\in \mathcal{F}_{B_1(0)}(\Gamma,\vec{n}_0,A)$,
there is a unique M\"obius transformation of $B_1(0)$ so that
$\Phi\circ \phi\in\mathcal{F}^*_{B_1(0)}(\Gamma,\vec{n}_0,A)$.
Moreover, the invariance by diffeomorphisms of the total curvature energy
implies that $E(\Phi)=E(\Phi\circ\phi)$, hence
\begin{align*}
\inf_{\Phi\in\mathcal{F}^*_{B_1(0)}(\Gamma,\vec{n}_0,A)}E(\Phi)
=\inf_{\Psi\in\mathcal{F}_{B_1(0)}(\Gamma,\vec{n}_0,A)}E(\Psi).
\end{align*}
\end{enumerate}
\end{remark}

The following lemma is a consequence of the Courant-Lebesgue lemma,
a key tool in the analysis of Plateau's problem (see
e.g. \cite[Lemma 4.14]{MR2780140}).

\begin{lemma}
\label{lemma:equicontinuity-from-Courant-Lebesgue}
For any sequence $(\Phi_k)_{k}\subset\mathcal{F}^*_{B_1(0)} (\Gamma,\vec{n}_0,A)$,
the sequence of boundary curves $(\Phi_k|_{\partial B_1(0)})_{k}$ is equicontinuous.
\end{lemma}
Equicontinuity of the boundary curves is equivalent to the equicontinuity of the $\sigma_k$'s.
As a consequence of lemma \ref{lemma:equicontinuity-from-Courant-Lebesgue}, we have:

\begin{lemma}
\label{lemma:technical-noncomplication}
Let $(\Phi_k)_{k}$ be a sequence in $\mathcal{F}^*_{B_1(0)} (\Gamma,\vec{n}_0,A)$
and let $x_0\in\partial B_1(0)$ be fixed.
Then, possibly passing to a subsequence, for any any $\varepsilon>0$,
there always exists an $r=r(\Gamma,\vec{n}_0)>0$ so that:
\begin{align}
\label{eq:uniform-estimate-e-k-g-three-pts}
\sup_{k\in\N}\,[\vec{e}_k]_{W^{1/2,2}(\partial B_1(0)\cap B_r(x_0))}
\leq \varepsilon
\quad\text{and}\quad
\sup_{k\in\N}
\|k_g(\sigma_k)\e^{\lambda_k}\|_{L^1(\partial B_1(0)\cap B_r(x_0))}\leq \varepsilon.
\end{align}
\end{lemma}
\begin{bfproof}[Proof of Lemma \ref{lemma:technical-noncomplication}]
   Possibly after extracting a subsequence, thanks to
lemma \ref{lemma:equicontinuity-from-Courant-Lebesgue} and the Arzel\`a-Ascoli
theorem, we may suppose that $\sigma_k$ converges uniformly on $\partial B_1(0)$ to some 
continuous map $\sigma$.
As a consequence,
we have the pointwise convergence away from the diagonal:
\begin{align*}
\lim_k\left(
 \frac{\sigma_k(\theta_1)-\sigma_k(\theta_2)}{\e^{i\theta_1}-\e^{i\theta_2}}
\right)
=
\frac{\sigma(\theta_1)-\sigma(\theta_2)}{\e^{i\theta_1}-\e^{i\theta_2}}
\quad\text{for }\theta_1\neq\theta_2,
\end{align*}
and, possibly after extracting another subsequence,
the bound:
\begin{align*}
\frac{|\sigma_k(\theta_1)-\sigma_k(\theta_2)|^2}
{|\e^{i\theta_1}-\e^{i\theta_2}|^2}
\leq 2
\frac{|\sigma(\theta_1)-\sigma(\theta_2)|^2}
{|\e^{i\theta_1}-\e^{i\theta_2}|^2}.
\end{align*}
Hence, integrating both sides,
we deduce that for every $\rho>0$ and $x\in\partial B_1(0)$
there holds:
\begin{align}
\label{eq:uniform-control-trace-sigma}
[\sigma_k]_{W^{1/2,2}(\partial B_1(0)\cap B_\rho(x))}^2
\leq 2 [\sigma]_{W^{1/2,2}(\partial B_1(0)\cap B_\rho(x))}^2.
\end{align}

Let us assume without loss of generality that $x_0=1$ and that $r<1$, so that we
may identify $[-\theta_0,\theta_0]\simeq \partial B_1(0)\cap B_r(x_0)$
for some $0<\theta_0<\pi$.
Writing in complex notation 
\begin{align*}
\e^{i\theta}(\vec{e}_{k,1}+i\vec{e}_{k,2})
=\star(\mathbf{t}\wedge\vec{n}_0)(\sigma_k)+i\mathbf{t}(\sigma_k)
\quad\text{on }\partial B_1(0),
\end{align*}
thanks to \eqref{eq:uniform-control-trace-sigma}
we deduce that
\footnote{
Recall the elementary inequality: 
\begin{align*}
[ab]_{W^{1/2,2}}^2\leq 2\left(\|a\|_{L^\infty}^2[b]^2_{W^{1/2,2}}
+\|b\|_{L^\infty}^2[a]^2_{W^{1/2,2}}\right).
\end{align*}
Also recall that if $a:(l_1,l_2)\to \C$ is a Lipschitz function
and $b:(-\theta_0,\theta_0)\to(l_1,l_2)$ is a function in $H^{1/2}((-\theta_0,\theta_0))$
we have:
\begin{align*}
[a\circ b]^2_{W^{1/2,2}((-\theta_0,\theta_0))}
&\leq [a]_{C^{0,1}((l_1,l_2))}^2[b]^2_{W^{1/2,2}((-\theta_0,\theta_0))}.
\end{align*}
}
that
\begin{align*}
[\vec{e}_k]^2_{W^{1/2,2}((-\theta_0,\theta_0))}
&\leq 2\left(
[\mathbf{t}(\sigma_k)]^2_{W^{1/2,2}((-\theta_0,\theta_0))}
+[\star(\mathbf{t}\wedge\vec{n}_0)(\sigma_k)]^2_{W^{1/2,2}((-\theta_0,\theta_0))}
\right) + 4\theta_0^2\\
&\leq2\left( [\mathbf{t}]^2_{C^{0,1}}
+[\star(\mathbf{t}\wedge\vec{n}_0)]^2_{C^{0,1}}\right)
[\sigma_k]^2_{W^{1/2,2}((-\theta_0,\theta_0))} +4\theta_0^2\\
&\leq C\left([\sigma]^2_{W^{1/2,2}((-\theta_0,\theta_0))} +\theta_0^2\right),
\end{align*}
and so the first inequality in
\eqref{eq:uniform-estimate-e-k-g-three-pts}
follow by choosing  a sufficiently small $\theta_0$.
As for the second inequality in \eqref{eq:uniform-estimate-e-k-g-three-pts}, if $a$ and $b$ denote the extrema of 
$\partial B_1(0)\cap B_r(x_0)$, from the point-wise convergence of $\sigma_k$ to $\sigma$
we have:
\begin{align*}
\lim_k\int_{\partial B_1(0)\cap B_r(x_0)}|k_g(\sigma_k)|\e^{\lambda_k(\e^{i\theta})}\,\dd\theta
=\lim_k\int_{\sigma_k(a)}^{\sigma_k(b)}|k_g(s)||\dot{\gamma}(s)|\,\dd s
= \int_{\sigma(a)}^{\sigma(b)}|k_g(s)||\dot{\gamma}(s)|\,\dd s,
\end{align*}
hence, possibly after extracting a subsequence, there holds: 
\begin{align*}
\|k_g(\sigma_k)\sigma_k'\|_{L^1(\partial B_1(0)\cap B_r(x_0))}
\leq 2\|k_g(\sigma)\sigma'\|_{L^1(\partial B_1(0)\cap B_r(x_0))},
\end{align*}
and the results then follows by choosing
$r$ sufficiently small.
This concludes the proof of the lemma.
\end{bfproof}

\begin{definition}[Weak Sequential convergence]
\label{def:weak-convergence}
Given a sequence $(\Phi_k)_{k}$ of conformal maps in $\mathcal{F}_{B_1(0)}(\Gamma,\vec{n}_0)$,
and a conformal  map $\Phi:B_1(0)\to\R^m$,
we say that $\Phi_k$ \emph{weakly converges}  to $\Phi$ if:
\begin{enumerate}[(i)]
\item $
\Phi_k\rightharpoonup\Phi\quad\text{in }W^{1,2}(B_1(0),\R^m)$ and a.e. on $B_1(0)$,
\item$
\Phi_k|_{\partial B_1(0)} \to 
\Phi|_{\partial B_1(0)}\quad\text{uniformly in }C^0(\partial B_1(0)),$
\item $\nabla\lambda_k\overset{*}{\rightharpoonup}
\nabla\lambda\text{ in }L^{(2,\infty)}(B_1(0)),$
\end{enumerate}
and there exists a finite, possibly empty set $\underline{\eta}=\{\eta_1,\ldots,\eta_N\}\subset\overline{B_1(0)}$
so that, for every open set $\Omega\subset\R^2$ with compact closure in
$\overline{B_1(0)}\setminus\underline{\eta}$, there holds:
\begin{enumerate}[(i)]
\setcounter{enumi}{3}
\item$
\lambda_k\overset{*}{\rightharpoonup}
\lambda\text{ in }L^\infty(\Omega,\R^n)$
and ,
\item$
\Phi_k\rightharpoonup\Phi\quad\text{in }W^{2,2}(\Omega,\R^n),$
\end{enumerate}
where $\lambda=\log(|\nabla\Phi|/\sqrt{2})$.
\end{definition}

We are now in the position to prove the following compactness result.
\begin{lemma}
\label{lemma:compactness}
Let $(\Phi_k)_k$ be a sequence
of conformal maps in $\mathcal{F}^*_{B_1(0)}(\Gamma,\vec{n}_0,A)$ 
with
$
\sup_{k}E(\Phi_k)<+\infty.
$
Then $(\Phi_k)_k$ contains a subsequence weakly converging in the sense of
definition \ref{def:weak-convergence}
to an element $\Phi\in\mathcal{F}^*_{B_1(0)}(\Gamma,\vec{n}_0,A)$.
\end{lemma}

\begin{bfproof}[Proof of Lemma \ref{lemma:compactness}]
\emph{Step 1.}
Since for every $k$ we have that $\|\nabla\Phi_k\|^2_{L^2(B_1(0))}=2A$
and the three-point condition \eqref{eq:three-pt-cond} holds,
we deduce that $\sup_k(\|\Phi_k\|_{W^{1,2}(B_1(0))})$ is finite and so by the Rellich-Kondrachov 
theorem, possibly passing to a subsequence, condition $(i)$
of definition \ref{def:weak-convergence} is satisfied.

\emph{Step 2.}
From lemma \ref{lemma:equicontinuity-from-Courant-Lebesgue},
by Arzel\`a-Ascoli theorem we deduce that, possibly passing to a subsequence, condition $(ii)$
of definition \ref{def:weak-convergence} is satisfied.

\emph{Step 3.}
From lemma \ref{lemma:uniform-control-quantities}-(iii),
we deduce that, possibly passing to a subsequence, condition $(iii)$
of definition \ref{def:weak-convergence} is satisfied.

\emph{Step 4.}
If $\{a_{k,1},\ldots a_{k,N_k}\}$ is the set of branch points
of $\Phi_k$, 
from lemma \ref{lemma:uniform-control-quantities},
possibly after extracting a subsequence, we may suppose that $N_k$
is independent of $k$ and that, for each $j=1,\ldots N$,
$\lim_{k\to \infty}a_{k,j}=a_j$ for some $a_j\in\overline {B_1(0)}$.\par

We say that a point $p$ belongs to $\underline{\eta}$ if either:
\begin{itemize}
\item $p=a_k$ for some $k=1,\ldots, N$, or
\item there holds
\begin{align}
\label{eq:energy-concentr-cond}
\liminf_{k\to\infty}
\left(\inf\left\{
r>0: 
\|\nabla\vec{n}_k\|^2_{L^2(D\cap B_r(p))}\geq \varepsilon_0
\right\}\right)=0,
\end{align}
with $\varepsilon_0$
is as in lemma \ref{lemma:lemma:interior-estimate-conf-fact} if $p\in B_1(0)$,
or
as in lemma \ref{lemma:boundary-estimate-conf-fact} if $p\in \partial B_1(0)$, or
\item $p=\e^{\frac{2\pi i}{3}j}$ for $j=1,2$ or $3$.
\end{itemize}
Note that the set of points satisfying \eqref{eq:energy-concentr-cond}
is, due to lemma \ref{lemma:uniform-control-quantities}, finite and uniformly bounded in $k$.
Let $\Omega\subset\R^2$ an open set with compact closure in
$\overline{B_1(0)}\setminus\underline{\eta}$ 
and let $\Omega'$ be a closed set contained in $\overline{B_1(0)}\setminus \underline{\eta}$
and with smooth boundary
so that $\overline{\Omega}\subset\Omega'$ and for some small $\delta>0$ there holds
\begin{align*}
\overline{B_1(0)}\setminus\cup_{p\in\underline{\eta}} B_{2\delta}(p)
\subset\subset\Omega'
\subset\subset
\overline{B_1(0)}\setminus\cup_{p\in\underline{\eta}} B_{\delta}(p).
\end{align*}
Possibly passing to a further subsequence,
we may suppose  that, for every $k$,
the set of branch points $\{a_{k,1},\ldots a_{k,N_k}\}$ of $\Phi_k$, lies in 
$\cup_{p\in\underline{\eta}} B_\delta(p)$.
Now, for every $x\in \overline{\Omega'}$, we can choose an
$r_x>0$ so that,
if $x\in B_1(0)$, then $\overline{B_{r_x}(x)}\subset\Omega'$ and 
$\|\nabla \vec{n}_k\|_{L^2( B_{r_x}(x))}^2<\varepsilon_0$,
and,
if $x\in\partial B_1(0)$, then
\begin{align*}
\|\nabla\vec{n}_k\|_{L^2(B_1(0)\cap B_{r_x}(x_0))}^2
+\|k_g(\sigma_k)\e^{\lambda_k}\|_{L^1(\partial B_1(0)\cap B_{r_x}(x))}&<\varepsilon_0,\\
[\vec{e}_k]_{W^{1/2,2}(\partial B_1(0)\cap B_{r_x}(x))}^2&<\varepsilon_0,
\end{align*}
(this can be done uniformly on $k$ thanks to lemma \ref{lemma:technical-noncomplication}).
The family $\{B_{r_x/4}(x)\}_{x\in\overline{\Omega'}}$
forms an open cover of $\overline{\Omega'}$, 
from which we may extract a finite sub-cover
$\{B_{r_j}(x_j)\}_{j=1}^M$.
From lemmas \ref{lemma:lemma:interior-estimate-conf-fact} and \ref{lemma:boundary-estimate-conf-fact}, we deduce 
that $\lambda\in C^0(\overline{B_1(0)\cap B_{r_j}(x_j)})
\cap W^{1,2}(B_1(0)\cap B_{r_j}(x_j))$ and there exists constants $l_k(x_j)$,  so that,
for $j=1,\ldots,M$,
\begin{align}
\label{eq:bound-conf-fact-covering}
\sup_{k}\bigg(
\|\nabla\lambda_k\|_{L^2(B_1(0)\cap B_{r_j}(x_j))} +
\|\lambda_k-l_k(x_j)\|_{L^\infty(B_1(0)\cap B_{r_j}(x_j))}\bigg)
\leq C.
\end{align}
Notice that, for every $i,j$ the bound $|l_k(x_i)-l_k(x_j)|\leq M C$ holds.
Indeed, if $B_\rho(x_i)\cap B_\rho(x_j)\neq\emptyset$, then from
\eqref{eq:bound-conf-fact-covering} and the triangle inequality
we have
$
|l_k(x_i)-l_k(x_j)|
\leq 2C.
$
and general $i$ and $j$ pick a collection from the covering which connects
$x_j$ and $x_j$ and reach a similar conclusion.
We can then assume that such constants do not depend on $x_j$ and consequently
that, for some $l_k$, there holds:
\begin{align}
\label{eq:bound-conf-fact-omega}
\sup_k\bigg(
\|\nabla\lambda_k\|_{L^2(\Omega')}+
\|\lambda_k-l_k\|_{L^\infty(\Omega')}
\bigg)\leq C.
\end{align}
We claim that the sequence $(l_k)_k$ is uniformly bounded on $k$.
To see that $(l_k)_k$ is bounded from below, note that, if we had $\limsup_{k\to\infty}l_k=+\infty$,
possibly after extracting a subsequence
condition \eqref{eq:bound-conf-fact-omega} would imply that $\lim_{k\to\infty}\lambda_k=+\infty$ 
uniformly on $\Omega$.
Consequently we would have
\begin{align*}
\lim_{k\to\infty}\int_{\Omega'}\e^{2\lambda_k}\,\dd x=+\infty,
\end{align*}
which contradicts condition $\|\nabla\Phi_k\|^2_{L^2(B_1(0))}=2A$.
Suppose now $(l_k)_k$ is not bounded from above, that is $\liminf_{k\to\infty}l_k=+\infty$.
Let $\alpha$ be an arbitrary closed, connected sub-arc of $\partial B_1(0)$
which does not contain any point in $\underline{\eta}$.
We add if necessary a finite number of balls $B_{r_j}(x_j)$, to the above finite cover of 
$\Omega'$
so that $\{B_{r_j}(x_j)\}_j$ also covers $\alpha$.
Since $\lambda_k$ is continuous, possibly passing to a subsequence, from \eqref{eq:bound-conf-fact-omega},
that $\lim_{k\to\infty}\lambda_k=-\infty$ uniformly in $\alpha$.
Then,
\begin{align*}
0=\lim_{k\to+\infty}\int_\alpha\e^{\lambda_k}\,\dd\sigma 
=\lim_{k\to+\infty}\hau^1(\Phi_k|_{\partial B_1(0)}(\alpha)),
\end{align*}
and thus, by the weak lower-semicontinuity of the Hausdorff
measure with respect to the uniform convergence, that
$\hau^1(\Phi(\alpha ))=0$.
Since the arc $\alpha$ was arbitrarily chosen,
from the Borel-regularity of the Hausdorff measure (see \cite{MR3409135})
we have $\hau^1(\Phi(\partial B_1(0)\setminus\underline{\eta}))=0$,
but then, $\Phi$ being continuous and $\underline{\eta}$ consisting of a finite set,
that $\hau^1(\Phi(\partial B_1(0)))=0$.
This contradicts the three-point normalization condition.
This proves that condition \eqref{eq:bound-conf-fact-omega} can actually be strengthened to
\begin{align}
\label{eq:bound-conf-fact-omega-better}
\sup_k\bigg(
\|\nabla\lambda_k\|_{L^2(\Omega')}+
\|\lambda_k\|_{L^\infty(\Omega')}
\bigg)\leq C.
\end{align}
and thus, possibly passing to a subsequence, condition (iv) of definition \ref{def:weak-convergence}
is satisfied.

\emph{Step 5.}
From \eqref{eq:bound-conf-fact-omega-better}, we can estimate
\begin{align}
\label{eq:bel-bound-1}
\|\nabla\Phi_k\|_{L^\infty(\Omega')}&\leq C,\\
\label{eq:bel-bound-2}
\|\Delta\Phi_k\|_{L^2(\Omega')}^2
\leq\frac{1}{4}\|\e^{2\lambda_k}\|_{L^\infty(\Omega')}W(\Phi_k)&\leq C,
\end{align}
moreover,
\begin{align*}
\|\e^{\lambda_k}\|_{W^{1/2,2}(\partial B_1(0)\cap\Omega')}
&\leq \|\e^{\lambda_k}\|_{W^{1/2,2}(\partial\Omega')}\\
&\leq C \|\e^{\lambda_k}\|_{W^{1,2}(\Omega')}\\
&\leq C\,\e^{\|\lambda_k\|_{L^\infty(\Omega')}}(1+\|\nabla\lambda_k\|_{L^2(\Omega')})
\leq C,
\end{align*}
hence
\begin{align}
\label{eq:bel-bound-3}
\|\partial_\tau\Phi_k\|_{W^{1/2,2}(\partial B_1(0)\cap\Omega')}
&=\|\e^{\lambda_k}\vec{e}_{k,1}(\sigma_k)\|_{W^{1/2,2}(\partial B_1(0)\cap\Omega')}\\
\notag&\leq C(\|\e^{\lambda_k}\|_{W^{1/2,2}(\partial B_1(0)\cap\Omega')}\\
\notag&\phantom{{}\leq C{}}
+\|\e^{\lambda_k}\|_{L^\infty(\partial B_1(0)\cap\Omega')}
[\vec{e}_{k,1}]_{W^{1/2,2}(\partial B_1(0)\cap\Omega')})\\
\notag&\leq C.
\end{align}
From \eqref{eq:bel-bound-1}-\eqref{eq:bel-bound-2}-\eqref{eq:bel-bound-3},
elliptic regularity theory yields $\sup_k\|\Phi_k\|_{W^{2,2}(\Omega)}<+\infty$.
Thus, possibly passing to a subsequence, also condition (iv) of definition \ref{def:weak-convergence}
holds.

\emph{Step 6.}
As in \cite[Lemma 5.1]{RivierePCMI},
we have that the Gauss map of $\Phi$ extends to a map in $W^{1,2}(B_1(0),\Gr_{m-2}(\R^m))$,
and consequently, from lemma \ref{lemma:interior-branch-points}
in appendix \ref{sec:branch-points},
the structure near points of $\underline{\eta}$ lying in $B_1(0)$ is
that of a (possibly removable) branch point.
Finally as shown in lemma \ref{lemma:boundary-singular-points}
singular points $a\in\underline{\eta}$
lying on the boundary
are always removable,
and the limiting map $\Phi$ extends to a conformal
Lipschitz immersion near $\partial B_1(0)$.
This concludes the proof of the lemma.
\end{bfproof}

The proof of the following lemma
can be easily deduced from its analogue in the closed case
(see \cite[Theorem 5.9]{RivierePCMI}).

\begin{lemma}
\label{lemma:lowersemicont-energy}
The total curvature energy functional $E$
is sequentially lower semi-continuous in $\mathcal{F}_{B_1(0)}(\Gamma,\vec{n}_0,A)$
with respect to weak convergence in the sense of definition \ref{def:weak-convergence},
that is,
if $(\Phi_k)_k$ is a sequence in $\mathcal{F}_{B_1(0)}(\Gamma,\vec{n}_0,A)$ weakly
converging to $\Phi$, then
\begin{align*}
\liminf_{k\to\infty} E(\Phi_k)\geq E(\Phi).
\end{align*}
\end{lemma}

\begin{bfproof}[Proof of theorem \ref{thm:existence-solution-Poisson-problem} (first part)]
Since 
the triple $(\Gamma,\vec{n}_0,A)$ is admissible, 
the set $\mathcal{F}_{B_1(0)}(\Gamma,\vec{n}_0,A)$ is not empty and
we can consider a sequence $(\Phi_k)_{k}$ in $\mathcal{F}_{B_1(0)}(\Gamma,\vec{n}_0,A)$
minimising the total curvature energy:
\begin{align*}
\lim_{k\to\infty}E(\Phi_k) 
= \inf\{E(\Psi): \Psi\in \mathcal{F}_{B_1(0)}(\Gamma,\vec{n}_0,A)\}.
\end{align*}
Thanks to lemma \ref{lemma:conformal-coordinates}, 
and remark \ref{rmk:moebius-equivalent},
we may assume that each $\Phi_k$ is conformal
and satisfies the three-point normalisation condition given by
definition \ref{def:three-point-cond}.
From lemma \ref{lemma:compactness} we can then extract a weakly converging subsequence 
in the sense of definition \ref{def:weak-convergence}
to a conformal map $\Phi$ in
 $\mathcal{F}_{B_1(0)}(\Gamma,\vec{n}_0,A)$.
Finally, because of lemma \ref{lemma:lowersemicont-energy} we have
\begin{align*}
E(\Phi) = \lim_{k\to\infty}E(\Phi_k)=
\inf\{E(\Psi): \Psi\in \mathcal{F}_{B_1(0)}(\Gamma,\vec{n}_0,A)\}.
\end{align*}
This concludes the proof of the theorem.
\end{bfproof}

\section{Regularity of Minimisers}
\label{sec:regularity}
This section is devoted to prove
that any element in $\mathcal{F}_{B_1(0)}(\Gamma,\vec{n}_0,A)$
which minimises the total curvature energy satisfies all the regularity
statements of theorem \ref{thm:regularity-solution-Poisson-problem},
thereby also
concluding the proof of theorem \ref{thm:existence-solution-Poisson-problem}.
\smallskip

We need first some preparatory results regarding
immersions with $L^2$-bounded second fundamental form
and estimates on suitable competitors for the Poisson problem.
In this section, we denote
with $D\Area(\Phi)w$ and $DE(\Phi)w$
the directional derivative at $\Phi$
along $w$ of the area and 
total curvature energy functional,
namely
\begin{align*}
D\Area(\Phi)w=\frac{\dd \Area(\Phi+t w)}{\dd t}\bigg|_{t=0},\quad
DE(\Phi)w=\frac{\dd E(\Phi+t w)}{\dd t}\bigg|_{t=0},
\end{align*}
and, if $\Omega$ is the domain of $\Phi$ (a ball or a half-ball) we set
\begin{align*}
\|D\Area(\Phi)\|
:=\sup\{
|D\Area(\Phi)w|:\, &w\in W^{1,\infty}(\Omega,\R^m),
\,\|w\|_{W^{1,\infty}(\Omega)}\leq 1,
\,\supp w\subset\subset \Omega\}.
\end{align*}
 
\subsection{Lemmas on Conformal Immersions}
\begin{lemma}[Interior Estimates for $\lambda_\Phi$]
\label{lemma:estimates-lambda}
There exists an $\varepsilon_0>0$ such that, if
$\Phi:B_1(0)\to\R^m$ is a conformal immersion 
with $L^2$-bounded second fundamental form satisfying
\begin{align*}
E(\Phi)=\int_{B_1(0)}|\vec{\II}_\Phi|^2_{g_\Phi}\,\dd vol_\Phi
<\varepsilon_0,
\end{align*}
then, 
\begin{enumerate}[(i)]
\item
for any $0<r<1$ there holds
\begin{align}
\label{eq:estimate-dirichlet-energy-lambda}
\int_{B_r(0)}|\nabla\lambda_\Phi|^2\,\dd x
\leq 
\bigg(
\frac{r^2}{2}+C\varepsilon_0
\bigg)
\int_{B_1(0)}|\nabla^2\Phi|^2_{g_\Phi}\,\dd  vol_\Phi,
\end{align}

\item for any compact set $K\subset\subset {B_1(0)}$, there holds
\begin{align}
\label{eq:sup-estimate-lambda}
\|\lambda_\Phi-(\lambda_\Phi)_{B_1(0)}\|_{L^\infty(K)}
\leq\frac{C}{\dist(K,\partial {B_1(0)})^2}\|\nabla\lambda_\Phi\|_{L^{(2,\infty)}({B_1(0)})}
+C\varepsilon_0,
\end{align}
where $(\lambda_\Phi)_{B_1(0)}=\dashint_{B_1(0)}\lambda_\Phi\,\dd x$
denotes the average of $\lambda_\Phi$ on ${B_1(0)}$
and $C>0$ is a constant independent of $\Phi$.
\end{enumerate}
\end{lemma}

\begin{bfproof}[Proof of Lemma \ref{lemma:estimates-lambda}]
Thanks to lemma \ref{lemma:Helein}, if $\varepsilon_0$ is small enough
we may find a Coulomb orthonormal frame $\vec{f}$ defined on ${B_1(0)}$
satisfying  the estimate
\begin{align*}
\|\nabla\vec{f}\|^2_{L^2({B_1(0)})}
\leq C\|\nabla\vec{n}_\Phi\|^2_{L^2({B_1(0)})}
\leq C\varepsilon_0.
\end{align*}
We write $\lambda_\Phi=\mu+h$, where $\mu$
is the solution be the solution to
\begin{align*}
\left\{
\begin{aligned}
-\Delta\mu&=\langle\nabla^\perp\vec{f}_1,\nabla\vec{f}_2\rangle &&\text{in }{B_1(0)},\\
\mu&=0 &&\text{on }\partial {B_1(0)},
\end{aligned}
\right.
\end{align*}
and $h$ is the harmonic rest.
By Wente's lemma, $\mu$ belongs to $C^0(\overline{B_r(0)})\cap W^{1,2}(B_r(0))$
with
\begin{align}
\label{eq:Wente-estimates-mu}
\|\nabla\mu\|_{L^2({B_1(0)})}+\|\mu\|_{L^\infty({B_1(0)})}&\leq
C\|\nabla \vec{n}_\Phi\|^2_{L^2({B_1(0)})}
\leq C\varepsilon_0.
\end{align}
As for $h$, since it is harmonic, for any $0<r<1$ it satisfies 
\begin{align}
\label{eq:inequality-dirichlet-energy-harmonic-rest}
\int_{B_{r}(0)}|\nabla h|^2\,\dd x
\leq r^2\int_{{B_1(0)}}|\nabla h|^2\,\dd x.
\end{align}
Using successively \eqref{eq:inequality-dirichlet-energy-harmonic-rest},
the Dirichlet principle, estimate \eqref{eq:Wente-estimates-mu}
 and identity \eqref{eq:pointwise-Hessian-identity} we then deduce
\begin{align*}
\int_{B_{r}(0)}|\nabla\lambda_\Phi|^2\,\dd x
&=\int_{B_{r}(0)}|\nabla(\mu+h)|^2\,\dd x\\
&\leq 2r^2\int_{{B_1(0)}}|\nabla h|^2 + 2\int_{{B_1(0)}}|\nabla\mu|^2\,\dd x\\
&\leq 2r^2\int_{{B_1(0)}}|\nabla\lambda_\Phi|^2\,\dd x
+C\|\nabla \vec{n}_\Phi\|^4_{L^2({B_1(0)})}\\
&\leq \frac{1}{2}r^2\int_{{B_1(0)}}|\nabla^2\Phi|^2_{g_\Phi}\,\dd  vol_\Phi
+C\|\nabla \vec{n}_\Phi\|^4_{L^2({B_1(0)})}\\
&=\bigg(\frac{r^2}{2}+C\varepsilon_0\bigg)
\bigg(\int_{{B_1(0)}}|\nabla^2\Phi|^2_{g_\Phi}\,\dd  vol_\Phi\bigg),
\end{align*}
which proves \eqref{eq:estimate-dirichlet-energy-lambda}.
As for \eqref{eq:sup-estimate-lambda}, we note that $h$ may be written
by means of the Poisson kernel
\footnote{
the explicit formula is  (see \cite[\textsection 2.2.4]{MR2597943}):
\begin{align*}
K(x,y) = \frac{1-|x|^2}{2\pi }\frac{1}{|x-y|^2}.
\end{align*}}
 as
\begin{align*}
h(x)=\int_{\partial {B_1(0)}}
K(x,y)\,\lambda(y)\,\dd\mathcal{H}^1(y),
\quad x\in {B_1(0)},
\end{align*}
consequently
we deduce that, for any $x\in K$,
using the trace theorem
and Poincar\'e's inequality,
there holds
\begin{align*}
|h(x)-(\lambda_\Phi)_{{B_1(0)}}|
&\leq\frac{1-|x|^2}{2\pi}
\int_{\partial {B_1(0)}}\frac{1}{|x-y|^2}
|\lambda_\Phi(y)-(\lambda_\Phi)_{B_1(0)}|\,\dd\mathcal{H}^1(y)\\
&\leq \frac{C}{\dist(K,\partial B_1(0))^2}\|\lambda_\Phi-(\lambda_\Phi)_{B_1(0)}\|_{L^1(\partial {B_1(0)})}\\
&\leq \frac{C}{\dist(K,\partial B_1(0))^2}\|\nabla\lambda_\Phi\|_{L^1({B_1(0)})}\\
&\leq \frac{C}{\dist(K,\partial B_1(0))^2}\|\nabla\lambda_\Phi\|_{L^{(2,\infty)}({B_1(0)})}.
\end{align*}
We may then conclude with \eqref{eq:Wente-estimates-mu} that
\begin{align*}
\|\lambda_\Phi-(\lambda_\Phi)_{B_1(0)}\|_{L^\infty(K)}
&\leq\|h-(\lambda_\Phi)_{B_1(0)}\|_{L^\infty(K)}+ \|\mu\|_{L^\infty({B_1(0)})}\\
&\leq \frac{C}{\dist(K,\partial B_1(0))^2}\|\nabla\lambda_\Phi\|_{L^{(2,\infty)}({B_1(0)})}
+ C\varepsilon_0,
\end{align*}
as desired.
This concludes the proof of the lemma.
\end{bfproof}

\begin{lemma}[Boundary Estimates for $\lambda_\Phi$]
\label{lemma:estimates-lambda-bndry}
There exists an $\varepsilon_0>0$ so that the following holds.
Let $\Phi:B_1^+(0)\to\R^m$ be a conformal immersion with 
$L^2$-bounded second fundamental form and let $\vec{e}=(\vec{e}_1,\vec{e}_2)$
be its coordinate frame.
If for some $p>1$ we have $\partial_\tau\vec{e}\in L^p(I,\R^m\times\R^m)$
and
\begin{align}
\label{eq:assumption-smallness-bndry}
\int _{B_1^+(0)}|\II_\Phi|^2_{g_\Phi}\,\dd vol_\Phi
+[\vec{e}\,]_{W^{1/2,2}(I)}^2<\varepsilon_0,
\end{align}
then:
\begin{enumerate}[(i)]
\item for any $0<r<1$ there holds
\begin{align}
\label{eq:estimate-dirichlet-bndry}
&
\int_{B_r^+(0)}|\nabla\lambda_\Phi|^2\,\dd x
\leq 
\bigg(\frac{r^2}{2}\bigg)
\int_{B_1^+(0)}|\nabla^2\Phi|^2_{g_\Phi}\,\dd vol_\Phi\\
&\notag
\phantom{{}+{}}+\left(
C\varepsilon_0 + C(p)\|\partial_\tau\vec{e}_1\|_{L^p(I)}
(1+\|\partial_\tau\vec{e}_2\|_{L^p(I)})
\right)
\bigg(\int_{B_1^+(0)}|\vec{\II}_\Phi|^2_{g_\Phi}\,\dd vol_\Phi
+\|\langle\partial_\tau\vec{e}_1,\vec{e}_2\rangle\|_{L^1(I)}\bigg),
\end{align}

\item for any compact set $K\subset \overline{B_1^+(0)}$
so that $\dist(K, S)>0$ there holds
\begin{align}
\label{eq:estimate-sup-norm-bndry}
\inf_{c\in\R}\|\lambda_\Phi -c\|_{L^\infty(K)}
\leq \frac{C}{\dist(K,S)^2}\|\nabla\lambda_\Phi \|_{L^{2,\infty}(B_1^+(0))}
+C(p)(\varepsilon_0 + \|\partial_\tau\vec{e}_1\|_{L^p(I)}
(1+\|\partial_\tau\vec{e}_2\|_{L^p(I)})),
\end{align}
where $C>0$ is an universal constant and 
$C(p)>0$ is a constant depending only on $p$.
\end{enumerate}
\end{lemma}

\begin{bfproof}[Proof of Lemma \ref{lemma:estimates-lambda-bndry}]
\emph{Step 1.} Thanks to lemma \ref{lemma:extension-in-S^1} (see also remark \ref{remark:extension-in-S^1}),
we can consider   an extension $\vec{f}$ of $\vec{e}|_{I}$ to
all of $\partial B_1^+(0)$ such  that 
\begin{align*}
[\vec{f}\,]_{W^{1/2,2}(\partial B_1^+(0))}
\leq 2 [\vec{f}\,]_{W^{1/2,2}(I)},
\end{align*}

\emph{Step 2.}  We choose $\varepsilon_0$ sufficiently small so that thanks to lemma
\ref{lemma:frame-prescribed-boundary}
we may find a frame $\vec{g}=(\vec{g}_1,\vec{g}_2)$ lifting $\vec{n}_\Phi$
on $B_1^+(0)$ that coincides with $\vec{f}$ on $\partial B_1^+(0)$
and satisfies
\begin{align*}
\|\nabla\vec{g}\|_{L^2(B_1^+(0))}
&\leq C\left(
\|\nabla\vec{n}_\Phi\|_{L^2(B_1^+(0))} + [\vec{f}]_{W^{1/2,2}(\partial B_1^+(0))}
\right) \leq C\varepsilon_0.
\end{align*}

\emph{Step 3.} We write $\lambda_\Phi = \mu + h$, where:
\begin{align*}
&\left\{
\begin{aligned}
-\Delta \mu &= \langle\nabla^\perp \vec{g}_1,\nabla \vec{g}_1\rangle &&\text{in } B_1^+(0),\\
\partial_\nu \mu&=\langle\partial_\tau \vec{g}_1,\vec{g}_2\rangle &&\text{on }I,\\
\mu&=0 &&\text{on }S,
\end{aligned}
\right.
&\text{and}&
&
\left\{
\begin{aligned}
-\Delta h &= 0 &&\text{in } B_1^+(0),\\
\partial_\nu h&=0 &&\text{on }I,\\
h&=\lambda &&\text{on }S.
\end{aligned}
\right.
\end{align*}

\emph{Step 4:  estimate of $\mu$.}
Since $\mu$ satisfies a homogeneous Dirichlet condition on $S$,
its extension to $\R^2_+$  by means of odd inversion along $S$
(given by the conformal map $x\mapsto x/|x|^2$):
\begin{align*}
\widehat{\mu}(x)=
\begin{cases}
\mu(x) &\text{for }x\in B_1^+(0),\\
-\mu(x/|x|^2) &\text{for }x\in \R^2_+\setminus B_1^+(0),
\end{cases}
\end{align*}
satisties (also thanks to the transformation law under conformal maps of the Laplace operator and
of the determinant):
\begin{align}
\label{eq:extended-Neumann-inversion}
\left\{
\begin{aligned}
-\Delta \widehat{\mu} 
&= \langle\nabla^\perp \widetilde{\vec{g}}_1,\nabla \widetilde{\vec{g}}_2\rangle &&\text{in } \R^2_+,\\
\partial_\nu \widehat{\mu}
&=\langle\partial_\tau \widetilde{\vec{g}}_1,\widetilde{\vec{g}}_2\rangle &&\text{on }\partial \R^2_+,
\end{aligned}
\right.
\end{align}
where for $i=1,2$,
\begin{align*}
\widetilde{\vec{g}}_i(x)=
\begin{cases}
\vec{g}_i(x) &\text{for }x\in B_1^+(0),\\
\vec{g}_i(x/|x|^2) &\text{for }x\in \R^2_+\setminus B_1^+(0),
\end{cases}
\end{align*}
denote extensions of $\vec{g}_i$ by means of even inversion along $S$.\par
We then consider the Cayley map
$\phi(z) = i\frac{-z+1}{z+1}$ mapping biholomorphically $B_1(0)$ onto $\R^2_+$ and by simplicity of notation  we continue to denote by
$\widehat{\mu}$ and $\widetilde{\vec{g}}_i$ the composition $\widehat{\mu}\circ\phi$ and $\widetilde{\vec{g}}_i\circ\phi.$
 By the conformal invariance of the problem \eqref{eq:extended-Neumann-inversion},  $\widehat{\mu}\circ\phi$ satisfies  
 \begin{align}
\label{eq:extended-Neumann-inversion}
\left\{
\begin{aligned}
-\Delta \widehat{\mu} 
&= \langle\nabla^\perp \widetilde{\vec{g}}_1,\nabla \widetilde{\vec{g}}_2\rangle &&\text{in } B_1(0),\\
\partial_\nu \widehat{\mu}
&=\langle\partial_\tau \widetilde{\vec{g}}_1,\widetilde{\vec{g}}_2\rangle &&\text{on }\partial B_1(0),
\end{aligned}
\right.
\end{align}
 Thanks to lemma  \ref{lemma:estimate-Wente-Neumann}
we have
\begin{align*}
\inf_{c\in\R}\|{\widehat{\mu}}-c\|_{L^\infty(B_1(0))}
&\leq C\|\nabla{\widetilde{\vec{g}}_1}\|_{L^2(B_1(0))}
(1+\|\nabla{\widetilde{\vec{g}}_2}\|_{L^2(B_1(0))})\\
&\phantom{{}\leq{}}+C(p)\|\partial_\tau {\widetilde{\vec{g}}_1}\|_{L^p(\partial B_1(0))}
(1+\|\partial_\tau {\widetilde{\vec{g}}_2}\|_{L^p(\partial B_1(0))}).
\end{align*}
Because of the conformal invariance of the Dirichlet energy there holds
\begin{align*}
\|\nabla{\widetilde{\vec{g}}_i}\|_{L^2(B_1(0))}
=\|\nabla\vec{g}_i\|_{L^2(B_1^+(0))}
\quad i=1,2,
\end{align*}
on the other hand a direct computation shows that
\begin{align*}
\|\partial_\tau {\widetilde{\vec{g}}_i}\|_{L^p(\partial B_1(0))}
\leq 2\|\partial_\tau \vec{g}_i\|_{L^p(I)}
\quad i=1,2,
\end{align*}
hence we deduce (assuming without loss of generality that $\varepsilon_0<1$
and recalling that $\vec{g}=\vec{e}$ on $I$):
\begin{align}
\label{eq:estimate-sup-mu-bndry}
\inf_{c\in\R}\|\mu-c\|_{L^\infty(B_1^+(0))}
\leq C\varepsilon_0 
+ C(p)\|\partial_\tau\vec{e}_1\|_{L^p(I)}
(1+\|\partial_\tau\vec{e}_2\|_{L^p(I)}).
\end{align}
Consequently, through integration by parts and by using \eqref{eq:estimate-sup-mu-bndry} we can estimate
\begin{align}
\label{eq:estimate-mu-bndry}
\int_{B_1^+(0)}|\nabla\mu|^2\,\dd x
&=-\int_{B_1^+(0)}(\mu-c)\Delta u\,\dd x
+\int_I (\mu-c)\partial_\nu\lambda\,\dd\mathcal{H}^1 \\
&\leq\inf\|\mu-c\|_{L^\infty(B_1^+(0))}
\bigg(
\|\langle \nabla^\perp\vec{g}_1,\nabla\vec{g}_2\rangle\|_{L^1(B_1(0))}
+\|\langle\partial_\tau\vec{e}_1,\vec{e}_2\rangle\|_{L^1(I)}
\bigg)\nonumber\\
&\leq( C\varepsilon_0 
+ C(p)\|\partial_\tau\vec{e}_1\|_{L^p(I)}
(1+\|\partial_\tau\vec{e}_2\|_{L^p(I)}))
\bigg(
\int_{B_1^+(0)}|\vec{\II}_\Phi|^2\,\dd vol_\Phi
+ \|\langle\partial_\tau\vec{e}_1,\vec{e}_2\rangle\|_{L^1(I)}
\bigg).\nonumber
\end{align}

\emph{Step 5: estimate of the harmonic rest $h$.}
We observe that the existence of $h$ can be deduced by variational methods. Since $h$ satisfies a homogeneous Neumann condition along $I$,
its  extension to $B_1(0)$ by even reflection along $I$:
\begin{align*}
\widetilde{h}(x)=
\begin{cases}
h(x^1,x^2) & \text{in } B_1^+(0),\\
h(x^1,-x^2)  & \text{in } B_1^-(0) = B_1(0)\cap \R^2_-,
\end{cases}
\end{align*}
will then satisfy
\begin{align*}
\left\{
\begin{aligned}
-\Delta \widetilde{h} &= 0 &&\text{in } B_1(0),\\
\tilde h&=\widetilde{\lambda}_\Phi &&\text{on }\partial B_1(0),
\end{aligned}
\right.
\end{align*}
where $\widetilde{\lambda}_\Phi$  similarly denotes  the
extension  of $\lambda_\Phi$ to $B_1(0)$ by even reflection along $I$.
From the classical estimate for harmonic function
we will then deduce
\begin{align*}
\int_{B_r(0)}|\nabla\widetilde{h}|^2\,\dd x
\leq r^2\int_{B_1(0)}|\nabla\widetilde{\lambda}_\Phi|^2\,\dd x,
\end{align*}
and consequently,
\begin{align}
\label{eq:estimate-h-bndry}
\int_{B_r^+(0)}|\nabla h|^2\,\dd x
\leq r^2\int_{B_1^+(0)}|\nabla \lambda_\Phi|^2\,\dd x.
\end{align}
By joining estimates \eqref{eq:estimate-mu-bndry}-\eqref{eq:estimate-h-bndry}
we then deduce
\begin{align*}
\int_{B_r^+(0)}|\nabla\lambda_\Phi|^2\,\dd x
&\leq 2\int_{B_1^+(0)}|\nabla\mu|^2\,\dd x 
+2\int_{B_r^+(0)}|\nabla h|^2\,\dd x\\
&\leq 
(C\varepsilon_0 
+ C(p)\|\partial_\tau\vec{e}_1\|_{L^p(I)}
(1+\|\partial_\tau\vec{e}_2\|_{L^p(I)}))
\bigg(
\int_{B_1^+(0)}|\vec{\II}_\Phi|^2\,\dd vol_\Phi
+ \|\langle\partial_\tau\vec{e}_1,\vec{e}_2\rangle\|_{L^1(I)}
\bigg)\\
&\phantom{{}\leq{}}
+r^2\int_{B_1^+(0)}|\nabla \lambda_\Phi|^2\,\dd x,
\end{align*}
which then yields estimate \eqref{eq:estimate-dirichlet-bndry}.\par

As far as the estimate
\eqref{eq:estimate-sup-norm-bndry} is concerned,
similarly as in lemma \ref{lemma:estimates-lambda},
we deduce that
\begin{align*}
\|\widetilde{h} - (\widetilde{\lambda}_\Phi)_{B_1(0)}\|_{L^\infty (\widetilde{K})}
\leq\frac{C}{\dist(\widetilde{K}, \partial_1B(0))}
\|\nabla\widetilde{\lambda}_\Phi\|_{L^{(2,\infty)}(B_1(0))},
\end{align*}
where we denoted by $\widetilde{K}=K \cup\{(x^1,-x^2): (x^1,x^2)\in K\}$,
and consequently that
\begin{align}
\label{eq:estimate-sup-h-bndry}
\|h - (\lambda_\Phi)_{B^+_1(0)}\|_{L^\infty (K)}
\leq\frac{C}{\dist(K, \partial_1B(0))^2}
\|\nabla\lambda_\Phi\|_{L^{(2,\infty)}(B_1^+(0))}.
\end{align}
We may then write
\begin{align*}
\inf_{c'\in\R}\|\lambda_\Phi -c'\|_{L^\infty(K)}
\leq
\|h - (\lambda_\Phi)_{B_1(0)}\|_{L^\infty (K)}
+\inf_{c\in \R}\|\mu-c\|_{L^\infty(B_1^+(0))},
\end{align*}
and with estimates \eqref{eq:estimate-sup-mu-bndry}-\eqref{eq:estimate-sup-h-bndry}
we deduce the validity of \eqref{eq:estimate-sup-norm-bndry}.
This concludes the proof of the lemma.
\end{bfproof}

\begin{lemma}[Affine approximation]
\label{lemma:conformal-affine-immersion}
For every $\delta>0$
there exists an $\varepsilon_0>0$ such that,
for every conformal immersion $\Phi:{B_1(0)}\to\R^m$
with $L^2$-bounded second fundamental form 
satisfying 
\begin{align*}
\int_{{B_1(0)}}|\nabla^2\Phi|^2_{g_\Phi}\,\dd  vol_\Phi<\varepsilon_0,
\end{align*}
there exists a conformal affine immersion
$L=\vec{L}_0+x^1\vec{X}_1+x^2\vec{X}_2
=\vec{L}_0+\langle x,\vec{X}\rangle$
so that
\begin{align}
\label{eq:W22-affine-approx}
\|\Phi - L\|_{W^{2,2}(B_{1/2}(0))}
<\delta\|\nabla\Phi\|_{L^2({B_1(0)})}
\end{align}
and, if $\e^\nu = |\vec{X}_1|=|\vec{X}_2|$ denotes the conformal
factor of $L$,
\begin{align}
\label{eq:logconfact-affine-approx}
\|\lambda_\Phi - \nu\|_{L^\infty(B_{1/2}(0))}
<\delta.
\end{align}
\end{lemma}

\begin{bfproof}[Proof of Lemma \ref{lemma:conformal-affine-immersion}]
We argue by contradiction and suppose that there exists
a $\delta>0$ such that, for every $k\in\N$, there is a conformal
Lipschitz immersion with $L^2$-bounded second fundamental form
$\Phi_k:{B_1(0)}\to\R^m$ such  that (writing as usual $\e^{\lambda_{\Phi_k}}=\e^{\lambda_k}$),
\begin{align}
\label{eq:absurd-hypotesis-1}
\int_{{B_1(0)}}|\nabla^2\Phi_k|^2_{g_{\Phi_k}}\dd  vol_{\Phi_k}
=\int_{{B_1(0)}}\e^{-2\lambda_k}|\nabla^2\Phi_k|^2\dd x\leq \frac{1}{k},
\end{align}
and  for every conformal affine immersion $L$ there holds
\begin{align}
\label{eq:contradiction-assumption}
\|\Phi_k - L\|_{W^{2,2}(B_{1/2}(0))}
>\delta\|\nabla\Phi_k\|_{L^2({B_1(0)})}.
\end{align}
or, if $\e^\nu$ denotes the conformal factor of $L$,
\begin{align}
\label{eq:contradiction-assumption-2}
\|\lambda_k - \nu\|_{L^\infty(B_{1/2}(0))}
\geq\delta.
\end{align}
Since \eqref{eq:absurd-hypotesis-1} is invariant under 
translations and dilations in $\R^m$, writing for short
$
c_k=-\dashint_{B_1(0)}\lambda_k\,\dd x 
$
if we set
\begin{align}
\label{eq:widetilde-Phi}
\widetilde{\Phi}_k(x) = \e^{c_k}(\Phi_k(x)-\Phi_k(0)),
\quad x\in {B_1(0)},
\end{align}
then $\widetilde{\Phi}_k:{B_1(0)}\to\R^m$ defines for every $k$ a conformal
Lipschitz immersion with $L^2$-bounded second fundamental form
such  that, if $\e^{\widetilde{\lambda}_k}$ denotes its conformal factor,
there holds
\begin{align*}
\widetilde{\lambda}_k&=\lambda_k+c_k,&
\widetilde{\Phi}_k(0)&=0,&
\dashint_{B_1(0)}\widetilde{\lambda}_k\,\dd x&=0,
\end{align*}
and
\begin{align}
\label{eq:absurd-hypotesis-1-tilde}
\int_{B_1(0)}\e^{-2\widetilde{\lambda}_k}
|\nabla^2\widetilde{\Phi}_k|^2\,\dd x\leq\frac{1}{k}.
\end{align}
From the identites \eqref{eq:pointwise-Hessian-identity}, Lemma \ref{lemma:estimates-lambda}
  and  
\eqref{eq:absurd-hypotesis-1-tilde} it follows
\begin{align*}
\|\widetilde{\lambda}_k
-(\widetilde\lambda_k)_{B_1(0)}\|_{L^\infty(B_{1/2}(0))}
&=\|\widetilde{\lambda}_k\|_{L^\infty(B_{1/2}(0))}\\
&\leq C\|\nabla\widetilde{\lambda}_k\|_{L^2({B_1(0)})}
+C\|\nabla\vec{\widetilde{n}}_k\|_{L^2({B_1(0)})}^2\leq\frac{C}{k}
\to 0\quad\text{as }k\to\infty,
\end{align*}
so
\begin{align}
\label{eq:uniorm-convergence-lambda-tilde}
\widetilde{\lambda}_k\to 0\quad\text{uniformly in }B_{1/2}(0),
\end{align}
and consequently we infer that
\begin{align}
\label{eq:convergence-gradient-Phi-tilde}
\|\nabla\widetilde{\Phi}_k\|_{L^2(B_{1/2}(0))}
&= \int_{B_{1/2}(0)}2\e^{2\widetilde{\lambda}_k}\,\dd x\to\pi/2
\quad\text{as }k\to\infty,\\
\|\nabla^2\widetilde{\Phi}_k\|_{L^2(B_{1/2}(0))}
&\leq\e^{C/k} \int_{B_{1/2}(0)}
\e^{-2\widetilde{\lambda}_k}|\nabla^2\widetilde{\Phi}_k|^2\,\dd x
\leq \frac{C}{k}\to 0\quad\text{as }k\to\infty,
\end{align}
from which we deduce, with Poincar\'e's inequality,
\begin{align*}
\|\Phi_k\|_{W^{2,2}(B_{1/2}(0))}
=\|\Phi_k-\Phi_k(0)\|_{W^{2,2}(B_{1/2}(0))}\leq C.
\end{align*}
We deduce that, up to extraction of subsequences,
for a map $\widetilde{\Phi}_\infty\in W^{2,2}(B_{1/2}(0),\R^m)$ we have
\begin{align*}
\widetilde{\Phi}_k&\rightharpoonup\widetilde{\Phi}_{\infty}
\quad\text{in }W^{2,2}(B_{1/2}(0),\R^m),\\
\nabla^2\widetilde{\Phi}_k&\to 0
\quad\text{in }L^{2}(B_{1/2}(0),\R^m),\\
\widetilde{\Phi}_k&\to\widetilde{\Phi}_{\infty}
\quad\text{in }W^{1,p}(B_{1/2}(0),\R^m)\text{ for every }1\leq p<\infty
\text{ and a.e. in }B_{1/2}(0).
\end{align*}
consequently, $\Phi_\infty:B_{1/2}(0)\to\R^m$ is a conformal map
and, from the uniform convergence of $\widetilde{\lambda}_k$
above, up to a further subsequence, its conformal factor is 1
(that is, $\tilde{\Phi}_\infty$ is a isometric linear immersion).
Being $\nabla^2\Phi_\infty = 0$,
there actually holds
\begin{align}
\label{eq:absurd-W-2-2-convergence}
\widetilde{\Phi}_k\to \widetilde{\Phi}_{\infty}
\quad\text{in }W^{2,2}(B_{1/2}(0),\R^m).
\end{align}
Note now that from the definition \eqref{eq:widetilde-Phi}, \eqref{eq:contradiction-assumption} is equivalent to
\begin{align*}
\|\e^{-c_k}\widetilde{\Phi}_k+\Phi_k(0)-L\|_{W^{2,2}(B_{1/2}(0))}
\geq \delta\|\nabla\Phi_k\|_{L^2(B_{1/2}(0))},
\end{align*}
hence
\begin{align*}
\|\widetilde{\Phi}_k+\e^{c_k}(\Phi_k(0)-L)\|_{W^{2,2}(B_{1/2}(0))}
\geq \delta\|\nabla(\e^{c_k}\Phi_k)\|_{L^2(B_{1/2}(0))}
=\delta\|\nabla(\widetilde{\Phi}_k)\|_{L^2(B_{1/2}(0))},
\end{align*}
Since $L$ is arbitrary,
we may consider the sequence $L_k=\Phi_k(0)+\e^{-c_k}\widetilde{\Phi}_\infty$
and deduce from \eqref{eq:convergence-gradient-Phi-tilde} that
\begin{align*}
\|\widetilde{\Phi}_k-\widetilde{\Phi}_\infty\|_{W^{2,2}(B_{1/2}(0))}
&\geq \delta \|\nabla(\widetilde{\Phi}_k)\|_{L^2(B_{1/2}(0))}\\
&=\delta(\pi/2+o(1))\quad \text{as }k\to\infty,
\end{align*}
which is a contradiction with \eqref{eq:absurd-W-2-2-convergence}.
Similarly, if \eqref{eq:contradiction-assumption-2} holds,
since the conformal factor of $L_k$ is  $\e^{-c_k}$,
we have
\begin{align*}
\|\lambda_k-\nu_k\|_{L^\infty(B_{1/2}(0))}
= \|\lambda_k+c_k\|_{L^\infty(B_{1/2}(0))}
= \|\widetilde\lambda_k\|_{L^\infty(B_{1/2}(0))}
\geq\delta
\quad\text{for every }k\in \N,
\end{align*}
which contradicts \eqref{eq:uniorm-convergence-lambda-tilde}.
\end{bfproof}

\subsection{Construction of Suitable Competitors}
\begin{lemma}[Interior Competitors]
\label{lemma:biharmonic-comparison}
Let $\Phi:{B_1(0)}\to\R^m$ be a
conformal immersion with $L^2$-bounded
second fundamental form.
For every $\delta>0$ there
there exists an $\varepsilon_0>0$ such that if 
\begin{align*}
\int_{B_{4r}(0)}
|\nabla^2\Phi|^2_{g_\Phi}\,\dd  vol_\Phi<\varepsilon_0,
\end{align*}
 for some $0<r\leq1/4$,
then there exists a $\rho\in[r/2,r]$ 
such  that the solution to
\begin{align}
\label{eq:biharmonic-extension}
\left\{
\begin{aligned}
\Delta^2 \psi &=0 &&\text{in }B_\rho(0),\\
\psi&= \Phi &&\text{on }\partial B_\rho(0),\\
\nabla \psi &=\nabla \Phi &&\text{on }\partial B_\rho(0),
\end{aligned}
\right.
\end{align}
defines an immersion which satisfies:
\begin{align}
\label{eq:biharmonic-estimate}
\int_{B_\rho(0)}|\nabla^2\psi|^2_{g_\psi}\,\dd vol_\psi
\leq C
\big(
1+C_0(\delta+o(\delta))
\big)
\int_{B_r(0)\setminus B_{r/2}(0)}|\nabla^2\Phi|^2_{g_\Phi}\,\dd  vol_\Phi,
\end{align}
and
\begin{align}
\label{eq:estimte-area-Phi-psi}
|\Area(\Phi|_{B_\rho(0)})-\Area(\psi)|
\leq
C_0(\delta+o(\delta))
\|\nabla\Phi\|_{L^2(B_{\rho}(0))}^2,
\end{align}
and
\begin{align}
\label{eq:estimte-D-area-Phi-psi}
&\phantom{{}={}}
\|D\Area(\Phi|_{B_\rho(0)})-D\Area(\psi)\|
\leq 
C_0(\delta+o(\delta))
\|\nabla\Phi\|_{L^2(B_{r}(0))},
\end{align}
where $C>0$ independent of $r$ and $\Phi$,
$C_0>0$ depends only on
$\|\nabla\lambda_\Phi\|_{L^{(2,\infty)}(B_{1}(0))}$ and 
$o(\delta)/\delta\to 0$ as $\delta\to 0$.
\end{lemma}

\begin{remark}
We note that:
\begin{enumerate}[(i)]
\item For every $\rho$, the existence and uniqueness of a solution to
\eqref{eq:biharmonic-extension} in $W^{2,2}$ is given, for example, by the fact that
such problem is the Euler-Lagrange equation for the biharmonic energy functional
$\mathcal{B}(\sigma)=\int_{B_\rho(0)}|\Delta\sigma|^2\,\dd\leb^2$
(or, equivalently, of the Hessian energy $\int_{B_\rho(0)}|\nabla^2\sigma|^2\,\dd\leb^2$),
subject to the prescribed boundary data.
Since $\Phi$ is of class $W^{2,2}$, existence and uniqueness by an
argument similar to the one for the Dirichlet problem.
\item An elementary fact that will be used in the proof of lemma \ref{lemma:biharmonic-comparison} is the following.
For a function $f\in L^1(B_R(0))$ ($R>0$ is arbitrary) and a constant $C>0$
we say that $\rho\in [R/2,R]$ defines a
\emph{$C$-good slice for $f$} (in $B_{R}(0)\setminus B_{R/2}(0))$
if $f|_{\partial B_\rho(0)}$ is in $L^1(\partial B_\rho(0))$
and there holds
\begin{align*}
\rho\int_{\partial B_\rho(0)}|f(r,\theta)|\,d\sigma
\leq C \int_{B_R(0)\setminus B_{R/2}(0)}|f|\,\dd x.
\end{align*}
The existence of $C$-good slices for some $C$
and any $f$ 
is a consequence of Fubini's theorem.
Moreover, one can check that,
for every $0<\delta<R/2$, there exists a $C_\delta>0$
so that, for every $f\in L^1(B_R(0))$, the radii $\rho\in[R/2,R]$
defining $C_\delta$- good slices for $f$ have
Lebesgue measure at least $R/2-\delta$.
\end{enumerate}
\end{remark}

\begin{bfproof}[Proof of Lemma \ref{lemma:biharmonic-comparison}]
It is sufficient to prove the thesis only for any
$\delta>0$ sufficiently small.
We first treat the case $r=1/4$ and
argue through rescalings at the end.
In what follows, we denote by $C$ a positive constant 
(possibly varying line to line) which is independent of $\Phi$,
and with $C_0$ a positive constant depending only on
$\|\nabla\lambda_1\|_{L^{(2,\infty)}(B_1(0))}$.

\emph{Step 1.}
For $\varepsilon_0$ sufficiently small as in lemma \ref{lemma:estimates-lambda},
we have that
\begin{align}
\label{eq:sup-estimate-lambda-3/4}
\|\lambda_\Phi-(\lambda_\Phi)_{B_1(0)}\|_{L^\infty(B_{3/4}(0))}
\leq C\left(
\|\nabla\lambda_\Phi\|_{L^{(2,\infty)}(B_1(0))} + \varepsilon_0
\right)=C_0,
\end{align}
where $(\lambda_\Phi)_{B_1(0)}$ denotes the average of $\lambda_\Phi$ over $B_1(0)$.
Also, for every $\delta>0$, if $\varepsilon_0$ is sufficiently small as in lemma \ref{lemma:conformal-affine-immersion},
then there exists a conformal affine immersion
$L$
whose conformal factor we denote by
$\e^\nu$,
satisfying the following estimate
\begin{align}
\label{eq:W22-affine-approx-1/4}
\|\Phi - L\|_{W^{2,2}(B_{1/4}(0))}
&<\delta\|\nabla\Phi\|_{L^2({B_{1/2}(0)})},\\
\label{eq:logconfact-affine-approx-1/4}
\|\lambda_\Phi - \nu\|_{L^\infty(B_{1/4}(0))}.
&<\delta.
\end{align}
By combining \eqref{eq:sup-estimate-lambda-3/4}-\eqref{eq:logconfact-affine-approx-1/4},
we deduce 
\begin{align}\label{estconffactor}
\|\lambda_\Phi-\nu\|_{L^\infty(B_{3/4}(0))}
&\leq\|\lambda_\Phi-(\lambda_\Phi)_{B_1(0)}\|_{L^\infty(B_{3/4}(0))}
+|(\lambda_\Phi)_{B_1(0)}-\nu|\\
&\leq \|\lambda_\Phi-(\lambda_\Phi)_{B_1(0)}\|_{L^\infty(B_{3/4}(0))}
+\|\lambda_\Phi-(\lambda_\Phi)_{B_1(0)}\|_{L^\infty(B_{1/2}(0))}
+\|\lambda_\Phi-\nu\|_{L^\infty(B_{1/2}(0))}\nonumber\\
&\leq C_0+\delta.\nonumber
\end{align}
It follows that
\begin{align}
\label{eq:poinwise-estimate-confact-3/4}
C_0^{-1}(1-\delta-o(\delta))\e^\nu
\leq\e^{\lambda_\Phi(x)}\leq
C_0(1+\delta+o(\delta))\e^\nu
\quad\text{for }x\in B_{3/4}(0),
\end{align}
and in particular
\begin{align}
\label{eq:estimate-area-3/4}
C_0^{-1}(1-\delta-o(\delta))\e^{2\nu}
\leq\|\nabla\Phi\|_{L^2(B_{3/4}(0))}^2\leq
C_0(1+\delta+o(\delta))\e^{2\nu}.
\end{align}
We then consider a good-slice choice $\rho\in[1/8,1/4]$  so that
$\Phi$ and $\Phi -L$ belong to $W^{2,2}(\partial B_\rho(0),\R^m)$ with
\begin{align*}
\|\Phi\|_{W^{2,2}(\partial B_\rho(0))}&\leq C\|\Phi\|_{W^{2,2}(B_{1/4}(0)\setminus B_{1/8}(0))},\\
\|\Phi-L\|_{W^{2,2}(\partial B_\rho(0))}
&\leq C\|\Phi-L\|_{W^{2,2}(B_{1/4}(0)\setminus B_{1/8}(0))},
\end{align*}
hence we consider the solution to \eqref{eq:biharmonic-extension}
for such choice of $\rho$.
Elliptic regularity theory
(see for instance \cite[Remark 7.2, Chapter 2]{MR0350177})
implies that
\begin{align*}
\|\psi-L\|_{W^{5/2,2}(B_\rho(0))}&\leq C\|\Phi-L\|_{W^{2,2}(\partial B_\rho(0))},
\end{align*}
while Sobolev embedding $W^{5/2,2}\hookrightarrow C^{1,\alpha}$ implies that, for every $0<\alpha<1/2$,
\begin{align*}
\|\psi-L\|_{C^{1,\alpha}(\overline{B_\rho(0)})}\leq C\|\psi-L\|_{W^{5/2,2}(B_\rho(0))}.
\end{align*}
Hence we have
\begin{align*}
\|\nabla\psi-\nabla L\|_{L^\infty(B_\rho(0))}
&\leq C\|\psi-L\|_{W^{5/2,2}(B_\rho(0))} &&\text{(Sobolev embedding)}\\
\notag&\leq C\|\Phi-L\|_{W^{2,2}(\partial B_\rho(0))} &&\text{(elliptic estimates)}\\
\notag&\leq C\|\Phi-L\|_{W^{2,2}(B_{1/4}(0)\setminus B_{1/8}(0))}&&\text{(good-slice choice)}\\
&\leq C\delta\|\nabla\Phi\|_{L^2(B_{1/4}(0))}
&&\text{(by \eqref{eq:W22-affine-approx-1/4})}\\
&\leq C_0\e^\nu(\delta+o(\delta))
&&\text{(by \eqref{eq:poinwise-estimate-confact-3/4} and \eqref{eq:estimate-area-3/4})}.
\end{align*}
Hence for $i=1,2$, we deduce the pointwise estimates in $B_\rho(0)$ 
\begin{align*}
||\partial_i\psi|^2-\e^{2\nu}|
&=||\partial_i\psi|-\e^{\nu}|\,||\partial_i\psi|+\e^\nu|\\
&\leq|\partial_i\psi-\partial_i L|\,||\partial_i\psi|+\e^\nu|\\
&\leq C_0\e^{2\nu}
(\delta+o(\delta)),
\end{align*}
and similarly
\begin{align*}
|\langle\partial_1\psi,\partial_2\psi\rangle|
&=|\langle\partial_1\psi,\partial_2\psi\rangle 
- \langle\partial_1L,\partial_2 L\rangle|\\
&=|\langle
(\partial_1\psi-\partial_1 L),\partial_2\psi\rangle
+\langle\partial_1L,(\partial_2\psi-\partial_2 L)\rangle|\\
&\leq ||\partial_2\psi|+|\partial_1L||
\big(|\partial_1\psi-\partial_1L|+|\partial_2\psi-\partial_2L|
\big)\\
&\leq C_0\e^{2\nu}(\delta+o(\delta)).
\end{align*}
This implies that, if $g_\psi = (\langle\partial_i\psi,\partial_j\psi\rangle)_{ij}$
denotes the metric associated with $\psi$,
for every vector $X=(X^1,X^2)\in\R^2$, we have
\begin{align}
\label{eq:metric-delta-nu}
\e^{2\nu}\left(1-C_0(\delta+ o(\delta))\right)|X|^2
\leq g_\psi(X,X)\leq 
\e^{2\nu}\left(1+C_0(\delta+o(\delta))\right)|X|^2,
\end{align}
where $|X|=\sqrt{(X^1)^2+(X^2)^2}$ is the Euclidean norm of $X$.
As a consequence, we deduce that
for $\delta>0$ small enough, $\psi$ defines an immersion,
and in such case we have
\begin{align}
\label{eq:determinant-delta-nu}
\e^{2\nu}\left(1-C_0(\delta+o(\delta))\right)
\leq \sqrt{\det g_\psi}
\leq \e^{2\nu}\left(1+C_0(\delta+o(\delta))\right),
\end{align}
and
\begin{align*}
\e^{-4\nu}\left(
1+C_0(\delta+o(\delta))\right)^{-2}
|\nabla^2\psi|^2
\leq |\nabla^2\psi|^2_{g_\psi}
\leq \e^{-4\nu}\left(
1-C_0(\delta+o(\delta))\right)^{-2}
|\nabla^2\psi|^2
\end{align*}
consequently wee that, point-wise in $B_\rho(0)$,
\begin{align}
\label{eq:estimate-II-psi-1}
|\nabla^2\psi|^2_{g_\psi}\sqrt{\det g_\psi}
\leq \frac{1+C_0(\delta+o(\delta))}
{\left(
1-C_0(\delta+o(\delta))\right)^{2}}\e^{-2\nu}
|\nabla^2\psi|^2.
\end{align}

\emph{Step 2: estimate for the curvature energy.}
Since $\psi$ solves \eqref{eq:biharmonic-extension},
from elliptic regularity theory, we have that 
for any affine function $M(x)=\vec{M}_0+\langle\vec{Y},x\rangle$
there holds
\begin{align*}
\|\nabla^2\psi\|_{L^2(B_\rho(0))}
&=\|\nabla^2(\psi-M)\|_{L^2(B_\rho(0))}\\
&\leq C\left(\|\Phi-M\|_{W^{2,2}(\partial B_\rho(0))} 
+ \|\nabla(\Phi-M)\|_{W^{1,2}(\partial B_\rho(0))}\right)\\
&\leq C\left(
\|\Phi-M\|_{L^2(\partial B_\rho(0))} 
+ \|\nabla(\Phi-M)\|_{L^2(\partial B_\rho(0))}
+ \|\nabla^2\Phi\|_{L^2(\partial B_\rho(0))}\right),
\end{align*}
and hence if we suitably choose $M$ so that
\begin{align*}
\|\nabla(\Phi-M)\|_{L^2(\partial B_\rho(0))}&\leq C\|\nabla^2\Phi\|_{L^2(\partial B_\rho(0))},\\
\|\Phi-M\|_{L^2(\partial B_\rho(0))}&\leq C\|\nabla(\Phi-M)\|_{L^2(\partial B_\rho(0))}
\end{align*}
(if $M(x)=\vec{M}_0+\langle Y,x\rangle$, it is sufficient to choose
$Y=(\nabla\Phi)_{\partial B_\rho(0)}$ and $M_0=(\Phi-M)_{\partial B_\rho(0)}$),
we actually deduce that
\begin{align*}
\|\nabla^2\psi\|_{L^2(B_\rho(0))}\leq C\|\nabla^2\Phi\|_{L^2(\partial B_\rho(0))},
\end{align*}
and so, from the choice of $\rho$ we made,
we have
\begin{align*}
\|\nabla^2\psi\|_{L^2(B_\rho(0))}\leq C\|\nabla^2\Phi\|_{L^2(\partial B_\rho(0))}
\leq C\|\nabla^2\Phi\|_{L^2(B_{1/4}(0)\setminus B_{1/8}(0))}.
\end{align*}
Note also that	
\begin{align}
\label{eq:estimate-II-psi-3}
\int_{B_{1/4}(0)\setminus B_{1/8}(0)}|\nabla^2\Phi|^2\,\dd x
&\leq\int_{B_{1/4}(0)\setminus B_{1/8}(0)}\e^{2\lambda_\Phi}\e^{-2\lambda_\Phi}|\nabla^2\Phi|^2\,\dd x\\
\notag&\leq\e^{2\delta}\e^{2\nu}
\int_{B_{1/4}(0)\setminus B_{1/8}(0)}\e^{-2\lambda_\Phi}|\nabla^2\Phi|^2\,\dd x\\
\notag&\leq \e^{2\nu}(1+2(\delta+o(\delta)))
\int_{B_{1/4}(0)\setminus B_{1/8}(0)}\e^{-2\lambda_\Phi}|\nabla^2\Phi|^2\,\dd x,
\end{align}
By joining estimates
\eqref{eq:estimate-II-psi-1}
--
\eqref{eq:estimate-II-psi-3},
 we deduce that
\begin{align*}
\int_{B_\rho(0)}|\nabla^2\psi|_{g_\psi}\,\dd vol_{g_\psi}
&\leq C
(1+C_0(\delta +o(\delta)))
\int_{B_{1/4}(0)\setminus B_{1/8}(0)}|\nabla^2\Phi|^2_{g_\Phi}\,\dd  vol_\Phi.
\end{align*}

\emph{Step 3: estimate on the area.}
From \eqref{eq:logconfact-affine-approx-1/4}
and
\eqref{eq:determinant-delta-nu}\footnote{
\eqref{eq:logconfact-affine-approx-1/4} implies
\begin{align*}
(1-2(\delta+o(\delta)))\e^{2\nu}
\leq \e^{2\lambda_\Phi}
\leq (1+2(\delta+o(\delta)))\e^{2\nu}
\quad\text{in }B_{1/4}(0).
\end{align*}
}
we deduce that
(recall that $\rho\in[1/8,1/4]$),
we have
\begin{align*}
|\sqrt{\det g_\psi} - \e^{2\lambda_\Phi}|
&\leq \e^{2\nu}C_0(\delta+o(\delta))
\leq \e^{2\lambda_\Phi}C_0(\delta+o(\delta))
\quad\text{in }B_{\rho}(0),
\end{align*}
hence by integrating over $B_\rho(0)$ we deduce
\begin{align*}
|\Area(\psi)-\Area(\Phi|_{B_\rho(0)})|
&=\bigg|\int_{B_\rho(0)}
\left(\sqrt{\det g_\psi}
-\e^{2\lambda_\Phi}\right)\dd x\bigg|\\
&\leq \int_{B_\rho(0)}
\left|\sqrt{\det g_\psi}
-\e^{2\lambda_\Phi}\right|\dd x\\
&\leq C_0\|\nabla\Phi\|_{L^2(B_\rho(0))}^2(\delta +o(\delta)).
\end{align*}

\emph{Step 4: estimate on the derivative of the area.}
For $w\in W^{1,\infty}(B_\rho(0),\R^m)$ with 
compact support in $B_\rho(0)$,
we have
\begin{align*}
D\Area(\Phi)w 
&=\int_{B_\rho(0)}\langle\nabla\Phi,\nabla w\rangle\,\dd x,\\
D\Area(\psi)w 
&=\int_{B_\rho(0)}g_\psi(\nabla\psi,\nabla w)
\dd vol_\psi
=\int_{B_\rho(0)}g_\psi^{ij}
\langle\partial_i\psi,\partial_jw\rangle
\sqrt{\det g_\psi}\,\dd x,
\end{align*}
and from \eqref{eq:metric-delta-nu} and \eqref{eq:determinant-delta-nu} we deduce that
\begin{align}\label{estdetg}
\delta^{ij}(1-C_0(\delta+o(\delta)))
\leq
g^{ij}_\psi
\sqrt{\det g_\psi}
&\leq \delta^{ij}(1+C_0(\delta+o(\delta))).
\end{align}
We also observe that
\begin{equation}\label{estpsiminuspsi}
\|\nabla\psi-\nabla \Phi\|_{L^2(B_{\rho}(0))}\le C\delta \|\nabla\Phi\|_{L^2(B_{1/4}(0))}.
\end{equation}
Moreover
\begin{align*}
g^{ij}_\psi \langle\partial_i\psi,\partial_jw\rangle
\sqrt{\det g_\psi}-\delta^{ij}\langle\partial_i\Phi,\partial_jw\rangle
&=g^{ij}_\psi \langle\partial_i\psi,\partial_jw\rangle
\sqrt{\det g_\psi}-g^{ij}_\psi \langle\partial_i\Phi,\partial_jw\rangle
\sqrt{\det g_\psi}\\&+
g^{ij}_\psi \langle\partial_i\Phi,\partial_jw\rangle
\sqrt{\det g_\psi}-\delta^{ij}\langle\partial_i\Phi,\partial_jw\rangle.
\end{align*}
By using \eqref{estdetg} we get that
\begin{align}\label{estprov1}
\int_{B_{\rho}(0)}\left| g^{ij}_\psi \langle\partial_i\Phi,\partial_jw\rangle
\sqrt{\det g_\psi}-\delta^{ij}\langle\partial_i\Phi,\partial_jw\rangle\right| \dd x 
&\leq C_0(\delta+o(\delta)) \int_{B_{\rho}(0)} |\langle\nabla\Phi,\nabla w\rangle |\dd x 
\end{align}
and
\begin{align}\label{estprov2}
&\int_{B_{\rho}(0)}\left|
g^{ij}_\psi \langle\partial_i\psi,\partial_jw\rangle
\sqrt{\det g_\psi}-g^{ij}_\psi \langle\partial_i\Phi,\partial_jw\rangle
\sqrt{\det g_\psi}\right| \dd x\\
&\le C_0\|\nabla \psi-\nabla\Phi\|_{L^2({B_{\rho}(0)})}\le C_0\delta \|\nabla\Phi\|_{L^2(B_{1/4}(0))}.\nonumber
\end{align}

By combining   estimates \eqref{estprov1} and \eqref{estprov2} we get
\begin{align*}
|D\Area(\psi)w-D\Area(\Phi)w|
&\leq C_0(\delta+o(\delta)) \int_{B_{\rho}(0)} |\langle\partial_i\Phi,\partial_i w\rangle | \dd x+
C_0\|\nabla \psi-\nabla\Phi\|_{L^2({B_{\rho}(0)})}\\
&\le C_0\delta \|\nabla\Phi\|_{L^2(B_{1/4}(0))}
\end{align*}

\emph{Step 5: the estimates for a general $r$.}
If $0<r<1/4$,
we may reduce to the case $r=1/4$:
indeed, if 
we consider the rescaling $\widetilde{\Phi}(x)=\Phi(4rx)$
for $x\in {B_1(0)}$,  by conformal invariance we 
have
\begin{align*}
\int_{B_{4r}(0)}|\nabla^2\Phi|^2_{g_\Phi}\,\dd  vol_\Phi
&=\int_{{B_1(0)}}|\nabla^2\Phi|^2_{g_\Phi}\,\dd  vol_\Phi,
\end{align*}
and the area functional, the Dirichlet energy,
solution to the problem \eqref{eq:biharmonic-extension}
and the $L^{(2,\infty)}$-seminorm of $\nabla\lambda_\Phi$
are invariant by rescalings as well
\footnote{to see this last fact,
note that, if $d_{\nabla \lambda_\Phi}(t)=\
\leb^2\left(\{x:|\nabla\lambda_\Phi(x)|>t\}\right\})$
denotes the distribution function of $\nabla\lambda_\Phi$, and $\sigma>0$,
${\lambda_\Phi}_{\Phi,\sigma}(x)=\lambda_\Phi(\sigma x)$
has distribution function 
$d_{\nabla \lambda_{\Phi,\sigma}}(t) = \sigma^{-2}d_{\nabla\lambda_\Phi}(t/\sigma)$.
Consequently
\begin{align*}
\|\nabla\lambda_{\Phi,\sigma}\|_{L^{(2,\infty)}}
=\sup_{t>0}t\sqrt{d_{\nabla\lambda_{\Phi,\sigma}}(t)}
=\sup_{u>0}u\sqrt{d_{\nabla\lambda_\Phi}(u)}
=\|\nabla\lambda_\Phi\|_{L^{(2,\infty)}}.
\end{align*}
}.
We may apply the previous steps to $\widetilde{\Phi}$,
estimate $\|\nabla\lambda_\Phi\|_{L^{(2,\infty)}(B_{4r}(0))}$
with $\|\nabla\lambda_\Phi\|_{L^{(2,\infty)}(B_1(0))}$ 
and then rescale back.
\end{bfproof}

\begin{definition}
\label{def:geometric-reflection-planar-boundary}
An immersion  $\Phi:B_1^+(0)\to\R^m$ 
is said to have \emph{flat geometric boundary data}
on the base diameter $I$ if
there holds
\begin{align*}
\Phi(x^1,0)\in
\spanop_\R\{\vec{\varepsilon}_1\}
\quad\text{and}\quad
\vec{n}_\Phi(x^1,0) = \vec{\varepsilon}_3\wedge\ldots\wedge\vec{\varepsilon}_{m-2}
\quad\text{for }(x^1,0)\in I.
\end{align*}
For a conformal immersion $\Phi:B_1^+(0)\to\R^m$ with flat geometric boundary data on $I$,
its \emph{geometric reflection} along $I,$
$\widehat{\Phi}:B_1(0)\to\R^3$  is defined as
\begin{align*}
\widehat\Phi(x^1,x^2)=
\begin{cases}
\Phi(x^1,x^2)&\text{if }x^2\geq 0,\\
\Phi^1(x^1,-x^2)\vec{\varepsilon}_1 - \sum_{l=2}^m\Phi^l(x^1,-x^2)\vec{\varepsilon}_l
&\text{if }x^2< 0.\\
\end{cases}
\end{align*}
\end{definition}
\medskip
Note that
if $\Phi$ as in definition \ref{def:geometric-reflection-planar-boundary} is conformal and   with $L^2$-bounded second fundamental form, there holds
\begin{align*}
\partial_1\Phi(x^1,0)
=\e^{\lambda(x^1,0)}\vec{\varepsilon}_1,
\quad\text{and}\quad
\partial_2\Phi(x^1,0)
=\e^{\lambda(x^1,0)}\vec{\varepsilon}_2.
\end{align*}
Hence, provided that
$\|\lambda_\Phi\|_{L^\infty(B_1^+(0))}<+\infty$,
the geometric nature of the reflection  and conformality imply that
$\widehat{\Phi}$ defines a conformal immersion of class
$(W^{1,\infty}\cap W^{2,2})(B_1(0),\R^3)$,
hence a Lipschitz immersion with $L^2$-bounded second fundamental form with
\begin{align*}
|\nabla\widehat\Phi(x^1,x^2)|^2&=
|\nabla\Phi(x^1,x^2)|^2\chi_{\{x^2\geq 0\}}
+|\nabla\Phi(x^1,-x^2)|^2\chi_{\{x^2<0\}},\\
|\nabla^2\widehat\Phi(x^1,x^2)|^2&=
|\nabla^2\Phi(x^1,x^2)|^2\chi_{\{x^2\geq 0\}}
+|\nabla^2\Phi(x^1,-x^2)|^2\chi_{\{x^2<0\}},\\
|\Delta\widehat\Phi(x^1,x^2)|^2&=
|\Delta\Phi(x^1,x^2)|^2\chi_{\{x^2\geq 0\}}
+|\Delta\Phi(x^1,-x^2)|^2\chi_{\{x^2<0\}},\\
\e^{\lambda_{\widehat\Phi}(x^1,x^2)}&=
\e^{\lambda_{\Phi}(x^1,x^2)}
\chi_{\{x^2\geq 0\}}
+\e^{\lambda_{\Phi}(x^1,-x^2)}
\chi_{\{x^2<0\}}.
\end{align*}
The following is a boundary analogue of lemma \ref{lemma:biharmonic-comparison},
where additionally we have flat geometric boundary data on $I$ in the sense of definition
\ref{def:geometric-reflection-planar-boundary}. \begin{lemma}
\label{lemma:biharmonic-comparison-bndry}
There exists an $\varepsilon_0$ with the following property.
Let $\Phi:{B_1^+(0)}\to\R^m$ be a conformal immersion with $L^2$-bounded
second fundamental form  and   flat geometric boundary data on $I$  such 
 that $\|\lambda_\Phi\|_{L^\infty(B_1^+(0))}<+\infty$.
For every $\delta>0$ there
there exists an $\varepsilon_0>0$ such that,
if \begin{align*}
\int_{B_{4r}^+(0)}
|\nabla^2\Phi|^2_{g_\Phi}\,\dd vol_\Phi
<\varepsilon_0,
\end{align*}
 for  $0<r\leq1/4$, 
then there exists a $\rho\in[r/2,r]$ 
and an immersion $\psi\in C^{1,\alpha}(\overline{B_\rho^+(0)},\R^m)$
which satisfies
\begin{align*}
\psi=\Phi\quad\text{on }\partial B_{\rho}(0)\cap B_1^+(0),\\
\nabla\psi=\nabla\Phi\quad\text{on }\partial B_{\rho}(0)\cap B_1^+(0),
\end{align*}
has flat geometric boundary data on $\rho I$
and satisfies
\begin{align}
\label{eq:biharmonic-estimate-bndry}
\int_{B^+_\rho(0)}|\nabla^2\psi|^2_{g_{\psi}}\,\dd vol_{\psi}
&\leq C
\big(
1+C_0(\delta+o(\delta))
\big)
\int_{B^+_r(0)\setminus B^+_{r/2}(0)}|\nabla^2\Phi|^2_{g_\Phi}\,\dd  vol_\Phi\\
&\notag
\phantom{{}\leq{}}
+C_0(\delta+o(\delta))
\Area(\Phi|_{B_\rho^+(0))}),
\end{align}
and
\begin{align}
\label{eq:estimate-area-Phi-psi-bndry}
|\Area(\Phi|_{B_\rho^+(0)})-\Area(\psi)|
\leq
C_0(\delta+o(\delta))
\Area(\Phi|_{B_\rho^+(0))}),
\end{align}
and
\begin{align}
\label{eq:estimate-D-area-Phi-psi-bndry}
&\phantom{{}={}}
\|D\Area(\Phi|_{B^+_\rho(0)})-D\Area(\psi)\|
\leq 
C_0(\delta+o(\delta))
\Area(\Phi|_{B_r^+(0))}),
\end{align}
where $C>0$ is independent of $r$ and $\Phi$,
$C_0>0$ may depend  on
$\|\nabla\lambda_\Phi\|_{L^{(2,\infty)}(B_{1}^+(0))}$ and 
$o(\delta)\to 0$ as $\delta\to 0$.
\end{lemma}

Parts of the proof of this lemma are similar to the proof
of lemma \ref{lemma:biharmonic-comparison}, so we will focus on the differences.
The overall idea is first to reflect geometrically $\Phi$ as in definition
\ref{def:geometric-reflection-planar-boundary}, then consider the biharmonic competitor
as in lemma \ref{lemma:biharmonic-comparison} and finally to smoothly ``correct'' it so that
it has flat geometric boundary data on $\rho I$.
Such ``correction'' will be essentially constructed
by means of the 1st order Taylor polynomials of such biharmonic comparison
at the points $(\rho,0)$ and $(-\rho,0)$ respectively.

\begin{bfproof}[Proof of Lemma \ref{lemma:biharmonic-comparison-bndry}]
As in the case of lemma \ref{lemma:biharmonic-comparison},
it sufficient to prove the thesis only for any
$\delta>0$ sufficiently small.
We first treat the case $r=1/4$ and
argue through rescalings at the end.
In what follows, we denote by $C$ a positive constant 
which is independent of $\Phi$,
and with $C_0$ a positive constant depending only on
$\|\nabla\lambda_\Phi\|_{L^{(2,\infty)}(B_1^+(0))}$.

\emph{Step 1.}
We consider the geometric reflection of $\Phi$, $\widehat{\Phi}:B_1(0)\to\R^m$,
according to definition \ref{def:geometric-reflection-planar-boundary}.
For $\varepsilon_0$ sufficiently small as in lemma \ref{lemma:estimates-lambda},
we have that
\begin{align}
\label{eq:sup-estimate-lambda-3/4-bndry}
\|\lambda_{\widehat\Phi}-(\lambda_{\widehat\Phi})_{B_1(0)}\|_{L^\infty(B_{3/4}(0))}
&\leq C\left(
\|\nabla\lambda_{\widehat\Phi}\|_{L^{(2,\infty)}(B_1(0))} + \varepsilon_0
\right)\\
&\notag
\leq C\left(
\|\nabla\lambda_\Phi\|_{L^{(2,\infty)}(B_1^+(0))} + \varepsilon_0
\right)\leq C_0,
\end{align}
where $(\lambda_{\widehat\Phi})_{B_1(0)}$ denotes the average of $\lambda_{\widehat\Phi}$ over $B_1(0)$.
Also, for every $\delta>0$, if $\varepsilon_0$ is sufficiently small as in lemma \ref{lemma:conformal-affine-immersion},
then there exists a conformal affine immersion $L$ whose conformal factor
$\e^\nu$ is so that the estimates
\begin{align}
\label{eq:W22-affine-approx-1/4-bndry}
\|\Phi - L\|_{W^{2,2}(B^+_{1/4}(0))}
&<\delta\|\nabla\Phi\|_{L^2({B^+_{1/2}(0)})}
& &\Longleftrightarrow &
\|\widehat\Phi - L\|_{W^{2,2}(B_{1/4}(0))}
&<\delta\|\nabla\widehat\Phi\|_{L^2({B_{1/2}(0)})}\\
\label{eq:logconfact-affine-approx-1/4-bndry}
\|\lambda_\Phi - \nu\|_{L^\infty(B_{1/4}^+(0))}
&<\delta
& &\Longleftrightarrow &
\|\lambda_{\widehat\Phi} - \nu\|_{L^\infty(B_{1/4}(0))}&<\delta
\end{align}
are satisfied.
By combining \eqref{eq:sup-estimate-lambda-3/4-bndry}
-\eqref{eq:logconfact-affine-approx-1/4-bndry},
we deduce 
\begin{align*}
\|\lambda_{\widehat\Phi}-\nu\|_{L^\infty(B_{3/4}(0))}
&\leq\|\lambda_{\widehat\Phi}-(\lambda_{\widehat\Phi})_{B_1(0)}\|_{L^\infty(B_{3/4}(0))}
+|(\lambda_{\widehat\Phi})_{B_1(0)}-\nu|\\
&\leq \|\lambda_{\widehat\Phi}-(\lambda_{\widehat\Phi})_{B_1(0)}\|_{L^\infty(B_{3/4}(0))}
+\|\lambda_{\widehat\Phi}-(\lambda_{\widehat\Phi})_{B_1(0)}\|_{L^\infty(B_{1/2}(0))}
+\|\lambda_{\widehat\Phi}-\nu\|_{L^\infty(B_{1/2}(0))}\\
&\leq C_0+\delta,
\end{align*}
consequently we point-wise estimate from above and below
\begin{align}
\label{eq:poinwise-estimate-confact-3/4-bndry}
C_0^{-1}(1-\delta-o(\delta))\e^\nu
\leq\e^{\lambda_{\widehat\Phi}(x)}\leq
C_0(1+\delta+o(\delta))\e^\nu
\quad\text{for }x\in B_{3/4}(0).
\end{align}
We then consider a good-slice choice $\rho\in[1/8,1/4]$  so that
$\widehat\Phi$ and $\widehat\Phi -L$ belong to $W^{2,2}(\partial B_\rho(0),\R^m)$
(equivalently, so that
$\Phi$ and $\Phi -L$ belong to $W^{2,2}(\partial B_\rho(0)\cap B_1^+(0),\R^m)$)
with
\begin{align}
\label{goodslicehatphi}
\|\widehat\Phi\|_{W^{2,2}(\partial B_\rho(0))}&\leq C\|\widehat\Phi\|_{W^{2,2}(B_{1/4}(0)\setminus B_{1/8}(0))},\\
\label{goodslicehatphi2}
\|\widehat\Phi-L\|_{W^{2,2}(\partial B_\rho(0))}
&\leq C\|\widehat\Phi-L\|_{W^{2,2}(B_{1/4}(0)\setminus B_{1/8}(0))},
\end{align}
hence we consider the solution
for such choice of $\rho$ to
\begin{align*}
\left\{
\begin{aligned}
\Delta^2 \psi_0 &=0 &&\text{in }B_\rho(0),\\
\psi_0&= \widehat\Phi &&\text{on }\partial B_\rho(0),\\
\nabla \psi_0 &=\nabla \widehat\Phi &&\text{on }\partial B_\rho(0),
\end{aligned}
\right.
\end{align*}
which satisfies, as in lemma \ref{lemma:biharmonic-comparison}, the estimates
\begin{align}
\label{eq:elliptic-est-1-bndry}
\|\psi_0-L\|_{W^{5/2,2}(B_\rho(0))}&\leq C_0\e^\nu(\delta+o(\delta)),\\
\|\nabla^2\psi\|_{L^2(B_\rho(0))}
&\label{eq:elliptic-est-3-bndry}
\leq C\|\nabla^2\widehat{\Phi}\|_{L^2(B_{1/4}(0)\setminus B_{1/8}(0))}
\leq C\|\nabla^2\Phi\|_{L^2(B_{1/4}(0)\setminus B_{1/8}(0))},
\end{align}
and consequently, by Sobolev embedding $W^{5/2,2}\hookrightarrow C^{1,\alpha}$
for every $0<\alpha<1/2$ we have,
\begin{align}
\label{eq:elliptic-est-2-bndry}
\|\nabla^2\psi_0\|_{C^{1,\alpha}(\overline{B_\rho(0)})}
&\leq C_0\e^\nu(\delta+o(\delta)).
\end{align}

The estimates \eqref{eq:elliptic-est-1-bndry} and \eqref{eq:elliptic-est-2-bndry}  will be crucial for what follows.
If $T^1_{\psi_0,(-\rho,0)}(x)$ and $T^1_{\psi_0,(\rho,0)}(x)$
denote the Taylor polynomial of $\psi_0$ at the points $(-\rho,0)$
$(\rho,0)$ respectively,
we may write, for $x\in B_\rho(0)$,
\begin{align*}
T^1_{\psi_0,(-\rho,0)}(x)
&=T^1_{\Phi,(-\rho,0)}(x)
:=
\Phi(-\rho,0) 
+\left\langle
\nabla\Phi(-\rho,0),
\begin{pmatrix}
x^1+\rho\\
x^2
\end{pmatrix}
\right\rangle,\\
T^1_{\psi_0,(\rho,0)}(x)
&=T^1_{\Phi,(\rho,0)}(x)
:=
\Phi(\rho,0) 
+\left\langle
\nabla\Phi(\rho,0),
\begin{pmatrix}
x^1-\rho\\
x^2
\end{pmatrix}
\right\rangle,
\end{align*}
where the expressions on the right-hand sides have a well-defined meaning
because of our good-slice choice of $\rho$.
Note moreover that $T^1_{\psi_0,(-\rho,0)}(x)$ and $T^1_{\psi_0,(\rho,0)}(x)$
are conformal, they define a parametrization of the plane $\spanop\{\vec{\varepsilon}_1,\vec{\varepsilon}_2\}$
and are so that $T^1_{\psi_0,(-\rho,0)}(x^1,0), T^1_{\psi_0,(\rho,0)}(x^1,0)\in
\spanop\{\vec{\varepsilon}_1\}$ (in particular they have flat geometric boundary data
on $\rho I$ according
to definition \ref{def:geometric-reflection-planar-boundary}).
For every $x\in B_\rho(0)$, 
by virtue of Taylor's theorem there exists $\xi\in ((-\rho, 0),x)\subset B_{\rho}(0)$
so that
\begin{align*}
\psi_0(x) - T^1_{\psi_0,(-\rho,0)}(x)
=\left\langle
(\nabla\psi_0(\xi)-\nabla\psi_0(-\rho,0)),
\begin{pmatrix}
x^1+\rho\\
x^2
\end{pmatrix}
\right\rangle,
\end{align*}
consequently we deduce
that for every $x\in B_\rho(0)$ there holds
\begin{align}
\label{eq:estimate-phi0-1st}
|\psi_0(x)-T^1_{\psi_0, (-\rho,0)}(x)|
&\leq[\nabla\psi_0]_{C^{0,\alpha}(\overline{B_\rho(0)})}
\left|
\begin{pmatrix}
x^1+\rho\\
x^2
\end{pmatrix}\right|^{1+\alpha},\\
\label{eq:estimate-phi0-2nd}
|\nabla\psi_0(x)-\nabla T^1_{\psi_0, (-\rho,0)}|
&\leq[\nabla\psi_0]_{C^{0,\alpha}(\overline{B_\rho(0)})}
\left|
\begin{pmatrix}
x^1+\rho\\
x^2
\end{pmatrix}\right|^{\alpha},
\end{align}
and similarly that
\begin{align}
\label{eq:estimate-phi0-3rd}
|\psi_0(x)-T^1_{\psi_0, (\rho,0)}(x)|
&\leq[\nabla\psi_0]_{C^{0,\alpha}(\overline{B_\rho(0)})}
\left|
\begin{pmatrix}
x^1-\rho\\
x^2
\end{pmatrix}\right|^{1+\alpha},\\
|\nabla\psi_0(x)-\nabla T^1_{\psi_0, (\rho,0)}|
&\label{eq:estimate-phi0-4th}
\leq[\nabla\psi_0]_{C^{1,\alpha}(\overline{B_\rho(0)})}
\left|
\begin{pmatrix}
x^1-\rho\\
x^2
\end{pmatrix}\right|^{\alpha}.
\end{align}
Note also that for every $x\in B_\rho(0)$  we may estimate
\begin{align*}
|T^1_{\Phi,(-\rho,0)}(x)-L(x)|
&=\left|T^1_{\Phi,(-\rho,0)}(x)-L(-\rho,0) - 
\left\langle
\nabla L(-\rho,0),
\begin{pmatrix}
x^1+\rho\\
x^2
\end{pmatrix}\right\rangle
\right|\\
&\leq |\Phi(-\rho,0) - L(-\rho,0)|
+|\nabla\Phi(-\rho,0)-\nabla L(-\rho,0)|
\left|
\begin{pmatrix}
x^1+\rho\\
x^2
\end{pmatrix}\right|,\\
|\nabla T^1_{\Phi, (-\rho,0)} - \nabla L|
&=|\nabla\Phi(-\rho,0)-\nabla L(-\rho,0)|,
\end{align*}
and similarly
\begin{align*}
|T^1_{\Phi,(\rho,0)}(x)-L(x)|
&\leq |\Phi(\rho,0) - L(\rho,0)|
+|\nabla\Phi(\rho,0)-\nabla L(\rho,0)|
\left|
\begin{pmatrix}
x^1-\rho\\
x^2
\end{pmatrix}\right|,\\
|\nabla T^1_{\Phi, (\rho,0)} - \nabla L(\rho,0)|
&=|\nabla\Phi(\rho,0)-\nabla L(\rho,0)|,
\end{align*}
hence thanks to \eqref{eq:W22-affine-approx-1/4-bndry}-\eqref{goodslicehatphi2},
\begin{align}
\label{eq:estimate-taylor-L-bndry-1}
\|T^1_{\Phi,(-\rho,0)}-L\|_{C^1(\overline{B_\rho(0)})}
&\leq 
C_0\e^{\nu}(\delta+o(\delta)),\\
\label{eq:estimate-taylor-L-bndry-2}
\|T^1_{\Phi,(\rho,0)}-L\|_{C^1(\overline{B_\rho(0)})}
&\leq 
C_0\e^{\nu}(\delta+o(\delta)).
\end{align}
Consequently, with \eqref{eq:elliptic-est-2-bndry} we deduce
\begin{align*}
\|\psi_0 - T^1_{\Phi,(-\rho,0)}\|_{C^{1,\alpha}(\overline{B_\rho(0)})}
&\leq \|\psi_0-L\|_{C^{1,\alpha}(\overline{B_\rho(0)})}
+\|T^1_{\Phi,(-\rho,0)} -L\|_{C^{1,\alpha}(\overline{B_\rho(0)})}\\
&\leq C_0\e^{\nu}(\delta+o(\delta)),
\end{align*}
and similarly
\begin{align*}
\|\psi_0 - T^1_{\Phi,(\rho,0)}\|_{C^{1,\alpha}(\overline{B_\rho(0)})}
&\leq C_0\e^{\nu}(\delta+o(\delta)).
\end{align*}

\emph{Step 2.}
We now let $f:\R\to\R$ be a non-negative smooth function so that
\begin{align*}
f(t)=
\begin{cases}
0 &\text{for }t\leq -\rho/2,\\
1 &\text{for }t\geq \rho/2,
\end{cases}
\end{align*}
and we set, for $x\in B_\rho(0)$,
\begin{align*}
\phi(x)=T_{\Phi,(-\rho,0)}^1(x)
+f(x^1)(T_{\Phi,(\rho,0)}(x)-T^1_{\Phi,(-\rho,0)}(x)).
\end{align*}
Such function has range in the plane $\spanop\{\vec{\varepsilon}_1,\vec{\varepsilon}_2\}$
and is so that $\phi(x^1,0)\in \spanop\{\vec{\varepsilon}_1\}$.
Moreover, since
\begin{align*}
\partial_1\phi(x)
&=\partial_1T_{\Phi,(-\rho,0)}^1(x) 
+f'(x^1)(T^1_{\Phi,(\rho,0)}(x) - T^1_{\Phi,(-\rho,0)}(x))\\
&\phantom{{}={}}
+f(x^1)(\partial_1T^1_{\Phi,(\rho,0)}(x) - \partial_1 T^1_{\Phi,(-\rho,0)}(x))\\
&=\e^{\lambda_\Phi(-\rho,0)}\vec{\varepsilon}_1
+f'(x^1)(T^1_{\Phi,(\rho,0)}(x) - T^1_{\Phi,(-\rho,0)}(x))\\
&\phantom{{}={}}
+f(x^1)(\e^{\lambda_\Phi(\rho,0)} - \e^{\lambda_\Phi(-\rho,0)})\vec{\varepsilon}_1,\\
&=\big(
\e^{\lambda_\Phi(-\rho,0)}
+f(x^1)(\e^{\lambda_\Phi(\rho,0)} - \e^{\lambda_\Phi(-\rho,0)})\big)\vec{\varepsilon}_1\\
&\phantom{{}={}}
+f'(x^1)(T^1_{\Phi,(\rho,0)}(x) - T^1_{\Phi,(-\rho,0)}(x))\\
\partial_2\phi(x)
&=\big(\e^{\lambda_\Phi(-\rho,0)}
+f(x^1)(\e^{\lambda_\Phi(\rho,0)} - \e^{\lambda_\Phi(-\rho,0)})\big)\vec{\varepsilon}_2,
\end{align*}
we have that, if $\delta>0$ is chosen small enough, it defines an immersion.
Indeed, on the one hand from
\eqref{eq:estimate-taylor-L-bndry-1}-\eqref{eq:estimate-taylor-L-bndry-2}
we can estimate
\begin{align*}
\|T^1_{\Phi,(-\rho,0)} - T^1_{\Phi,(\rho,0)}\|_{C^0(\overline{B_\rho(0)})}
&\leq \|T^1_{\Phi,(-\rho,0)} - L\|_{C^0(\overline{B_\rho(0)})}
+\|T^1_{\Phi,(\rho,0)} - L\|_{C^0(\overline{B_\rho(0)})}\\
&\leq C_0\e^{\nu}(\delta+o(\delta)),
\end{align*}
on the other hand, from  \eqref{eq:logconfact-affine-approx-1/4-bndry} 
and the mean value theorem we may estimate
\begin{align*}
|\e^{\lambda_\Phi(-\rho,0)} - \e^{\lambda_\Phi(\rho,0)}|
&\leq|\lambda_\Phi(-\rho,0)-\lambda_\Phi(\rho,0)| 
\sup\{\e^\xi:\xi\in [\lambda_\Phi(\pm\rho,0),\lambda_\Phi(\mp\rho]\}\\
&\leq 2\delta\sup\{\e^\xi:\xi\in [\nu-\delta,\nu+\delta]\}\\
&\leq C\delta \e^{\nu},
\end{align*}
hence
we have the estimates, uniformly in $x\in B_\rho(0)$,
\begin{align*}
|\partial_1\phi(x)|&\geq
\e^{\lambda_\Phi(-\rho,0)}
-C\delta\e^{\nu} - C_0\|f'\|_{L^\infty((-\rho,\rho))}\e^{\nu}(\delta+o(\delta))\\
&\geq \e^{\nu}\big(\e^{-\delta} - C\delta 
-  C_0\|f'\|_{L^\infty((-\rho,\rho))}(\delta+o(\delta))\big),
\end{align*}
and similarly
\begin{align*}
|\partial_2\phi(x)|&\geq \e^{\nu}\big(\e^{-\delta} - C\delta \big),
\end{align*}
and
\begin{align*}
|\langle\partial_1(x),\partial_2\phi(x)\rangle|
&\leq\e^{2\nu} \big(
\|f'\|_{L^\infty(-\rho,\rho)}C_0(\delta+o(\delta))\big)
(\e^\delta + C\delta).
\end{align*}
These inequalities imply the immersive nature of $\phi$ if
$\delta>0$ is chosen small enough.
Note also that thanks to \eqref{eq:estimate-taylor-L-bndry-1}-\eqref{eq:estimate-taylor-L-bndry-2}
there holds
\begin{align*}
|\nabla^2\phi(x)|
&\leq |f''(x^1)||T^1_{\Phi,(\rho,0)}(x) - T^1_{\Phi,(-\rho,0)}(x)|
+2|f'(x^1)||\nabla T^1_{\Phi,(\rho,0)} - \nabla T^1_{\Phi,(-\rho,0)}|\\
&\leq C_0\e^\nu(\delta+o(\delta)).
\end{align*}
Since we may write
\begin{align*}
\psi_0(x)-\phi(x)
=(1-f(x^1))(T^1_{\Phi,(-\rho,0)}(x)-\psi_0(x))
+f(x^1) (T^1_{\Phi,(\rho,0)}(x)-\psi_0(x)),
\end{align*}
we deduce
thanks to 
\eqref{eq:estimate-phi0-1st}
--
\eqref{eq:estimate-phi0-2nd}
--
\eqref{eq:estimate-phi0-3rd}
--
\eqref{eq:estimate-phi0-4th}
that
\begin{align}
\label{eq:psi-phi}
|\psi_0(x)-\phi(x)|
&\leq [\nabla\psi_0]_{C^{0,\alpha}(\overline{B_\rho(0)})}
\bigg( 
(1-f(x^1))
\left|
\begin{pmatrix}
x^1+\rho\\
x^2
\end{pmatrix}
\right|^{1+\alpha}
+f(x^1)
\left|
\begin{pmatrix}
x^1-\rho\\
x^2
\end{pmatrix}
\right|^{1+\alpha}
\bigg),\\
\label{eq:nabla-psi-phi}
|\nabla\psi_0(x)-\nabla\phi(x)|
&\leq [\nabla\psi_0]_{C^{0,\alpha}(\overline{B_\rho(0)})}
\bigg(
|f'(x^1)|
\left|
\begin{pmatrix}
x^1+\rho\\
x^2
\end{pmatrix}
\right|^{1+\alpha}
+(1-f(x^1))\left|
\begin{pmatrix}
x^1+\rho\\
x^2
\end{pmatrix}
\right|^{\alpha}\\
\notag
&\phantom{[\nabla\psi_0]_{C^{0,\alpha}(\overline{B_\rho(0)})}
\bigg(
|f'(x^1)|\bigg(}
+|f'(x^1)|
\left|
\begin{pmatrix}
x^1-\rho\\
x^2
\end{pmatrix}
\right|^{1+\alpha}
+f(x^1)
\left|
\begin{pmatrix}
x^1-\rho\\
x^2
\end{pmatrix}
\right|^{\alpha}
\bigg),\\
\label{eq:hessian-psi-phi}
|\nabla^2\psi_0(x)-\nabla^2\phi(x)|
&\leq |\nabla^2\psi_0(x)| +C_0\e^\nu(\delta+o(\delta)).
\end{align}

\emph{Step 3.}
In this step we construct a function $\chi:B_\rho(0)\to\R$,
which we will use in a moment, with the following properties:
it is supported in $B_{\rho}(0)\setminus\{(-\rho,0),(\rho,0)\}$,
it is
smooth away from $(-\rho,0),(\rho,0)$ and is so that
\begin{align}
\notag
\chi&\equiv 1
\text{ in a neighbourhood of } (-\rho,\rho)\times\{0\}
\text{ in }B_\rho(0)
\text{ which shrinks at } (\pm\rho,0),\\
\label{eq:nabla-cutoff}
|\nabla\chi(x)|
&\sim
\left|
\begin{pmatrix}
x^1+\rho\\
x^2
\end{pmatrix}
\right|^{-1}
\text{ as }x\to(-\rho,0),
\quad
|\nabla\chi(x)|
\sim 
\left|
\begin{pmatrix}
x^1-\rho\\
x^2
\end{pmatrix}
\right|^{-1}
\text{ as }x\to(\rho,0),\\
\label{eq:hessian-cutoff}
|\nabla^2\chi(x)|
&\sim
\left|
\begin{pmatrix}
x^1+\rho\\
x^2
\end{pmatrix}
\right|^{-2}
\text{ as }x\to(-\rho,0),
\quad
|\nabla^2\chi(x)|
\sim 
\left|
\begin{pmatrix}
x^1-\rho\\
x^2
\end{pmatrix}
\right|^{-2}
\text{ as }x\to(\rho,0).
\end{align}

Such function may be constructed as follows.
Let $k_0:S^1\to\R$ be a smooth function 
so that, for an angle $\beta$ to be specified below, it satisfies
\begin{align*}
k_0(\theta)=
\begin{cases}
1 &\text{if }-\beta\leq\theta\leq\beta,\\
0 &\text{if } \theta\in (-\pi,\pi]\setminus [-\beta,\beta].
\end{cases}
\end{align*}

\begin{align*}
\chi_0(x) = r\,k_0\left(\frac{x}{|x|}\right),
\quad x\in \R^2\setminus\{0\}.
\end{align*}
\begin{figure}
\begin{center}
\begin{tikzpicture}[scale=3, baseline={(0,0)}]
\fill[black!15] (-1,0) -- (-3/4,0) arc [start angle=0, end angle=26.57, radius=1/4] -- cycle;
\draw (-1,0) -- (-3/4,0) arc [start angle=0, end angle=26.57, radius=1/4] -- cycle;
\draw (-4/5,0) node [above right]{$\beta$};
\draw[thick] (0,0) circle [radius=1];
\draw (-1,0) node{$\scriptstyle+$} node[below left]{$-\rho$};
\draw (1,0) node{$\scriptstyle+$} node[below right]{$\rho$};
\draw (-1/2,0) node{$\scriptstyle+$} node[below left]{$-\frac{\rho}{2}$};
\draw (1/2,0) node{$\scriptstyle+$} node[below]{$\frac{\rho}{2}$};
\draw (0,1/4) node{$\scriptstyle+$} node[right]{$\frac{\rho}{4}$};
\draw (0,-1/4) node{$\scriptstyle+$} node[right]{$-\frac{\rho}{4}$};
\draw (0,1/2) node{$\scriptstyle+$} node[right]{$\frac{\rho}{2}$};
\draw (0,-1/2) node{$\scriptstyle+$} node[right]{$-\frac{\rho}{2}$};
\draw[dashed] (-1,0) -- (-1/2,1/2) -- (1/2,1/2);
\draw[dashed] (-1,0) -- (-1/2,1/4) -- (1/2,1/4);
\draw[dashed] (-1,0) -- (-1/2,-1/2) -- (1/2,-1/2);
\draw[dashed] (-1,0) -- (-1/2,-1/4) -- (1/2,-1/4);
\draw[dashed] (-1/2,-3/4) -- (-1/2,3/4);
\draw[thick] (-1,0) -- (-1/2,-1/4);
\draw (-1/2,-1/4) node [below right]{$r$};
\draw (-1,0) circle [radius=sqrt(5)/4];
\draw[->] (-1.2,0)--(1.2,0) node[right]{$x^1$};
\draw[->] (0,-1.2)--(0,1.2) node[above]{$x^2$};
\end{tikzpicture}
\caption{definition of $\alpha$ and $r$ in 
the construction of $\chi$}
\label{fig:construction-cutoff}
\end{center}
\end{figure}
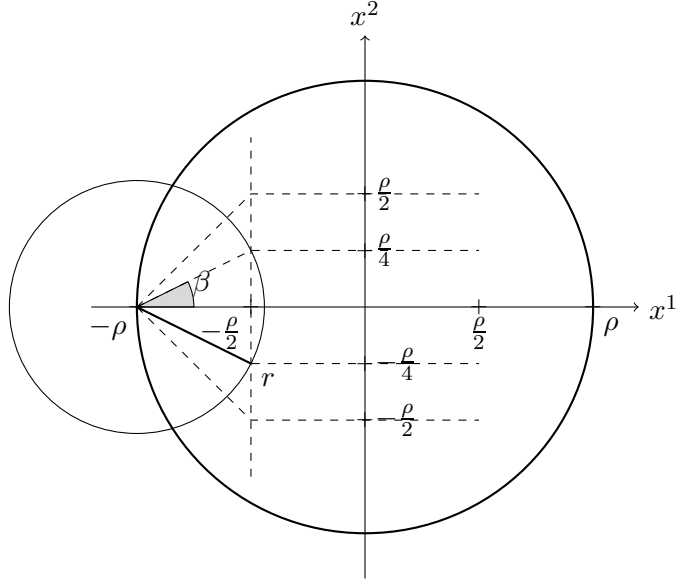
We extend it by homogeneity to $\R^2\setminus \{0\}$,
we choose $\beta = \arccos(2/\sqrt{5})$  and we rescale it of a factor $r=\rho\sqrt{5}/4$ 
so to match the construction indicated in Fig.
\ref{fig:construction-cutoff}:

Such function will satisfy
\begin{align*}
|\nabla\chi_0(x)|
&\leq \left|\chi_0'\left(\frac{x}{|x|}\right)\right|
\left|\nabla\left(\frac{x}{|x|}\right)\right|\leq \frac{C}{|x|},\\
|\nabla^2\chi_0(x)|
&\leq \left|\chi_0''\left(\frac{x}{|x|}\right)\right|
\left|\nabla\left(\frac{x}{|x|}\right)\right|^2
+\left|\chi_0'\left(\frac{x}{|x|}\right)\right|
\left|\nabla^2\left(\frac{x}{|x|}\right)\right|
\leq \frac{C}{|x|^2}.
\end{align*}
We then let $f_1,f_2,f_3:\R\to\R$ be a triple of smooth functions
so that $f_1+f_2+f_3\equiv 1$ and
\begin{align*}
f_1(t)&=
\begin{cases}
1 &\text{if }t\leq -3\rho/4,\\
0&\text{if } t\geq -\rho/2,
\end{cases}
&
f_2(t)&=
\begin{cases}
1 &\text{if }-\rho/2\leq t\leq \rho/2\\
0&\text{if } t\geq 3\rho/4,
\end{cases}
&
f_3(t)&=
\begin{cases}
0 &\text{if }t\leq 3\rho/4,\\
1&\text{if } t\geq \rho,
\end{cases}
\end{align*}
and let $\eta:\R^2\to\R$ be a smooth function so that
\begin{align*}
\eta(x^1,x^2)
=\begin{cases}
1 &\text{if }x^2\in[-\rho/4,\rho/4],\\
0 &\text{if }x^2\in\R^2\setminus [-\rho/2,\rho/2].
\end{cases}
\end{align*}
The required function $\chi$ is then given by
\begin{align*}
\chi(x)
=f_1(x^1)\chi_0(x+(\rho,0))
+f_2(x^1)\eta(x)
+f_3(x^1)\chi_0(x-(\rho,0)),
\quad x\in B_\rho(0).
\end{align*}

\emph{Step 4.}
We claim that the function $\psi$ in the statement of the lemma is
the given by (the restriction to $B_1^+(0)$ of)
the interpolation between $\psi$ and $\psi_0$
through $\chi$, namely
\begin{align*}
\psi(x)
=\chi(x)\phi(x) + (1-\chi(x))\psi_0(x),
\quad x\in B_{\rho}(0).
\end{align*}
Indeed, by construction we have that $\psi$
has flat boundary data on $\rho I$ according to definition \ref{def:geometric-reflection-planar-boundary}
and that
\begin{align*}
\psi = \hat\Phi,
\quad \nabla\psi = \nabla\hat\Phi\quad\text{on } \partial B_{\rho}(0).
\end{align*}
As the proof of lemma \ref{lemma:biharmonic-comparison} shows, to prove
all the estimates in the statement, we have to verify that
\begin{align}
\label{eq:a}
\|\nabla\psi - \nabla L\|_{L^{\infty}(B_\rho(0))}
&\leq C_0\e^{\nu}(\delta+o(\delta)),\\
\label{eq:b}
\|\nabla^2\psi\|_{L^{2}(B_\rho(0))}
&\leq C\|\nabla^2\Phi\|_{L^2(B_{1/4}(0)\setminus B_{1/8}(0))}
+C_0(\delta+o(\delta))\|\nabla\Phi\|_{L^2(B_{\rho}(0))}.
\end{align}
We may write 
\begin{align*}
\psi(x) - L(x)
&=\psi(x) - \psi_0(x) +\psi_0(x)-L(x)
=\chi(x)(\phi(x)-\psi_0(x))+\psi_0(x)-L(x).
\end{align*}
To see that \eqref{eq:a} holds, 
note first of all that from the definition of $\chi$ we
have
\begin{align*}
\nabla\chi(x)
&=\phantom{{}+{}}f_1'(x^1)\chi_0(x+(\rho,0))\vec{\varepsilon}_1
+f_1(x^1)\nabla\chi_0(x+(\rho,0))\\
&\phantom{{}={}}
+f_2'(x^1)\eta(x))\vec{\varepsilon}_1
+f_2(x^1)\nabla\eta(x)\\
&\phantom{{}={}}
+f_3'(x^1)\chi_0(x-(\rho,0))\vec{\varepsilon}_1
+f_3(x^1)\nabla\chi_0(x-(\rho,0)),
\end{align*}
while from the definition of the functions $f, f_1$ and  $f_3$
we have for every $x^1\in [-\rho,\rho]$ that
\begin{align*}
f_1(x^1) f(x_1)&\equiv 0, &
f_1'(x^1) f(x_1)&\equiv 0, &
f_1(x^1)(1- f(x_1))&\equiv f_1(x^1), &
f_1'(x^1)(1- f(x_1))&\equiv f_1'(x^1),\\
f_3(x^1) f(x_1)&\equiv f_3(x^1), &
f_3'(x^1) f(x_1)&\equiv f_3'(x^1), &
f_3(x^1)(1- f(x_1))&\equiv 0, &
f_3'(x^1)(1- f(x_1))&\equiv 0,
\end{align*}
consequently from \eqref{eq:psi-phi}
we deduce the estimate
\begin{align*}
|\nabla\chi(x)(\phi(x)-\psi_0(x))|
\leq C[\nabla\psi_0]_{C^{0,\alpha}(\overline{B_\rho(0)})}
&\bigg(
|f_1'(x^1)|
\left|
\begin{pmatrix}
x^1+\rho\\ x^2
\end{pmatrix}
\right|^{1+\alpha}
+f_1(x^1)\left|
\begin{pmatrix}
x^1+\rho\\
x^2
\end{pmatrix}
\right|^{\alpha}\\
&\phantom{{}={}}+|f_2'(x^1)|
\left|
\begin{pmatrix}
x^1+\rho\\ x^2
\end{pmatrix}
\right|^{1+\alpha}
+f_2(x^1)\left|
\begin{pmatrix}
x^1+\rho\\
x^2
\end{pmatrix}
\right|^{1+\alpha}\\
&\phantom{{}={}}+|f_2'(x^1)|
\left|
\begin{pmatrix}
x^1-\rho\\ x^2
\end{pmatrix}
\right|^{1+\alpha}
+f_2(x^1)\left|
\begin{pmatrix}
x^1-\rho\\
x^2
\end{pmatrix}
\right|^{1+\alpha}\\
&\phantom{{}={}}+|f_3'(x^1)|
\left|
\begin{pmatrix}
x^1-\rho\\ x^2
\end{pmatrix}
\right|^{1+\alpha}
+f_3(x^1)\left|
\begin{pmatrix}
x^1-\rho\\
x^2
\end{pmatrix}
\right|^{\alpha}
\bigg),
\end{align*}
which then implies
\begin{align*}
\|\nabla\chi (\phi-\psi_0)\|_{L^\infty(B_\rho(0))}
\leq C[\nabla\psi_0]_{C^{0,\alpha}{(\overline{B_\rho(0)}})}.
\end{align*}
From estimate \eqref{eq:nabla-psi-phi} we immediately deduce that
\begin{align*}
\|\chi(\nabla\phi - \nabla\psi_0)\|_{L^\infty(\overline{B_\rho(0)})}
\leq C[\nabla\psi_0]_{C^{0,\alpha}{(\overline{B_\rho(0)}})},
\end{align*}
consequently using \eqref{eq:elliptic-est-2-bndry} we can estimate
\begin{align*}
\|\nabla\psi-\nabla L\|_{L^\infty(B_\rho(0))}
&=\|\nabla\chi(\phi - \psi_0)\|_{L^\infty(B_\rho(0))}
+\|\chi(\nabla\phi-\nabla\psi_0)\|_{L^\infty(B_\rho(0))}
+\|\psi_0-L\|_{L^\infty(B_\rho(0))}\\
&\leq C[\nabla\psi_0]_{C^{0,\alpha}(\overline{B_\rho(0)})}
+C_0\e^{\nu}(\delta+o(\delta))\\
&\leq C_0\e^{\nu}(\delta+o(\delta)),
\end{align*}
which proves \eqref{eq:a}.
To establish \eqref{eq:b}, the argument is similar:
from the properties of $\chi$ one deduces the estimates
\begin{align*}
|\nabla^2\chi(x)(\phi(x)-\psi_0(x))|
&\leq C [\nabla\psi_0]_{C^{0,\alpha}(\overline{B_\rho(0)})}
\bigg(
f_1(x^1)
\left|
\begin{pmatrix}
x^1+\rho\\
x^2
\end{pmatrix}
\right|^{-1+\alpha}
+f_3(x^1)
\left|
\begin{pmatrix}
x^1-\rho\\
x^2
\end{pmatrix}
\right|^{-1+\alpha}
\bigg),\\
|\nabla\chi(x)(\nabla\phi(x)-\nabla\psi_0(x))|
&\leq C [\nabla\psi_0]_{C^{0,\alpha}(\overline{B_\rho(0)})}
\bigg(
f_1(x^1)
\left|
\begin{pmatrix}
x^1+\rho\\
x^2
\end{pmatrix}
\right|^{-1+\alpha}
+f_3(x^1)
\left|
\begin{pmatrix}
x^1-\rho\\
x^2
\end{pmatrix}
\right|^{-1+\alpha}
\bigg),
\end{align*}
and, since $\alpha>0$, the right-hand sides 
of these last two inequalities are in $L^2(B_\rho(0))$.
Consequently also thanks to \eqref{eq:elliptic-est-2-bndry}
and \eqref{eq:hessian-psi-phi}
we estimate
\begin{align*}
\|\nabla^2\psi\|_{L^2(B_\rho(0))}
&\leq \|\nabla^2\chi(\phi-\psi_0)\|_{L^2(B_\rho(0))}
+2\|\nabla\chi(\nabla\phi-\nabla\psi_0)\|_{L^2(B_\rho(0))}
+\|\chi(\nabla^2\phi-\nabla^2\psi_0)\|_{L^2(B_\rho(0))}\\
&\leq\|\nabla^2\psi_0\|_{L^2(B_\rho(0))}
+ C[\nabla\psi_0]_{C^{0,\alpha}(\overline{B_\rho(0)})} + 
C_0\e^{\nu}(\delta+o(\delta))\\
&\leq\|\nabla^2\psi_0\|_{L^2(B_\rho(0))}
+C_0\e^{\nu}(\delta+o(\delta)),
\end{align*}
which then implies, thanks to \eqref{eq:elliptic-est-3-bndry},
the estimate \eqref{eq:b}.
The rest of the proof now follows the same lines as in
the proof of lemma \ref{lemma:biharmonic-comparison}.
\end{bfproof}

\subsection{Morrey-Type Estimates and Conclusion}
For a conformal map $\Phi:B_1(0)\to\R^m$
which is a minimiser for the total curvature energy in the class
$\mathcal{F}_{B_1(0)}(\Gamma,\vec{n}_0,A)$
as that in theorem \ref{thm:existence-solution-Poisson-problem},
there are two possibilities.
\par
The first one is that $\Phi$ is a minimal surface, that is
\begin{align*}
D\Area(\Phi)w=0
\quad\text{for all }w\in C^{\infty}_c({B_1(0)},\R^m),
\end{align*}
and this implies that $\Phi$ satisfies
\begin{align*}
\Delta\Phi^i=0\quad\text{in }\mathcal{D}'({B_1(0)})\quad\text{for }i=1,\ldots m.
\end{align*}
The regularity up to the boundary
for the first case is classic, and is essentially the one 
for the Plateau problem,
for which we refer to \cite{MR2760441,MR2566897}.\par
We now study the second possibility.

\begin{lemma}[Interior Morrey-type Estimates]
\label{lemma:interior-Morrey-estimate}
Let $\Phi:B_1(0)\to\R^m$
be a conformal Lipschitz immersion with
$L^2$-bounded second fundamental form
which is an interior minimiser for the total curvature energy in
$\mathcal{F}_{B_1(0)}$ at a fixed area value, that is
\begin{align*}
E(\Phi)\leq E(\Psi)
\end{align*}
for every map $\Psi\in\mathcal{F}_{B_1(0)}$
coinciding with $\Phi$ outside some compact subset of $B_1(0)$
and so that $\Area(\Psi)=\Area(\Phi)$.
Assume that $\Phi$ is not a minimal surface and set
\begin{align}
\label{eq:gamma-D-area}
\zeta=\|D\Area(\Phi)\|>0.
\end{align}
Then, there exists some $r_0>0$ such that 
\begin{align*}
\sup
\bigg\{
r^{-\gamma}\int_{B_r(p)}
(|\nabla^2\Phi|^2_{g_\Phi}+1)
\dd  vol_\Phi :
p\in\overline{B_{1/2}(0)},\, 0<r<r_0
\bigg\}\leq C_0,
\end{align*}
for constants $\gamma>0$ and $C_0>0$ depending only on
$\zeta$ and $\|\nabla\lambda_\Phi\|_{L^{(2,\infty)}(B_1(0))}$.
\end{lemma}
\begin{bfproof}[Proof of lemma \ref{lemma:interior-Morrey-estimate}]
In what follows, we denote by $C$ a positive constant 
(possibly varying line to line) which is independent of $\Phi$,
and with $C_0$ a positive constant
which  depend on  $\|\nabla\lambda_\Phi\|_{L^{(2,\infty)}(B_1(0))}$
and $\zeta$.

\emph{Step 1: constructing a suitable competitor.}
Since $\Phi$ is not a minimal surface,
we may choose some non-zero
$w\in C^\infty_c({B_1(0)},\R^m)$
such that $D\mbox{Area}(\Phi)w>\zeta/2$.
We let $\delta,\varepsilon_0>0$ to be as in
lemma \ref{lemma:biharmonic-comparison}
and whose size is specified in what follows.
Let $r_0>0$ be sufficiently small so that
\begin{align*}
\sup\bigg\{\int_{B_{4r_0}(p)}|\nabla^2\Phi|^2_{g_\Phi}\,\dd vol_\Phi :
p\in\overline{B_{1/2}(0)}
\bigg\}<\varepsilon_0.
\end{align*}

We now fix arbitrary $p\in \overline{B_{1/2}(0)}$
and $0<r<r_0$ and 
(similarly as done as in the proof of lemma \ref{lemma:biharmonic-comparison})
for $\varepsilon_0$ sufficently small as in lemma \ref{lemma:estimates-lambda},
we have that
\begin{align}
\label{eq:sup-estimate-lambda-r-p}
\|\lambda_\Phi-(\lambda_\Phi)_{B_{4r}(p)}\|_{L^\infty(B_{3r}(p))}
\leq C_0,
\end{align}
where $(\lambda_\Phi)_{B_{4r}(p)}$ denotes the average of $\lambda_\Phi$ over $B_{4r}(p)$.
If $\varepsilon_0$ is sufficiently small as in lemma \ref{lemma:conformal-affine-immersion},
then there exists a conformal affine immersion
$L$ whose conformal factor we denote by $\e^\nu$,
such  that the estimates
\begin{align}
\label{eq:W22-affine-approx-r-p}
\|\Phi - L\|_{W^{2,2}(B_{r}(p))}
&<\delta\|\nabla\Phi\|_{L^2({B_{2r}(p)})},\\
\label{eq:logconfact-affine-approx-r-p}
\|\lambda_\Phi - \nu\|_{L^\infty(B_{r}(p))}
&<\delta,
\end{align}
are satisfied.\par
By combining \eqref{eq:sup-estimate-lambda-r-p}-\eqref{eq:logconfact-affine-approx-r-p},
we deduce 
\begin{align*}
\|\lambda_\Phi-\nu\|_{L^\infty(B_{3r}(p))}
&\leq\|\lambda_\Phi-(\lambda_\Phi)_{B_{4r}(p)}\|_{L^\infty(B_{3r}(p))}
+|(\lambda_\Phi)_{B_{4r}(p)}-\nu|\\
&\leq \|\lambda_\Phi-(\lambda_\Phi)_{B_{4r}(p)}\|_{L^\infty(B_{3r}(p))}
+\|\lambda_\Phi-(\lambda_\Phi)_{B_{4r}(p)}\|_{L^\infty(B_{r}(p))}
+\|\lambda_\Phi-\nu\|_{L^\infty(B_{r}(p))}\\
&\leq C_0+\delta,
\end{align*}
consequently we pointwise estimate from above and below
\begin{align}
\label{eq:poinwise-estimate-confact-r-p}
C_0^{-1}(1-\delta-o(\delta))\e^\nu
\leq\e^{\lambda_\Phi(x)}\leq
C_0(1+\delta+o(\delta))\e^\nu
\quad\text{for }x\in B_{3r}(p).
\end{align}
Hence we consider
\begin{align*}
\Psi=
\begin{cases}
\psi &\text{in } B_{\rho}(p),\\
\Phi &\text{in }{B_1(0)}\setminus B_{\rho}(p).
\end{cases}
\end{align*}
where $\rho\in[r/2,r]$ and $\psi$ are given as in lemma
\ref{lemma:biharmonic-comparison}.\par
Thanks to \eqref{eq:estimte-D-area-Phi-psi},
we have that, if $\delta$  is chosen sufficiently small, there holds
\begin{align*}
|D\Area(\Psi)w-D\Area(\Phi)w|
=|D\Area(\psi)w-D\Area((\Phi|_{B_\rho(p)})w|<\frac{\zeta}{4},
\end{align*}
and consequently,
\begin{align*}
D\Area(\Psi)w>\frac{\zeta}{4}>0.
\end{align*}
We consider the function given by
\begin{align*}
a(t)=\Area(\Psi+t w),
\quad t\in\R.
\end{align*}
Let $\varepsilon>0$ be sufficiently small so that,
for every $t\in[-\varepsilon, \varepsilon]$,  $\Psi+tw$
defines a Lipschitz immersion (with $L^2$-bounded
second fundamental form).
Then $a$ is continuously differentiable in $[-\varepsilon,\varepsilon]$
with
\begin{align*}
a'(t)=D\Area(\Psi+t w)w
&=-2\int_{B_1(0)}
\langle \vec{H}_{\Psi+tw},w\rangle,
\,\dd vol_{\Psi+t w},\quad\text{for }t\in[-\varepsilon,\varepsilon],
\end{align*}
and in particular
\begin{align*}
a'(0)&>\frac{\zeta}{4}>0.
\end{align*}
By the inverse function theorem, we deduce that, after
possibly shrinking $\varepsilon$, 
$a$ defines a $C^1$-diffeomorphism
of $[-\varepsilon,\varepsilon]$ onto 
$[\Area(\Psi)-\varepsilon,\Area(\Psi)+\varepsilon]$
and
\begin{align}
\frac{\zeta}{8}\leq a'(t)\leq \frac{\zeta}{2}
\quad\text{for }t\in[-\varepsilon,\varepsilon].
\end{align}
Thanks to \eqref{eq:estimte-D-area-Phi-psi},
we have that, if $\delta$  is chosen sufficiently small, there holds
\begin{align*}
|\Area(\Psi)-\Area(\Phi)|
=|\Area(\psi)-\Area((\Phi|_{B_\rho(p)})|
\leq\frac{\varepsilon}{2},
\end{align*}
so we may find a unique 
$\overline{t}\in [-\varepsilon,\varepsilon]$
so that
\begin{align*}
\Area(\Psi+\overline{t}w)
=\Area(\Phi).
\end{align*}
We then set
\begin{align*}
\overline{\Psi}=
\Psi(x)+\overline{t}w(x)
\quad
 \text{for $x$ in } B_1(0).
\end{align*}
Then $\overline{\Psi}$ is a Lipschitz immersion with $L^2$-bounded second fundamental form,
and by construction there holds $\Area(\Psi)=\Area(\Phi)$.
\medskip

\emph{Step 2: comparison of $\Phi$ with $\overline{\Psi}$.}
By the minimality of $\Phi$ we then have
\begin{align}
\label{eq:minimality-Phi}
\int_{{B_1(0)}}|\vec{\II}_\Phi|^2_{g_\Phi}\,\dd vol_\Phi
\leq \int_{{B_1(0)}}
 |\vec{\II}_{\overline{\Psi}}|^2_{g_{\overline{\Psi}}}
\,\dd vol_{\overline{\Psi}},
\end{align}
Following a computation analogous to (\cite[Lemma A.5]{MR2989995}),
the term on the right-hand-side can be expanded to
\begin{align*}
 \int_{{B_1(0)}}
 |\vec{\II}_{\overline{\Psi}}|^2_{g_{\overline{\Psi}}}
\,\dd vol_{\overline{\Psi}}
=E(\Psi+\overline{t}w)
=E(\Psi)+\overline{t}DE(\Psi)w +R^{\Psi}_w(\overline{t}),
\end{align*}
where $R^\Psi_w(\overline{t})$ is a remainder term satisfying
\begin{align*}
|R^{\Psi}_w(\overline{t})|\leq
C_{\Psi,w}\overline{t}^2,
\end{align*}
and, 
since $\Phi$ is a minimiser for the total curvature energy
with prescribed area,
we may write (we use the divergence form of the Willmore
equation, valid for weak immersions, introduced in \cite{MR2430975}), there holds
for some $c\in\R$,
\begin{align*}
DE(\Psi)w
&=\frac{1}{4}DW(\Psi)w\\
&=\frac{1}{4}\int_{B_1(0)}
\langle
\nabla\vec{H}_{\Psi} -3\pi_{\vec{n}_{\Psi}}(\nabla\vec{H}_{\Psi})
+\star(\nabla^\perp\vec{n}_{\Psi}\wedge\vec{H}_{\Psi}),
\nabla w\rangle\,\dd x\\
&=
\frac{1}{4}\int_{B_1(0)\setminus B_\rho(p)}
\langle
\nabla\vec{H}_{\Phi} -3\pi_{\vec{n}_{\Phi}}(\nabla\vec{H}_{\Phi})
+\star(\nabla^\perp\vec{n}_{\Phi}\wedge\vec{H}_{\Phi}),
\nabla w\rangle\,\dd x\\
&\phantom{{}=}
\frac{1}{4}\int_{B_\rho(p)}
\langle
\nabla\vec{H}_{\psi} -3\pi_{\vec{n}_{\psi}}(\nabla\vec{H}_{\psi})
+\star(\nabla^\perp\vec{n}_{\psi}\wedge\vec{H}_{\psi}),
\nabla w\rangle\,\dd x\\
&=
c\int_{B_1(0)\setminus B_\rho(p)}
\langle
\nabla\Phi,
\nabla w\rangle\,\dd x\\
&\phantom{{}=}
\frac{1}{4}\int_{B_\rho(p)}
\langle
\nabla\vec{H}_{\psi} -3\pi_{\vec{n}_{\psi}}(\nabla\vec{H}_{\psi})
+\star(\nabla^\perp\vec{n}_{\psi}\wedge\vec{H}_{\psi}),
\nabla w\rangle\,\dd x,\\
\end{align*}
so that we can simply estimate:
$
|DE(\Psi)w|
\leq C_{\Phi,w}
$.\par
By the mean value theorem and the estimates
\eqref{eq:estimte-area-Phi-psi}
and \eqref{eq:poinwise-estimate-confact-r-p} it holds
\begin{align*}
|\overline{t}|
&=|a^{-1}(a(\overline{t}))-a^{-1}(a(0))|\\
&\leq \sup_{\xi\in J}|(a^{-1})'(\xi)|
|a(\overline{t})-a(0)|\\
&\leq \frac{8}{\zeta} |\Area(\Psi+tw)-\Area(\Psi)|\\
&\leq C|\Area(\Phi)-\Area(\Psi)|\\
&\leq C(\delta+o(\delta))\Area(\Phi|_{B_\rho(0)})\\
&\leq C_0(\delta+o(\delta))\e^{2\nu},
\end{align*}
 It follows that
 \begin{align*}
 \int_{{B_1(0)}}
 |\vec{\II}_{\overline{\Psi}}|^2_{g_{\overline{\Psi}}}
\,\dd vol_{\overline{\Psi}}
\leq \int_{B_1(0)}
|\vec{\II}_{\Psi}|^2_{g_{\Psi}}
\,\dd vol_{\Psi}
+C_0\e^{2\nu}(\delta+o(\delta)).
\end{align*}
Thanks to \eqref{eq:biharmonic-estimate},
the above estimate and \eqref{eq:minimality-Phi} then imply
\begin{align}
\label{eq:comparison estimate}
\int_{B_{\rho(p)}}
|\vec{\II}_\Phi|^2_{g_\Phi}\,\dd vol_\Phi
 \leq C_0\int_{B_r(p)\setminus B_{r/2}(p)}|\nabla^2\Phi|^2_{g_\Phi}\,\dd  vol_\Phi
+C_0\e^{2\nu}(\delta+o(\delta)).
\end{align}

\emph{Step 3: monotonicity of Area.}
Thanks to \eqref{eq:poinwise-estimate-confact-r-p},
for every
$0<s<r$ we can estimate
\begin{align}
\label{eq:monotonicity-gradient}
\Area(\Phi|_{B_s(p)})
=\int_{B_{s}(p)}|\nabla\Phi|^2\,\dd x
=2\int_{B_{s}(p)}\e^{2\lambda_\Phi}\,\dd x
\leq C_0\e^{2\nu}s^2\leq C_0\e^{4\nu}\frac{s^2}{r^2}\int_{B_{r}(p)}|\nabla\Phi|^2 \dd x.
\end{align}

\emph{Step 4: obtaining the Morrey decrease.}
For any $0<\eta<1/2$, 
thanks to \eqref {eq:pointwise-Hessian-identity} and lemma
\ref{lemma:estimates-lambda},
there exists $C>0$ independent of $p$ and $r$
so that
\begin{align*}
\int_{B_{\eta r}(p)}
(|\nabla^2\Phi|^2_{g_\Phi}\,\dd  vol_\Phi+1)\,\dd vol_\Phi
&\leq 
\left(
\frac{\eta^2}{2}+C\varepsilon_0\right)
\int_{B_{r}(p)}
|\nabla^2\Phi|^2_{g_\Phi}\,\dd  vol_\Phi\\
&\phantom{{}\leq{}}
+\int_{B_{\eta r}(p)}|\vec{\II}_\Phi|^2_{g_\Phi}\,\dd vol_\Phi
+\Area(\Phi|_{B_{\eta r}(p)}),
\end{align*}
and so by estimate \eqref{eq:comparison estimate} we
deduce that there holds
\begin{align*}
\int_{B_{\eta r}(p)}
(|\nabla^2\Phi|^2_{g_\Phi}+1)\dd  vol_\Phi
&\leq
\left(
\frac{\eta^2}{2}+C\varepsilon_0\right)
\int_{B_{r}(p)}
|\nabla^2\Phi|^2_{g_\Phi}\,\dd  vol_\Phi\\
&\phantom{{}\leq{}} +C_0
\int_{B_{2\eta r}(p)
\setminus B_{\eta r}(p)}
|\nabla^2\Phi|^2_{g_\Phi}\,\dd  vol_\Phi\\
&\phantom{{}\leq{}}
+C_0\left(
(\delta+o(\delta))+\eta^2r^2\right)
\e^{2\nu}.
\end{align*}
By using \eqref{eq:poinwise-estimate-confact-r-p} and \eqref{eq:monotonicity-gradient}
and adding 
$C_0\int_{B_{\eta r}(p)}
|\nabla^2\Phi|^2_{g_\Phi}\,\dd  vol_\Phi$ to both hand-sides
and dividing by $1+C_0$ yields
\begin{align*}
\int_{B_{\eta r}(p)}
(|\nabla^2\Phi|^2_{g_\Phi}+1)\dd  vol_\Phi
&\leq
\left(
\frac{\eta^2/2+C\varepsilon_0+C_0}{C_0+1}
\right)
\int_{B_{r}(p)}
|\nabla^2\Phi|^2_{g_\Phi}\,\dd  vol_\Phi\\
&\phantom{{}\leq{}}
+\left(\frac{C_0(\delta+o(\delta)+\eta^2)}{C_0+1}
\right)
\Area(\Phi|_{B_{r}(p)}).
\end{align*}
If $\eta$   and $\delta$ are chosen  sufficiently small so that
\begin{align*}
\beta:=\max\left\{
\frac{\eta^2/2+C\varepsilon_0+C_0}{C_0+1},
\frac{C_0(\delta+o(\delta)+\eta^2)}{C_0+1}
\right\}<1,
\end{align*}
we deduce that
\begin{align*}
\int_{B_{\eta r}(p)}
(|\nabla^2\Phi|^2_{g_\Phi}+1)\,\dd  vol_\Phi
\leq \beta
\int_{B_{r}(p)}
(|\nabla^2\Phi|^2_{g_\Phi}+1)\,\dd  vol_\Phi
\end{align*}
for any $p\in\overline{B_{1/2}(0)}$
and any $0<r<r_0$
where $0<\beta<1$ does not depend on $r$ or $p$.
This inequality can be now iterated and interpolated
to yield
\begin{align*}
\int_{B_{r}(p)}
(|\nabla^2\Phi|^2_{g_\Phi}+1)\,\dd  vol_\Phi
\leq r^{\log_{1/\eta}(1/\beta)}
\frac{1}{r_0^{\log_{1/\eta}(1/\beta)}}
\int_{B_{r_0}(p)}
(|\nabla^2\Phi|^2_{g_\Phi}+1)\,\dd  vol_\Phi
\end{align*}
for any $r<\eta r_0$,
where $\beta$ and $\eta$ depend only on $\zeta$ and
$\|\nabla\lambda_\Phi\|_{L^{(2,\infty)}(B_1(0))}$.
After relabelling $r_0$, and setting $\gamma:=\log_{1/\eta}(1/\beta)$ we can  concludes
of the proof of the lemma.
\end{bfproof}

For boundary points, we have the following.
Recall that we write $\partial B_1^+(0)= I + S$,
where $I$ is the base diameter and $S$ is the upper semi-circle,
and when we say that a Lipschitz immersion ``has geometric boundary data of class $C^{1,1}$''
if its boundary curve and boundary Gauss map are $C^{1,1}$ up to re-parametrization
and not in a point-wise sense, see definition \ref{def:weak-immersions}--(ii).
\begin{lemma}[Boundary Morrey-type Estimates]
\label{lemma:Morrey-estimate-bndry}
Let $\Phi:B_1(0)^+\to\R^m$
be a conformal Lipschitz immersion with
$L^2$-bounded second fundamental form so that
\begin{align*}
E(\Phi)\leq E(\Psi)
\end{align*}
for every map $\Psi\in\mathcal{F}_{B_1^+(0)}$ that coincides with $\Phi$ outside
some subset $K$ of $\overline{B_1^+(0)}$ with $\dist(K,S)>0$
and having the same geometric boundary data of $\Psi$ in $I$, that is
\begin{align*}
\Psi(I)=\Phi(I)
\quad\text{and}\quad
\vec{n}_\Psi(I)=\vec{n}_\Phi(I),
\end{align*}
and the same area, that is $\Area(\Phi)=\Area(\Phi)$.
Assume that the boundary geometric data of $\Phi$ on $I$ are of class $C^{1,1}$
and that $\Phi$ is not a minimal surface.
Let $\zeta>0$ be as in \eqref{eq:gamma-D-area}.
Then,
there exist some $0<\overline{r}<r_0<1$
so that
\begin{align*}
\sup
\bigg\{
r^{-\gamma}\int_{B_r^+(p)}|\nabla^2\Phi|^2_{g_\Phi}
\,\dd vol_\Phi :
p\in r_0 I,\, 0<r<\overline{r}
\bigg\}< +\infty,
\end{align*}
for some constant $\gamma>0$.
\end{lemma}

\begin{remark}
\label{rmk:facts-for-Morrey-estimate-bndry}
Two elementary facts that will be used in the proof of lemma
\ref{lemma:Morrey-estimate-bndry} are the following:
\begin{enumerate}[(i)]
\item If $\vec{e}=(\vec{e}_1,\vec{e}_2)$ denotes the
coordinate frame of the map $\Phi$, then we have
\begin{align*}
\vec{e}(x) = (\mathbf{t},\star(\mathbf{t}\wedge\vec{n}_0)(\sigma(x^1)),
\quad x^1\simeq(x^1,0)\in I,
\end{align*}
where $\mathbf{t}$ denotes the tangent vector of the boundary curve
and $\sigma_\Phi$ is some homeomorphism with $\sigma'(x^1)=\e^{\lambda_\Phi(x^1,0)}$.
In particular, since the boundary data are assumed of class $C^{1,1}$, we see that
for $i=1,2$ we can estimate, for every $1<p<\infty$,
\begin{align*}
\|\partial_\tau\vec{e}_i\|_{L^p(I)}
\leq C_0\|\e^{\lambda_\Phi(\cdot,0)}\|_{L^p(I)},
\end{align*}
where $C_0$ depends only on the geometric boundary data.
Furthermore, for every $1<p<\infty$,
if we set $\Phi_r(x)=\Phi(rx)$, $x\in B_1^+(0)$
one can compute that
$\|\e^{\lambda_{\Phi_r}(\cdot,0)}\|_{L^p(I)} = r^{\frac{p-1}{p}}\|\e^{\lambda_{\Phi}(\cdot,0)}\|_{L^p(rI)}$
and thus deduce that the $L^p$-norm of $\e^{\lambda_\Phi(\cdot,0)}$
is \emph{decreasing} with respect to rescaling in the domain, for $0<r<1$.
\item For a generic immersion of an open domain $\Omega\subset\R^m$,
 $X:\Omega\to\R^m$ and a diffeomorphism 
$f:\R^m\to\R^m$, denoting for brevity
denoting
$
\vartheta=
O(\|\nabla f- \mathbbm{1}_{m\times m}\|_{L^\infty(\R^m)})$ and
$
\eta=O(\|\nabla^2 f\|_{L^\infty(\R^m)})
$,
we can deduce the point-wise estimates
\begin{align*}
\left(
1-\vartheta
\right)
g_{X}
&\leq g_{f\circ X}\leq(1+\vartheta)g_{X},\\
(1-\vartheta)\dd vol_X
&\leq
\dd vol_{{f\circ X}}\leq(1+\vartheta)\dd vol_X,\\
|\nabla^2(f\circ X)|_{g_{f\circ X}}^2
&\leq (1+\vartheta)
(\eta + \vartheta |\nabla ^2 X|^2_{g_X})
\leq C_0(1+|\nabla^2 X|^2_{g_X}),
\end{align*}
where $C_0$ is a constant depending only on $f$.
 
\end{enumerate}
\end{remark}

\begin{bfproof}[Proof of Lemma \ref{lemma:Morrey-estimate-bndry}]
In what follows, we denote by $C$ a positive constant 
(possibly varying line to line) which is independent of $\Phi$,
and with $C_0$ a positive constant
dependeding only on $\|\nabla\lambda_\Phi\|_{L^{(2,\infty)}(B_1(0))}$,
and on the geometric boundary data at $I$ of $\Phi$.

\emph{Step 1: preliminaries and reductions.}
We fix    $p>1$ and a suitably small $\varepsilon_0$ that will be specified below.
Since the boundary data of $\Phi$ along $I$
are of class $C^{1,1}$, we may find some $0<r_0<1$ 
and a  $C^{1,1}$-homeomorphism $f:\R^m\to\R^m$
(that is, $f$ and its inverse belong to $C^{1,1}(\R^m,\R^m)$),
so that $(f\circ\Phi)|_{B_{r_0}^+(0)}$ has flat boundary data along $r_0 I$
as in the sense of definition \ref{def:geometric-reflection-planar-boundary}. 
Up to further reducing $r_0$, we also assume that
\begin{align}
\label{eq:initial-smallness}
\int_{B_{r_0}^+(0)}
(|\nabla^2\Phi|^2_{g_{\Phi}}
+1)\,\dd vol_{\Phi}
+\|\e^{\lambda_{\Phi}}\|_{L^p(r_0I)} <\varepsilon_0.
\end{align}
Note that, since 
the Lagrangian $\Phi\mapsto \int (1+|\nabla^2\Phi^2|_{g_\Phi}^2)\,\dd vol_\Phi$
is invariant with respect to re-parametrizations
and the conformal factor is decreasing with respect 
to rescalings (see remark \ref{rmk:facts-for-Morrey-estimate-bndry}--(i)),
\eqref{eq:initial-smallness} implies in particular that,
having set for brevity $\Phi_{r_0}(x)=\Phi(r_0 x)$,
there holds
\begin{align*}
\sup\bigg\{\int_{B_{4r}^+(x)}
(|\nabla^2\Phi_{r_0}|^2_{g_{\Phi_{r_0}}}+1)\,\dd vol_{\Phi_{r_0}}
+\|\e^{\lambda_{\Phi_{r_0}}}\|_{L^p(x+rI)}:
x\in \frac{1}{2}I,\; 0<r\leq \frac{1}{16}
\bigg\}<\varepsilon_0.
\end{align*}
We will work from now on, omitting the subscript,
with $\Phi_{r_0}(x)=\Phi(r_0 x)$ in place of $\Phi$.
Also, since $\Phi$ is not a minimal surface, there exists some
$w\in C^{\infty}_c(B_1^+(0),\R^m)$ so that
\begin{align}
\label{eq:non-minimal-surface}
D\Area(\Phi)w \geq \zeta/2>0.
\end{align}
Finally, thanks to lemma \ref{lemma:conformal-coordinates} there exist some bi-Lipschitz homeomorphism
$\phi:B_1^+(0)\to B_1^+(0)$ so
that $f\circ\Phi\circ \phi$ is conformal
moreover up to further composing with a conformal
self-map of $B_1^+(0)$ we may suppose that $\phi(\pm 1)=\pm 1$ and $\phi(0)=0$,
hence that the geometric boundary data on $I$ are sent (globally) onto themselves: $f(\Phi(I))=f(\Phi(\phi(I)))$.

\emph{Step 2}.
For simplicity of notation we will prove the Morrey-type decay at
$x=0$; the one for other points in $(1/2)I$ is analogous.
Since $\phi$ is bi-Lipschitz, we may find a sufficiently big $N\in \N$
and a sufficiently big $M=M(N)$ so that, for every $0<r<1$, we have
\begin{align}
\label{eq:bi-lip-M-N}
B_{r/2^M}^+(0)\subset\phi(B_{r/{2^{N}}}(0))\subset B_{r}^+(0).
\end{align}
Let $0<r\leq 1/16$ be fixed.
For a sufficiently small $\varepsilon_0$, thanks to lemma 
\ref{lemma:biharmonic-comparison-bndry} we may find a 
$\rho\in [r/2^{N+1},r/2^N]$
and an immersion
$\psi\in C^{1,\alpha}(\overline{B_\rho^+(0)},\R^m)$
which satisfies
\begin{align*}
\begin{aligned}
\psi&=f\circ\Phi\circ\phi &&\text{on }\partial B_{\rho}(0)\cap B_1^+(0),\\
\nabla\psi&=\nabla(f\circ\Phi\circ\phi) &&\text{on }\partial B_{\rho}(0)\cap B_1^+(0),
\end{aligned}
\end{align*}
has flat boundary data on $\rho I$
and satisfies the estimates
\eqref{eq:biharmonic-estimate-bndry},
\eqref{eq:estimate-area-Phi-psi-bndry} and 
\eqref{eq:estimate-D-area-Phi-psi-bndry}
with $r/2^N$ in place of $r$ and $f\circ\Phi\circ\phi$ in place of $\Phi$,
and in particular
\begin{align}
\label{eq:most-important-estimate-psi}
\int_{B^+_{r/2^{N+1}}(0)}
|\nabla^2\psi|^2_{g_\psi}\,\dd vol_{\psi}
\leq C_0\int_{B^+_{r/2^N}(0)\setminus B_{r/2^{N+1}}(0)}
|\nabla^2(f\circ\Phi\circ\phi)|^2\dd vol_{f\circ\Phi\circ\phi}
+C_0\int_{B^+_{r/2^N}(0)}\dd vol_{f\circ\Phi\circ\phi}
\end{align}

Hence we set
\begin{align*}
\Psi=
\begin{cases}
f^{-1}\circ\psi\circ\phi^{-1} &\text{in }\phi(B_\rho^+(0)),\\
\Phi &\text{in }B_1^+(0)\setminus\phi(B_\rho^+(0))
\end{cases}
\end{align*}

From \eqref{eq:estimate-area-Phi-psi-bndry} and
\eqref{eq:estimate-D-area-Phi-psi-bndry}
we deduce
\begin{align*}
|\Area(\Phi|_{\phi(B_\rho^+(0))})-\Area(f^{-1}\circ\psi\circ\phi^{-1})|
&\leq C_0(\delta+o(\delta))\Area(\Phi|_{\phi(B_\rho^+(0))}),\\
\|D\Area(\Phi|_{\phi(B_\rho^+(0))})-D\Area(f^{-1}\circ\psi\circ\phi^{-1})\|
&\leq C_0(\delta+o(\delta))\Area(\Phi|_{\phi(B^+_{r/2^N}(0))}),
\end{align*}
where $C_0>0$ depends on $\|\nabla\lambda_\Phi\|_{L^{2\infty}(B_1^+(0))}$
and on $f$.
Thanks to \eqref{eq:non-minimal-surface},
we have that, if $\delta$ (and accordingly $\varepsilon_0$)
is chosen sufficiently small, there holds
\begin{align*}
|D\Area(\Phi|_{\phi(B_\rho^+(0))})w-D\Area(f^{-1}\circ\psi\circ\phi^{-1})w|
<\frac{\zeta}{4},
\end{align*}
and consequently,
$
D\Area(\Psi)w>\zeta/4.
$
We then  consider the function given by
\begin{align*}
a(t)=\Area(\Psi+t w),
\quad t\in\R.
\end{align*}
Let $\varepsilon>0$ be sufficiently small so that,
for every $t\in[-\varepsilon,\varepsilon]$,  $\Psi+tw$
defines a Lipschitz immersion (with $L^2$-bounded
second fundamental form).
Then $a$ is continuously differentiable in $[-\varepsilon,\varepsilon]$
with
\begin{align*}
a'(t)=D\Area(\Psi+t w)w
&=-2\int_{B_1^+(0)}
\langle \vec{H}_{\Psi+tw},w\rangle,
\,\dd vol_{\Psi+t w},\quad\text{for }t\in[-\varepsilon,\varepsilon],
\end{align*}
and in particular
$
a'(0)>\zeta/4>0.
$
By the inverse function theorem, we deduce that, after
possibly shrinking $\varepsilon$, 
$a$ defines a $C^1$-diffeomorphism
of $[-\varepsilon,\varepsilon]$ onto 
$[\Area(\Psi)-\varepsilon,\Area(\Psi)+\varepsilon]$
and
\begin{align}
\frac{\zeta}{8}\leq a'(t)\leq \frac{\zeta}{2}
\quad\text{for }t\in[-\varepsilon,\varepsilon].
\end{align}
By choosing $\delta$ sufficiently small we may suppose that  
\begin{align*}
|\Area(\Phi)-\Area(\Psi)|
=|\Area(\Phi|_{\phi(B_\rho^+(p))})-\Area(f^{-1}\circ\psi\circ\phi^{-1})|
\leq\frac{\varepsilon}{2},
\end{align*}
so we may find a unique 
$\overline{t}\in [-\varepsilon,\varepsilon]$
so that
$
\Area(\Psi+\overline{t}w)
=\Area(\Phi).
$
We then set
\begin{align*}
\overline{\Psi}=
\Psi(x)+\overline{t}w(x)
\quad
 \text{for $x$ in } B_1^+(0).
\end{align*}
Then $\overline{\Psi}$ is a Lipschitz immersion with $L^2$-bounded second fundamental form,
and by construction there holds $\Area(\overline\Psi)=\Area(\Phi)$.
Similarly as done in the proof of lemma \ref{lemma:interior-Morrey-estimate},
following a computation analogous to (\cite[Lemma A.5]{MR2989995}),
the total curvature energy of $\overline\Psi$ can be expanded with respect to $\overline{t}$
as
\begin{align*}
 \int_{{B_1^+(0)}}
 |\vec{\II}_{\overline{\Psi}}|^2_{g_{\overline{\Psi}}}
\,\dd vol_{\overline{\Psi}}
=E(\Psi+\overline{t}w)
=E(\Psi)+\overline{t}DE(\Psi)w +R^{\Psi}_w(\overline{t}),
\end{align*}
with $|DE(\Psi)w|\leq C_{\Phi,w}$ and
$R^\Psi_w(\overline{t})$ satisfies
$
|R^{\Psi}_w(\overline{t})|\leq
C_{\Psi,w}\overline{t}^2.
$
By the mean value theorem,
we have the estimate
\begin{align*}
|\overline{t}|
&=|a^{-1}(a(\overline{t}))-a^{-1}(a(0))|\\
&\leq \sup_{\xi\in J}|(a^{-1})'(\xi)|
|a(\overline{t})-a(0)|\\
&\leq \frac{8}{\zeta} |\Area(\overline\Psi)-\Area(\Psi)|\\
&= \frac{8}{\zeta} |\Area(\Phi)-\Area(\Psi)|\\
&\leq C_0(\delta+o(\delta))\Area(\Phi|_{\phi(B_\rho^+(0))}),
\end{align*}
  where $C_0$ depends on 
$\|\nabla\lambda_\Phi\|_{L^{2\infty}(B_1^+(0))}$, on $f$
and on $\zeta$ and
this yields the estimate
\begin{align*}
 \int_{{B_1(0)}}
 |\vec{\II}_{\overline{\Psi}}|^2_{g_{\overline{\Psi}}}
\,\dd vol_{\overline{\Psi}}
\leq \int_{B_1(0)}
|\vec{\II}_{\Psi}|^2_{g_{\Psi}}
\,\dd vol_{\Psi}
+C_0(\delta+o(\delta))\Area(\Phi|_{\phi(B_\rho^+(0))}).
\end{align*}
We write
\begin{align*}
\int_{\phi(B_{r/2^{N+1}}(0))}|\nabla^2\Phi|^2_{g_\Phi}
\,\dd vol_\Phi
=\int_{\phi(B_{r/2^{N+1}}(0))}|\vec{\II}_\Phi|^2_{g_\Phi}
\,\dd vol_\Phi
+4\int_{\phi(B_{r/2^{N+1}}^+(0))}
|\nabla\lambda_\Phi|^2\,\dd x.
\end{align*}
We have, on the one hand, thanks to lemma \ref{lemma:estimates-lambda-bndry}
and the choice of $N$,
the estimate
\begin{align}\label{morrey-lambda}
&\phantom{{}\leq{}}
\int_{\phi(B_{r/2^{N+1}}^+(0))}
|\nabla\lambda_\Phi|^2\,\dd x
\leq
\int_{B_{r/2}^+(0)}
|\nabla\lambda_\Phi|^2\,\dd x\\
&\leq \frac{1}{8}\int_{B_r^+(0)}|\nabla^2\Phi|^2_{g_\Phi}
\,\dd vol_\Phi
+C_0\varepsilon_0 
\bigg(
\int_{B_r^+(0)}
|\nabla^2\Phi|^2_{g_\Phi}\,\dd vol_{\Phi}
+C_0\|\e^{\lambda_\Phi}\|_{L^p(rI)}
\bigg),\nonumber
\end{align}
where $C_0$ depends only on $\zeta$ and $\vec{n}_0$;
on the other hand, the comparison of $\Phi$ with $\overline\Psi$
yields (we use the pointwise a.e. estimates in 
remark \ref{rmk:facts-for-Morrey-estimate-bndry}--(ii))

\begin{align}
&\phantom{{}\leq{}}
\int_{\phi(B_{r/2^{N+1}}^+(0))}
|\vec{\II}_\Phi|^2_{g_\Phi}\,\dd vol_\Phi\nonumber\\
&\leq
\int_{\phi(B_{r/2^{N+1}}^+(0))}
|\vec{\II}_{f^{-1}\circ\psi\circ\phi}|^2_{g_{f^{-1}\circ\psi\circ\phi}}
\,\dd vol_{{f^{-1}\circ\psi\circ\phi}}\nonumber\\
&= \int_{B_{r/2^{N+1}}^+(0)}
|\vec{\II}_{f^{-1}\circ\psi}|^2_{g_{f^{-1}\circ\psi}}
\,\dd vol_{f^{-1}\circ\psi}\nonumber\\
&\leq\int_{B_{r/2^{N+1}}^+(0)}
|\nabla^2(f^{-1}\circ\psi)|^2_{g_{f^{-1}\circ\psi}}
\,\dd vol_{f^{-1}\circ\psi}\nonumber\\
&\leq C_0\int_{B_{r/2^{N+1}}^+(0)}
(1+|\nabla^2\psi|^2_{g_{\psi}})
\,\dd vol_\psi &&\text{(by \eqref{eq:most-important-estimate-psi})}
\nonumber\\
&\leq C_0\int_{B_{r/2^{N+1}}^+(0)}
\,\dd vol_\psi \nonumber\\
&\phantom{{}\leq{}}
+ C_0\int_{B_{r/2^{N}}^+(0)\setminus B_{r/2^{N+1}}^+(0)}
|\nabla^2(f\circ\Psi\circ\phi)|^2_{g_{f\circ\Psi\circ\phi}}
\,\dd vol_\psi
+ C_0\int_{B_{r/2^{N}}^+(0)}
\,\dd vol_{f\circ\Psi\circ\phi}\nonumber\\
&\leq C_0\int_{\phi(B_{r/2^{N}}^+(0))}
\,\dd vol_\Phi
+ C_0\int_{\phi(B_{r/2^{N}}^+(0)\setminus B_{r/2^{N+1}}^+(0))}
|\nabla^2\Phi|^2_{g_{\Phi}}
\,\dd vol_\Phi\nonumber,
\end{align}
and so all in all,
\begin{align}
\label{morrey-general}
\int_{\phi(B_{r/2^{N+1}}^+(0))}
|\vec{\II}_\Phi|^2_{g_\Phi}\,\dd vol_\Phi
\leq C_0\int_{\phi(B_{r/2^{N}}^+(0))}
\,\dd vol_\Phi
+ C_0\int_{\phi(B_{r/2^{N}}^+(0)\setminus B_{r/2^{N+1}}^+(0))}
|\nabla^2\Phi|^2_{g_{\Phi}}
\,\dd vol_\Phi
\end{align}
By combining \eqref{morrey-lambda} and  \eqref{morrey-general} we deduce that
\begin{align*}
\int_{\phi(B_{r/2^{N+1}}^+(0))}
|\nabla^2\Phi|^2_{g_\Phi}\,\dd vol_\Phi
&\leq 
 C_0\int_{\phi(B_{r/2^{N}}^+(0))}
\,\dd vol_\Phi
+ C_0\int_{\phi(B_{r/2^{N}}^+(0)\setminus B_{r/2^{N+1}}^+(0))}
|\nabla^2\Phi|^2_{g_{\Phi}}
\,\dd vol_\Phi\\
&\phantom{{}\leq{}}
+\frac{1}{2}\int_{B_r^+(0)}|\nabla^2\Phi|^2_{g_\Phi}
\,\dd vol_\Phi
+C_0\varepsilon_0 
\bigg(
\int_{B_r^+(0)}
|\nabla^2\Phi|^2_{g_\Phi}\,\dd vol_{\Phi}
+C_0\|\e^{\lambda_\Phi}\|_{L^p(rI)}
\bigg)
\end{align*}
and so adding 
$
C_0
\int_{\phi(B_{r/2^{N+1}}^+(0))}
|\nabla^2\Phi|^2_{g_\Phi}\,\dd vol_\Phi
$
to both hand-sides 
\begin{align*}
\int_{\phi(B_{r/2^{N+1}}^+(0))}
|\nabla^2\Phi|^2_{g_\Phi}\,\dd vol_\Phi
&\leq
\bigg(
\frac{1/2 + C_0\varepsilon_0+C_0}{C_0+1}
\bigg)
\int_{B_r^+(0)}
|\nabla^2\Phi|^2_{g_\Phi}\,\dd vol_{\Phi}\\
&\phantom{{}\leq{}}
+C_0\int_{B_r^+(0)}
\,\dd vol_\Phi
+C_0\|\e^{\lambda_\Phi}\|_{L^p(rI)}.
\end{align*}
Choosing  $\varepsilon_0$ sufficiently small so that 
$C_0\varepsilon_0\leq 1/4$, with \eqref{eq:bi-lip-M-N}
and the above inequality
we deduce that, for every $0<r<1/4$ there holds
\begin{align*}
\int_{B_{r/2^{M+1}}^+(0)}
|\nabla^2\Phi|^2_{g_\Phi}\,\dd vol_\Phi
\leq \beta
\int_{B_r^+(0)}
|\nabla^2\Phi|^2_{g_\Phi}\,\dd vol_{\Phi}
+C_0\|\e^{\lambda_\Phi}\|_{L^\infty(B_1^+(0))}r,
\end{align*}
 where $0<\beta<1$ depends only on $C_0$.
As in the proof of lemma \ref{lemma:interior-Morrey-estimate},
this equality can now be iterated and interpolated to yield
the existence of some $\gamma>0$ so that
\begin{align*}
\sup_{r<r<\overline{r}}
\bigg\{
r^{-\gamma}\int_{B_{r}^+(0)}
|\nabla^2\Phi|^2_{g_\Phi}\,\dd vol_\Phi
\bigg\}<+\infty,
\end{align*}
for some suitably small $\overline{r}<1/16$.
After going back to the original scale,
this yields to the conclusion of the proof of the lemma.
\end{bfproof}

\begin{bfproof}[Proof of Theorem
\ref{thm:existence-solution-Poisson-problem} (conclusion)
and of Theorem
 \ref{thm:regularity-solution-Poisson-problem}]
Let $\Phi$ be a conformal map which is a  minimiser for the total curvature energy in
$\mathcal{F}_{B_1(0)}(\Gamma,\vec{n}_0,A)$,
and let $a_1,\ldots,a_N$ be its branch points 
(recall that none of them lays on the boundary).
For every sufficiently small $\delta>0$, the conformal factor
of $\Phi$, $\e^{\lambda_\Phi}=|\nabla\Phi|^2/\sqrt{2}$
is then uniformly bounded from above and below
in $\overline{B_1(0)}\setminus\cup_{i=1}^NB_\delta(a_i)$,
and consequently, covering $\overline{B_1(0)}\setminus\cup_{i=1}^NB_\delta(a_i)$
with finitely many balls, thanks to
lemma \ref{lemma:interior-Morrey-estimate}
and lemma
\ref{lemma:Morrey-estimate-bndry},
we deduce that its Hessian $\nabla^2\Phi$ belongs to the Morrey space 
$L^{2,a}(B_1(0)\setminus\cup_{i=1}^NB_\delta(a_i))$
for some $a>0$ (see e.g. \cite{MR717034,MR3099262}),
and consequently, by Morrey's Dirichlet growth theorem (\cite{MR0202511}, see also 
\cite[Theorem 5.7]{MR3099262}) that 
$\Phi\in C^{1,a/2}(\overline{B_1(0)}\setminus\cup_{i=1}^NB_\delta(a_i))$.

Similarly as done in \cite{MR2430975}, $\Phi$
satisfies the Euler-Lagrange equation for the Poisson problem, i.e.
the Willmore equation in divergence form
plus a Lagrange multiplier term for the area constraint,
in the sense of distributions away from the branch points.
Being however $\vec{H}_\Phi=\e^{-2\lambda_\Phi}\Delta\Phi/2$
in $L^2(B_1(0),\R^m)$, the contributions at the branch points can only consist on
sum of Dirac masses and thus we deduce that
\begin{align*}
\divop\left(
\nabla\vec{H}_{\Phi} -3\pi_{\vec{n}_{\Phi}}(\nabla\vec{H}_{\Phi})
+\star(\nabla^\perp\vec{n}_\Phi\wedge\vec{H}_{\Phi})
+ c\nabla \Phi
\right)
=\sum_{i=1}^N\vec{\alpha}_i\delta_{a_i},
\quad\text{in }\mathcal{D}'(B_1(0),\R^m),
\end{align*}
for some $\alpha_N,\ldots\alpha_N\in\R^m$
and $c\in\R$.
Since however $\Phi$ is a variational immersion,
i.e. it is obtained as a critical point of the total curvature energy,
a simple implicit function theorem argument similar to that
used in the proof of lemmas 
\ref{lemma:interior-Morrey-estimate}
and
\ref{lemma:Morrey-estimate-bndry}
allows to construct suitable comparisons for $\Phi$
thanks to which one deduces that
\begin{align*}
\vec\alpha_i=0
\quad\text{for every\quad}
\quad i=1,\ldots,N.
\end{align*}
(We refer the reader to \cite{PalmThesis} for the details of the above statements).
This yields that $\Phi$ satisfies the area-constrained
Willmore equation through the branch points,
and thus, from the analysis of singularities for Willmore
surfaces \cite{MR2430975,MR3096502,MR2119722,MR2318282},
this implies that $\Phi$ is of class $C^{1,\alpha}$ for every $0<\alpha<1$
through the branch points.

We have thus proved that
$\Phi$ is of class $C^{1,\alpha}$-up to the boundary, for some $0<\alpha<1$
(and accordingly the Gauss map $\vec{n}_\Phi$ extends to a map of class 
$C^{0,\alpha}$-up to the boundary)
and this concludes the proof of the theorems.
\end{bfproof}

\appendix
\section{Appendix}
\subsection{Existence of Immersed Disks with Given Boundary Data}
\label{app:contruction-competitors}
In this section we prove,
along with some other facts and comments,
the following result.
\begin{lemma}
\label{lemma:existence-competitors}
Let $\Gamma\subset\R^m$ be a simple, closed curve of class
$C^{k,\alpha}$ for $k\in \N_{\geq 1}$ and $\alpha\in(0,1]$
whose unit tangent vector we denote by $\mathbf{t}$
and let $\vec{n}_0$ be a unit-normal $(m-2)$-vector field along
$\Gamma$ of class $C^{k,\alpha}$.
There exists a possibly branched immersed disk $\Phi:B_1(0)\to\R^m$
of class $C^{k,\alpha}$ and boundary  $\Gamma$ and 
whose Gauss map along $\Gamma$ is $\vec{n}_0$.
In particular, a branch-point-free immersion  $\Phi$ 
can be produced when either $m>3$ or when
$m=3$ and the map
\begin{align*}
x\mapsto
(\mathbf{t}\times\vec{n}_0,\mathbf{t},\vec{n}_0)(x),\quad x\in S^1,
\end{align*}
defines a non-nullhomotopic loop in $SO(3)$.
\end{lemma}
Note that, when $k=\alpha=1$, or when $k\geq 2$,
then $\Phi$ in the above lemma satisfies the
assumptions of lemma \ref{lemma:conformal-coordinates} 
and can thus be conformally re-parametrised.

We treat the case $m=3$;
when  $m\geq 3$
see the final remark \ref{remark:higher-codim}.
For the elementary concepts of algebraic topology here mentioned
we refer the reader to \cite{MR1867354,MR807945, MR1138462}.

Any (non-branched) immersion $\Phi:\overline{B_1(0)^2}\to\R^3$
naturally defines a map into the space 
invertible of matrices with positive determinant, $E=E_\Phi:\overline{B_1(0)^2}\to GL^+(3,\R)$,
by
\begin{align*}
E(x)=(\partial_1\Phi(x),\partial_2\Phi(x),\vec{n}_\Phi(x)),
\quad x\in\overline{B_1(0)^2},
\end{align*}
where $n_\Phi$ denotes the Gauss map of $\Phi$.
The classical Gram-Schmidt algorithm
gives the existence of a deformation retraction of $GL^+(3,\R)$  to
the 3-dimensional special orthogonal group $SO(3)$, and in particular the map $E$ is homotopic in $GL^+(3,\R)$
to the coordinate frame map
\begin{align*}
e(x)=e_\Phi(x) = (e_1(x),e_2(x),  e_3(x)),
\quad x\in \overline{B_1(0)^2},
\end{align*}
where
\begin{align*}
e_1(x)&=\frac{\partial_1\Phi(x)}{|\partial_1\Phi(x)|},\\
e_2(x)&=\frac{\partial_2\Phi(x)}{|\partial_2\Phi(x)|}
-\bigg\langle\frac{\partial_2\Phi(x)}{|\partial_2\Phi(x)|},e_1(x)\bigg\rangle e_1(x),\\
e_3(x)&=e_1\times e_2(x)=\vec{n}_\Phi(x).
\end{align*}
We can similarly define the \emph{polar frame map}
defined by means of polar coordinates
$x=r\e^{i\theta}$ in $\overline{B_1(0)^2}\setminus\{0\}$ as
$p(x)=(p_1(x),p_2(x),p_3(x))$,
where
\begin{align*}
p(r\e^{i\theta})
=(e_1,e_2,e_3)(r\e^{i\theta})
\begin{pmatrix}
\cos\theta & -\sin\theta & 0\\
\sin\theta & \cos\theta &0\\
0 &0 & 1
\end{pmatrix}.
\end{align*}
We recall that the fundamental group of $SO(3)$ consists precisely 
of two components:
\begin{align*}
\pi_1(SO(3))=\Z/2\Z,
\end{align*}
the non-trivial class being represented for instance
by the family realising a complete rotation
around the $z$-axis:
\begin{align}
\label{eq:family-rotations-z-axis}
R(\theta,\hat{z})
=\begin{pmatrix}
\cos\theta &-\sin\theta &0\\
\sin\theta &\cos\theta &0\\
0&0&1
\end{pmatrix},\quad\theta\in[0,2\pi].
\end{align} 
Recall moreover that, being $SO(3)$ a topological group,
the matrix product operation is compatible with the one of $\pi_1(SO(3))$.
Since the restriction of the coordinate frame map $e$ to $S^1=\partial B_1(0)$ defines a nullohomotopic loop
in $SO(3)$, the homotopy being induced by the immersion:
\begin{align*}
e_t(x)=e\left(tx\right),\quad x\in\partial B_1(0),\,t\in[0,1].
\end{align*}
the polar frame defines then a non-contractible loop in $SO(3)$.
This argument implies that, given an immersed curve $\gamma: S^1\to\R^3$
and a unit-normal vector field $\vec{n}_0:S^1\to S^2$ along $\gamma$,
a necessary condition for the existence of an immersion $\Phi:\overline{B_1(0)}\to\R^3$
bounding $\gamma$ and so that $\vec{n}_\Phi=\vec{n}_0$ on $\partial B_1(0)$
is that
\begin{align*}
x\mapsto(\mathbf{t}\times \vec{n}_0,\mathbf{t},\vec{n}_0)(x),
\quad x\in S^1,
\end{align*}
is \emph{not} a nullhomotopic loop in $SO(3)$.
Examples of couples $(\gamma, n_0)$ that do not satisfy this condition are
easy to produce.

\begin{example}[Dirac Belt]
\label{example:dirac-belt}
Let $\gamma:[0,2\pi]\to\R$ be the unit circle:
\begin{align*}
\gamma(\theta)=(\cos\theta,\sin\theta,0).
\end{align*}
We consider a rotation of angle $\theta$ around the
tangent vector of $\gamma$, namely
$\mathbf{t}(\theta)=(-\sin\theta,\cos\theta,0)$.
Its matrix is given by
\begin{align*}
R(\theta,\mathbf{t}(\theta))
=B^T(\theta)R(\theta,\hat{z})B(\theta),
\end{align*}
where
\begin{align*}
B(\theta)=
\begin{pmatrix}
0 & 0 & 1\\
\cos\theta &\sin\theta & 0\\
-\sin\theta &\cos\theta & 0
\end{pmatrix},
\end{align*}
hence we consider the polar frame map given by
\begin{align*}
p(\theta)
=\big(R(\theta,\mathbf{t}(\theta))\hat{z}
\times\mathbf{t}(\theta),
\mathbf{t}(\theta),
R(\theta,\mathbf{t}(\theta))\hat{z}\big),
\end{align*}
where $\hat{z}=(0,0,1)$ is the $z$-versor,
which corresponds to 
rotating the polar frame map of the standard
unit disk:
\begin{align*}
p(\theta)=R(\theta,\mathbf{t}(\theta))
\begin{pmatrix}
\cos\theta & -\sin\theta & 0\\
\sin\theta & \cos\theta & 0\\
0&0&1
\end{pmatrix}.
\end{align*}
Using the compatibility of the product operations between
the $SO(3)$ and $\pi_1(SO(3))$,
we see that
\begin{align*}
[p(\theta)]
&=[R(\theta,\mathbf{t}(\theta))]+[R(\theta,\hat{z})]\\
&=[R(\theta,\mathbf{t}(\theta))]+1.
\end{align*}
To prove that also $R(\theta,\mathbf{t}(\theta))$
belongs to the non-trivial class of $\pi_1(SO(3))$,
we can use its quaternion representation:
\begin{align*}
R(\theta,\mathbf{t}(\theta))
&=(\cos(\theta/2),\sin(\theta/2)\,\mathbf{t}(\theta))\\
&=\cos(\theta/2)-\sin(\theta/2)(-\sin\theta \mathbf{i}
+\cos\theta\mathbf{j}).
\end{align*}
With this representation, since there holds
\begin{align*}
R(0,\mathbf{t}(0))=1 \quad\text{and}\quad R(2\pi,\mathbf{t}(2\pi))=-1,
\end{align*}
the lift of $R(\theta,\mathbf{t}(\theta))$
to the universal cover $S^3$ is not a closed loop,
and this means that the base loop is not nullhomotopic.
We then conclude that
\begin{align*}
[p(\theta)]=1+1=0\quad\text{in }\Z/2\Z.
\end{align*}
\end{example}
This example demonstrates also that a couple $(\gamma, \vec{n}_0)$
needs not to bound an immersion of $\Phi:\overline{B_1(0)^2}\to\R^3$
not even if $\gamma$ is planar and injective.
We now want to prove that the only additional requirement for $\Phi$
to exist is to have a branch point.
The key step to prove it,
obtained through an elementary application of the
so-called \emph{$h$-principle} \cite{MR1909245, MR864505},
is the following lemma.
In what follows, let us denote, for $0\leq r\leq R<\infty$,
\begin{align*}
A[R,r]=\overline{B_R(0)}\setminus B_r(0)
\end{align*}
the annulus of radii $R$ and $r$ centered at 0.

\begin{lemma}
\label{lemma:existence-connecting strip}
Let $\gamma_1,\gamma_2:S^1\to\R^3$
be regular, closed curves of class $C^{k,\alpha}$ for
$k\in\N_{\geq 1}\cup\{\infty\}$ and $\alpha\in(0,1]$,
whose unit tangent vectors we denote
by $\mathbf{t}_1$ and $\mathbf{t}_2$, 
and let $\vec{n}_0,\vec{n}_1:S^1\to S^2$ be unit normal vector fields along 
$\gamma_1$ and $\gamma_2$ respectively of class $C^{k,\alpha}$.
There exists a regular, immersed strip of class $C^{k,\alpha}$
\begin{align*}
\Phi:A[2,1]\to\R^3,
\end{align*}
satisfying
\begin{align}
\label{eq:boundary-prescription-strip}
\Phi|_{\partial B_1(0)}=\gamma_1,\quad \vec{n}_{\Phi}|_{\partial B_1(0)}=\vec{n}_1
\quad\text{and}\quad
\Phi|_{\partial B_1(0)}=\gamma_2,\quad \vec{n}_{\Phi}|_{\partial B_1(0)}=\vec{n}_2
\end{align}
if and only if the maps
\begin{align*}
p_1(x) = (\mathbf{t}_1\times \vec{n}_1, \mathbf{t}_1,\vec{n}_1)(x)
\quad\text{and}\quad
p_2(x) = (\mathbf{t}_2\times \vec{n}_2, \mathbf{t}_2,\vec{n}_2)(x),
\quad x\in S^1,
\end{align*}
are homotopic in $SO(3)$.
\end{lemma}
\begin{bfproof}[Proof of Lemma \ref{lemma:existence-connecting strip}]
The necessity of the condition is clear, we prove the sufficiency.

\emph{Step 1.}
Set, for $\delta>0$, $K_\delta=S^1\times(-\delta,\delta)^2$
and define the following maps for $i=1,2$:
\begin{align*}
\phi_i(\xi,u,v)&=\gamma_i(\xi)+u\,(\mathbf{t}_i\times \vec{n}_i)(\xi)+v\,\vec{n}_i(\xi),
(\xi,u,v)\in \overline{K_\delta}.
\end{align*}
If $\delta$ is chosen small enough, $\phi_1$ and $\phi_2$
define regular immersions (i.e. 
the Jacobian matrix $D\phi(x)$ has rank 3 for every $x\in \overline{K}_\delta$)
of class $C^{k,\alpha}$.
Since $S^1$ and $SO(3)$ are strong deformation retracts of $\overline{K_\delta}$
and $GL^+(3,\R)$ respectively,
an homotopy between $p_1$ and $p_2$ in $SO(4)$ induce an homotopy in $GL^{+}(3,\R)$
between $D\phi_1$ and $D\phi_2$.
Let $(x,t)\mapsto m(x,t)$, $(x,t)\in\overline{K_\delta}\times[0,1]$ be such an homotopy.

\emph{Step 2.}
Let $J^1(\R^3,\R^3)$ be the 1-jet space of maps from $\R^3$ to itself
(see \cite[Chapter 1]{MR1909245}) and let us consider the (local) section
\begin{align*}
F:S^1\times\left[-\delta,\frac{5}{2}\delta\right]\times[-\delta,\delta]\to J^1(\R^3,\R^3),
\quad x\mapsto(x,\phi(x), M(x)),
\end{align*}
where
\begin{align*}
&\phantom{{}={}}
\phi(x)=\phi(\xi,u,v)\\
&=\begin{cases}
\phi_1(\xi,u,v) &\text{if }u\in\left[-\delta,\frac{\delta}{2}\right],\\
\frac{2}{\delta}(\delta-u)\phi_1(\xi,u,v)
+\frac{2}{\delta}\left(\frac{\delta}{2}-u\right)
\phi_2\left(\xi,u-\frac{3}{2}\delta,v\right)
&\text{if }u\in\left[-\frac{\delta}{2},\delta\right],\\
\phi_2\left(\xi,u-\frac{3}{2}\delta, v\right)
&\text{if }u\in\left[\delta,\frac{5}{2}\delta\right],
\end{cases}
\end{align*}
and
\begin{align*}
&\phantom{{}={}}
M(x)
=M(\xi,u,v)\\
&=\begin{cases}
D\phi_1(\xi,u,v) &\text{if }u\in\left[-\delta,\frac{\delta}{2}\right],\\
m\left(\xi,u-3\left(u-\frac{\delta}{2}\right),v,
\frac{2}{\delta}\left(1-\frac{\delta}{2}\right)\right)
&\text{if }u\in\left[-\frac{\delta}{2},\delta\right],\\
D\phi_2\left(\xi,u-\frac{3}{2}\delta,v\right)
&\text{if }u\in\left[\delta,\frac{5}{2}\delta\right],
\end{cases}
\end{align*}

Performing a normalisation of the parameters, we obtain a section
$G:\overline{K_{1/2}}=S^1\times[-1,1]\times[-1,1]\to J^1(\R^3,\R^3)$
which is holonomic in the set
$S^1\times[-1,-1/4]\times[-1,1]
\cup S^1\times[1/4,1]\times[-1,1]$.

\emph{Step 3.}
By the relative version of the Holonomic Approximation theorem
(\cite[Theorems 3.1.1, 3.2.1]{MR1909245})
with
\begin{align*}
A=S^1\times\left[-\frac{3}{4},\frac{3}{4}\right]\times\{0\},
\quad B=S^1\times\left\{-\frac{1}{2},\frac{1}{2}\right\}\times\{0\},
\end{align*}
we may obtain, for every $\varepsilon_1>0$, a diffeomorphism $h:\R^3\to\R^3$
with $\|h-\Id_{\R^3}\|_{C^0(\R^3)}\leq \varepsilon_1$
and satisfying $h\equiv \Id_{\R^3}$ on a open neighbourhood $U$ of $B$,
a holonomic section $\tilde{G}:V\to J^1(\R^3,\R^3)$, where $V\supseteq U$ is an
open neighbourhood of $h(A)$, satisfying $\tilde{G}\equiv G$ on $U$
and $\|\tilde{F}-F\|_{C^0(V)}\leq \varepsilon_1$.
By choosing $\varepsilon_1$ small enough, $\tilde{G}$ is then the 
1-jet extension of an immersion coinciding 
with the one of $\phi_0$ in a open neighbourhood of $S^1\times\{-1/2\}\times\{0\}$
and with the one of $\phi_1$ in a open neighbourhood of $S^1\times\{1/2\}\times\{0\}$;
in particular it is of class $C^{k,\alpha}$ in such neighbourhoods.

Possibly reducing $U$ and $V$, 
the existence of a diffeomorphism $g:K_{1/2}\to V$ that shrinks $K_{1/2}$ into $V$ wile keeping $U$ fixed is ensured.
If we consider the restriction of holonomic section $H=G\circ g:\overline{K_{1/2}}\to J^1(\R^3,\R^3)$ to $ S^1\times [-1/2,1/2]\times\{0\}$,
and denote by $\Psi: S^1\times [-1/2,1/2]\simeq A[2,1]\to\R^3$ the base
map,
we have that $\Psi$ is realises a regular immersion of class $C^1$
which, in a small neighdourhood of the
boundary, containing, say, $S^1\times[1/2+\varepsilon_2,1/2-\varepsilon_2]$
for some small $\varepsilon_2>0$,
and  of class $C^{k,\alpha}$  and satisfies the desired boundary prescriptions
\eqref{eq:boundary-prescription-strip}.

\emph{Step 4.}
To ensure the global  $C^{k,\alpha}$ regularity,
we us use a localised mollification
for $\Psi$:
\begin{align*}
\Phi(x)= \Psi * \rho_{\varepsilon(x)}(x)\quad
x\in S^1\times\left[-\frac{1}{2},\frac{1}{2}\right],
\end{align*}
where $\rho$ is the standard mollification kernel and
$\varepsilon(x)= \varepsilon_3\chi(x)$, where $\chi$ is smooth
cut-off function identically 1 on $ S^1\times[-1/2+\varepsilon_3,1/2-\varepsilon_3]$
and compactly supported on $S^1\times (-1/2,1/2)$
and $\varepsilon_3>0$ has be chosen so small that $\Phi$ has maximum rank.
We conclude that the map $\Phi$ thus defined is the one we have been looking for.
\end{bfproof}

From this lemma we can  prove lemma \ref{lemma:existence-competitors},
in the case $m=3$, that is  we can construct,
given any couple curve-normal vector field $(\gamma,n_0)$ of
class $C^{k,\alpha}$, a possibly branched immersion $\Phi:\overline{B_1(0)}\to\R^3$
of class $C^{k,\alpha}$ assuming such data at the boundary.
\begin{bfproof}[Proof of Lemma \ref{lemma:existence-competitors}]
We consider the the loop
$p(x)=(\mathbf{t}\times n_0,\mathbf{t},n_0)(x)$ in $SO(3)$
induced by $Gamma$ and $\vec{n}_0$.
If it is not nullhomotopic, we can connect, in a $C^1$ way, the couple $(\gamma,n_0)$
and the flat immersion of the disk $z\mapsto(z,0)$ by means of a regular strip of
class $C^{k,\alpha}$.
If it is nullhomotopic instead we can do the same with the 
branched immersion $z\mapsto(z^2,0)$.
if necessary, we smooth out the immersion near the junction
as done in Step 4 of the proof of lemma \ref{lemma:existence-connecting strip};
a final reparametrization of gives then the immersion $\Phi:\overline{B_1(0)^2}\to\R^3$
we have been looking for.

For the last statement,
it is enough to note that when $k=\alpha=1$ or $k\geq 2$, this $\Phi$
satisfies all the assumptions of lemma \ref{lemma:conformal-coordinates}.
This concludes the proof of the lemma.
\end{bfproof}

\begin{remark}[Higher codimension case]
\label{remark:higher-codim}
For general $m\geq 3$,  a regular curve
$\gamma:S^1\to\R^m$ and $(m-2)$-unit normal vector field
$\vec{n}_0:S^1\to \Gr_{m-2}(\R^m)$ along $\gamma$
uniquely determine a loop into the set of couples of ortho-normal
vectors or $\R^m$:
\begin{align*}
x\mapsto (\star(\mathbf{t}\wedge \vec{n}_0), \mathbf{t})(x)
\quad x\in S^1,
\end{align*}
that is to say, into the Stiefel manifold $V_2(\R^m)$
(see e.g. \cite{MR1867354}), which
for the case $m=3$ we could identify with $SO(3)$.
As it is well-known, $\pi_1(V_2(\R^m))=0$ for $m>3$
and hence, the higher-dimensional version of lemma \ref{lemma:existence-connecting strip}
basically says that a regular strip bounding any two couples $(\gamma_1,\vec{n}_1)$
and $(\gamma_2,\vec{n}_2)$ can always be constructed.
As a consequence, with the aid of the Holonomic approximation,
in a similar fashion as the one just described,
we may always find a regular immersion bounding $\gamma$ and $\vec{n}_0$.
This is perhaps not so surprising, as higher codimension gives us
more freedom.
\end{remark}

\subsection{Results on Integrability by Compensation}
\label{app:integrability-compensation}
\begin{theorem}[Wente's Inequality, \cite{MR0336563, MR592104, MR733715}]
\label{thm:Wente}
Let $a,b\in W^{1,2}(D)$ and let $u\in W^{1,1}_0(B_1(0))$ be the solution
to
\begin{align*}
\left\{
\begin{aligned}
-\Delta u &= \langle\nabla^\perp a,\nabla b\rangle &&\text{in } B_1(0),\\
u&=0 &&\text{on }\partial B_1(0).
\end{aligned}
\right.
\end{align*}
Then $u\in C^0(\overline{B_1(0)})\cap W^{1,2}(B_1(0))$ and there is a constant $C>0$ independent of $u,a$ and $b$ 
so that
\begin{align}
\label{eq:Wente-estimate}
\|u\|_{L^\infty(B_1(0))}+\|\nabla u\|_{L^2(B_1(0))}
\leq C\|\nabla a\|_{L^2(B_1(0))}\|\nabla b\|_{L^2(B_1(0))}.
\end{align}
\end{theorem}

For the following lemma we use a result obtained in \cite{MR3764916}.

\begin{lemma}
\label{lemma:estimate-Wente-Neumann}
Let $a,b\in W^{1,2}(B_1(0))$ be so that their traces $a|_{\partial B_1(0)}, b|_{\partial B_1(0)}$
belong to $W^{1,p}(\partial B_1(0))$ for some $p>1$.
Then there exist a  constant $C>0$ independent of $u,a$ and $b$and constant $C(p)>0$ depending only on $p$ so that 
every solution $u\in W^{1,2}(B_1(0))$ of the problem
\begin{align}
\label{eq:wente-neumann-compatible}
\left\{
\begin{aligned}
-\Delta u&=\langle\nabla^\perp a,\nabla b\rangle &&\text{in }B_1(0),\\
\partial_\nu u &=\partial_\tau a\,b &&\text{on }\partial B_1(0),
\end{aligned}
\right.
\end{align}
belongs to $C^0(\overline{B_1(0)})$ and  satisfies the estimate
\begin{align}
\label{eq:estimate-wente-neumann}
\|\nabla u\|_{L^2(B_1(0))}+
\inf_{c\in\R}\|u-c\|_{L^\infty(B_1(0))}
&\leq C\| \nabla a\|_{L^2(B_1(0))}\|b\|_{W^{1,2}(B_1(0))}\\
\notag &\phantom{{}\leq{}}
+C(p)\left(\left\|\partial_\tau a|_{\partial B_1(0)}\right\|_{L^{p}(\partial B_1(0))}\|b|_{\partial B_1(0)}\|_{W^{1,p}(\partial B_1(0))}\right).
\end{align}
\end{lemma}

\begin{remark}
Any weak solution of \eqref{eq:wente-neumann-compatible}
is naturally in $W^{1,2}(B_1(0))$, since it is found by considering
the minimisation of the the functional
\begin{align*}
w\mapsto\int_{B_1(0)}|\nabla w - b\nabla^\perp a|^2\,\dd x,
\end{align*}
among all the functions in $W^{1,2}(B_1(0))$ with a fixed average.
Any two solutions of \eqref{eq:wente-neumann-compatible} differ by a constant.
\end{remark}

\begin{bfproof}[Proof of Lemma \ref{lemma:estimate-Wente-Neumann}]
\emph{Step 1.}
We establish first the $L^{\infty}$ estimate for the analogous problem on the upper half-plane
$\R^2_+=\{(x^1,x^2):x^2>0\}$ for $a,b$ smooth and compactly supported
in some ball $B_R(0)$ of fixed radius $R>0$,
namely
\begin{align}
\label{eq:wente-neumann-compatible-half-plane}
\left\{
\begin{aligned}
-\Delta u&=\langle\nabla^\perp a,\nabla b\rangle &&\text{in }\R^2_+,\\
\partial_\nu u &=\partial_\tau a\,b &&\text{on }\partial \R^2_+.
\end{aligned}
\right.
\end{align}
Because of the ``compatible'' boundary condition (using a terminology from 
\cite{MR3764916}), 
the only solution of \eqref{eq:wente-neumann-compatible-half-plane}
belonging to $W^{1,2}(\R^2_+)$
is given by the representation formula
\begin{align}
\label{eq:rep-formula-neumann-half-plane}
u(x)
&=-\int_{\R^2_+}\langle
\nabla_y\mathcal{G}_{\R^2_+}(x,y),b(y)\nabla^\perp a(y)\rangle\,\dd y,
\quad x\in \R^2_+,
\end{align}
where $\mathcal{G}_{\R^2_+}$ is the Green function for the Neumann
problem in the half-plane given by
\begin{align*}
\mathcal{G}_{\R^2_+}(x,y)=-\frac{1}{2\pi}
\left(\log|x-y| +\log|\overline{x}-y|\right),
\quad (x,y)\in \R^2_+\times\R^2_+,
\end{align*}
(here $\overline{x}$ denotes the complex conjugate of $x$,
that is if $x=(x^1,x^2)$ then $\overline{x}=(x^1,-x^2)$).
Note that such formula implies that $u\in C^{\infty}(\overline{\R^2_+})$
\footnote{this need \emph{not} to be the case,
not even for smooth $a$ and $b$, for different Neumann
(for example homogeneous) boundary conditions.}.
As proved in \cite[\S 3.1]{MR3764916}, the trace of the solution
$u$ of \eqref{eq:wente-neumann-compatible-half-plane} given
by \eqref{eq:rep-formula-neumann-half-plane} can be written as
\begin{align*}
u(x^1,0)=A(a,b)(x^1)+B(a,b)(x^2),
\quad x^1\in\R\simeq\partial \R^2_+,
\end{align*}
where the function $A(a,b)$ satisfies
\begin{align}
\label{eq:estimate-A}
\|A(a,b)\|_{L^\infty(\R)}
\leq C
\|\nabla a\|_{L^2(\R^2_+)}\|\nabla b\|_{L^2(\R^2_+)},
\end{align}
and the function $B(a,b)$ is given by
\begin{align}
\notag B(a,b)(x^1)&=
\frac{1}{2\pi}\pv\int_{\R}
\frac{1}{t}\left(a(x^1+t,0)+a(x^1-t,0)\right)
\left(b(x^1-t,0)-b(x^1+t,0)\right)\,\dd t\\
\label{eq:rep-form-B}&=\frac{1}{\pi}\pv\int_{\R}
\frac{1}{t}a(x^1+t,0)b(x^1-t,0)\,\dd t
+\frac{1}{\pi}\pv\int_{\R}
\frac{1}{t}a(x^1-t,0)b(x^1-t,0)\,\dd t.
\end{align}
We recognise that the second term in the above expression is
the Hilbert transform $H(ab)$ of $a b$ (see e.g.\cite{MR0290095} or \cite[\S 5.1]{MR3243734}).
We recall that the Hilbert transform $H$ is a tempered distribution
which  maps $L^p(\R)$ into itself for $1<p<\infty$;
moreover it is a convolution-type distribution and consequently
it commutes with derivatives.
This implies that, for a constant $C(p)>0$ depending only on $p$,
there holds
\begin{align*}
\|H(f)\|_{W^{1,p}(\R)}\leq C(p) \|f\|_{W^{1,p}(\R)},
\quad\forall\, f\in W^{1,p}(\R).
\end{align*}
On the other hand, by the Sobolev embedding we have
$W^{1,p}(\R)\hookrightarrow C^{0,\gamma}(\R)$ 
continuously with $\gamma=1-1/p$.
This implies in particular that $W^{1,p}(\R)$ is an algebra,
that is
\begin{align*}
\|f\,g\|_{W^{1,p}(\R)}\leq C(p)\|f\|_{W^{1,p}(\R)}\|g\|_{W^{1,p}(\R)}
\quad\forall\, f, g\in W^{1,p}(\R),
\end{align*}
where $C(p)>0$ is a constant depending only on $p$.
Using these facts we may estimate the second term in
expression \eqref{eq:rep-form-B}
uniformly in $x^1$ as 
\begin{align}
\label{eq:estimate-B-2nd}
\left|\frac{1}{\pi}\pv\int_{\R}
\frac{1}{t}a(x^1-t,0)b(x^1-t,0)\,\dd t\right|
\leq\|H(a b)\|_{C^{0,\alpha}(\R)}
\leq C(p)\|a\|_{W^{1,p}(\R)}\|b\|_{W^{1,p}(\R)}.
\end{align}
The first term in expression \eqref{eq:rep-form-B}
may be similarly estimated as follows:
fix $x^1$ and define the function $h_{x^1}(t)=a(x^1+t,0)b(x^1-t,0)$.
We then have:
\begin{align}
\label{eq:estimate-B-1st}
\bigg|\frac{1}{\pi}
\pv\int_{\R}
\frac{1}{t}a(x^1+t,0)b(x^1-t,0)\,\dd t\bigg|
&=|H(h_{x^1})(0)|\\
\notag&\leq \|H(h_{x^1})\|_{C^{0,\gamma}}\\
\notag&\leq C(p)\|h_{x^1}\|_{W^{1,p}(\R)}\\
\notag&\leq C(p)\|a(x^1+\cdot)\|_{W^{1,p}(\R)}\|a(x^1-\cdot)\|_{W^{1,p}(\R)}\\
\notag&=C(p)\|a\|_{W^{1,p}(\R)}\|b\|_{W^{1,p}(\R)}.
\end{align}
We may then join estimates \eqref{eq:estimate-A}, \eqref{eq:estimate-B-2nd}
and \eqref{eq:estimate-B-1st} to deduce the following bound on the trace of $u$:
\begin{align}
\label{eq:estimate-trace-half-plane}
\|u|_{\partial \R^2_+}\|_{L^\infty(\partial \R^2_+)}
&\leq 
C\|\nabla a\|_{L^2(\R^2_+)}\|\nabla b\|_{L^2(\R^2_+)}\\
\notag&\phantom{{}\leq{}}+C(p)\left(\|a|_{\partial \R^2_+}\|_{W^{1,p}(\partial \R^2_+)}
\|b|_{\partial \R^2_+}\|_{W^{1,p}(\partial \R^2_+)}\right).
\end{align}

\emph{Step 2.}
We now come to the problem on the disk, that is, we consider a solution
$u$ of the problem \eqref{eq:wente-neumann-compatible}.
By an approximation argument we may suppose that $a,b\in C^{\infty}(\overline {B_1(0)})$
and, as in step 1, because of the ``compatible'' boundary condition, the solution is
represented by
\begin{align*}
u(x)
=-\int_{B_1(0)}
\langle\nabla_y\mathcal{G}(x,y),b(y)\nabla^\perp a(y)\rangle
\,\dd y\quad x\in B_1(0),
\end{align*}
where $\mathcal{G}$ is the Green function for the Neumann problem
(see \eqref{eq:Green-for-Neumann}), and in particular we have
$u\in C^\infty(\overline{B_1(0)})$.
Let us fix a function $\chi\in C^\infty(\overline{B_1(0)})$
so that $\chi=1$ in $\overline{B_1(0)}\cap \{x^1>0\}$ and
whose support is contained in $\overline{B_1(0)}\cap \{x^1\geq-1/2\}$.
We then write $u=u_1+u_2$, where
\footnote{note that $u_1$ and $u_2$ are determined up to a constant.}
\begin{align*}
&\left\{
\begin{aligned}
-\Delta u_1&=\langle\nabla^\perp a,\nabla (\chi b)\rangle &&\text{in }B_1(0),\\
\partial_\nu u_1 &=[\partial_\tau a]\,(\chi b) &&\text{on }\partial B_1(0),
\end{aligned}
\right.
&\text{and}&
&\left\{
\begin{aligned}
-\Delta u_2&=\langle\nabla^\perp a),\nabla ((1-\chi) b)\rangle &&\text{in }B_1(0),\\
\partial_\nu u_2 &=[\partial_\tau a]\,((1-\chi) b) &&\text{on }\partial B_1(0).
\end{aligned}
\right.
\end{align*}
If we fix
a function $\psi\in C^{\infty}(\overline{B_1(0)})$ so that
$\psi =1$ in $\overline{B_1(0)}\cap \{x^1\geq -1/2\}$ and whose support
is contained in $\overline{B_1(0)}\cap \{x^1\geq -3/4\}$,
we may write
\begin{align*}
\left\{
\begin{aligned}
-\Delta u_1&=\langle\nabla^\perp(\psi a),\nabla (\chi b)\rangle &&\text{in }B_1(0),\\
\partial_\nu u_1 &=[\partial_\tau (\psi a)]\,(\chi b) &&\text{on }\partial B_1(0),
\end{aligned}
\right.
\end{align*}
If $\varphi:\R^2_+\overset{\sim}{\to} B_1(0)$ denotes the fractional linear transformation 
sending $0$ to $1$, $i$ to $0$ and $-1$ to $\infty$,
given in complex coordinates by $\varphi(z)=-(z-i)/(z+i)$,
 $(\psi a)\circ\varphi$ and $(\chi b)\circ\varphi$ have compact support
and with estimates
\begin{align*}
\|\nabla((\psi a)\circ\varphi)\|_{L^2(\R^2_+)}&\leq C\| a\|_{W^{1,2}(B_1(0))},&
\|((\psi a)\circ\varphi)|_{\partial\R^2_+}\|_{W^{1,p}(\partial\R^2_+)},
\leq C(p)\|a|_{\partial B_1(0)}\|_{W^{1,p}(\partial B_1(0))},\\
\|\nabla((\chi b)\circ\varphi)\|_{L^2(\R^2_+)}&\leq C\|b\|_{W^{1,2}(B_1(0))},&
\|((\chi a)\circ\varphi)|_{\partial\R^2_+}\|_{W^{1,p}(\partial\R^2_+)},
\leq C(p)\|a|_{\partial B_1(0)}\|_{W^{1,p}(\partial B_1(0))}.
\end{align*}
The invariance of Dirichlet energy and of the problem \eqref{eq:wente-neumann-compatible}
implies that, if
we compose $u_1$ $\psi a$ and $\chi b$ with $\varphi$,
then there must be a constant $c_1\in\R$ so that
\begin{align*}
u_1\circ\varphi(x) - c_1
&=-\int_{\R^2_+}\langle\nabla_y\mathcal{G}_{\R^2_+}(x,y),
((\chi b)\circ\varphi)(y)\nabla^\perp((\psi a)\circ\varphi)\rangle(y)\,\dd y
\quad x\in\R^2_+.
\end{align*}
We then deduce from estimate \eqref{eq:estimate-trace-half-plane},
that the trace $u_1|_{\partial B_1(0)}$ can be estimated as
\begin{align}
\label{eq:estimate-trace-u1}
\|u_1|_{\partial B_1(0)}-c_1\|_{L^\infty(\partial B_1(0))}
&=\|(u_1\circ\varphi)|_{\partial B_1(0)}-c_1\|_{L^\infty(\partial \R^2_+)}\\
\notag&\leq C\|a\|_{W^{1,2}(B_1(0))}\|b\|_{W^{1,2}(B_1(0))}\\
\notag&\phantom{{}\leq{}}+C(p)\left(
\|a|_{\partial_{B_1(0)}}\|_{W^{1,p}(\partial B_1(0))}
\|b|_{\partial_{B_1(0)}}\|_{W^{1,p}(\partial B_1(0))}
\right).
\end{align}
In the same way, 
by working with $1-\chi$, $\psi(-z)$ and $\varphi(-z)$
ins place of $\chi$, $\psi$ and $\phi$ respectively,
we obtain an analogous estimate for $u_2|_{\partial_{B_1(0)}}-c_2$ for some $c_2\in\R$.
Since the problem \eqref{eq:wente-neumann-compatible}
is invariant under translations with respect to $a$, that is, for any $c\in\R$
we have
$\langle\nabla^\perp(a-c),\nabla b\rangle
=\langle\nabla^\perp a,\nabla b\rangle$ in $B_1(0)$ and
$\partial_\tau(a-c)\, b
=\partial_\tau a\, b$ on $\partial B_1(0)$,
we conclude that
\begin{align}
\label{eq:estimate-trace-u}
\inf_{c\in\R}\|u|_{\partial B_1(0)}-c\|_{L^\infty(\partial B_1(0))}
&\leq C\inf_{c\in\R}\|a-c\|_{W^{1,2}(B_1(0))}\|b\|_{W^{1,2}(B_1(0))}\\
&\notag
\phantom{{}\leq{}}
+C(p)\left(
\inf_{c\in\R}\|a|_{\partial B_1(0)}-c\|_{W^{1,p}(\partial B_1(0))}
\|b|_{\partial B_1(0)}\|_{W^{1,p}(\partial B_1(0))}
\right).
\end{align}
holds. By the maximum principle for harmonic functions,
Poincar\'e's inequality 
and Wente's inequality \eqref{eq:Wente-estimate} we 
finally deduce the the validity of 
 {the $L^\infty$-part of} estimate
of 
\eqref{eq:estimate-wente-neumann}.
{
The estimate on the Dirichlet energy of $u$ follows
integrating by parts
(recall that, by approximation, we are working with
$u\in C^{\infty}(\overline{B_1(0)})$), namely:
\begin{align*}
\int_{B_1(0)}|\nabla u|^2\,\dd x
=-\int_{B_1(0)}(u-c)\Delta u\,\dd x
+\int_{B_1(0)}(u-c)\partial_\nu u\,\dd \mathcal{H}^1
\quad \text{for every} c\in\R,
\end{align*}
hence estimating the terms on the right-hand-side with
H\"older and Sobolev inequalities:
\begin{align*}
\bigg|
\int_{B_1(0)}(u-c)\Delta u\,\dd x
\bigg|
&\leq \|u-c\|_{L^\infty(B_1(0))}\|\nabla a\|_{L^2(B_1(0))}\|\nabla b\|_{L^2(B_1(0))},\\
\bigg|
\int_{\partial B_1(0)}(u-c)\partial_\nu u\,\dd \hau^1
\bigg|
&\leq \|u-c\|_{L^\infty(\partial B_1(0))}
\|\partial_\tau a\|_{L^1(\partial B_1(0))}
\| b\|_{L^\infty(\partial B_1(0))}\\
&\leq
C(p)\|u-c\|_{L^\infty( B_1(0))}
\|\partial_\tau a\|_{L^p(\partial B_1(0))}
\| b\|_{W^{1,p}(\partial B_1(0))}.
\end{align*}
and using finally the $L^\infty$-estimate above obtained.}
This concludes the proof of the lemma.
\end{bfproof}

With a standard argument,  we may localise the
above result as follows.

\begin{lemma}
\label{lemma:estimate-localised-Wente-Neumann}
Let $f\in L^1(B_1(0))$, $g\in L^1(\partial B_1(0))$,
$a_1,\ldots a_N$ be points in $B_1(0)$
and $\alpha_1,\ldots,\alpha_N$ be real numbers 
satisfying $\int_{B_1(0)} f\,\dd x +\sum_i \alpha_i =- \int_{\partial_{B_1(0)}}g\,\dd\sigma$.
Let $u\in W^{1,1}(B_1(0))$ be a weak solution to the problem
\begin{align}
\left\{
\begin{aligned}
-\Delta u&=f+\sum_{i=1}^N\alpha_i\delta_{a_i} &&\text{in }B_1(0),\\
\partial_\nu u &=g &&\text{on }\partial B_1(0).
\end{aligned}
\right.
\end{align}
Assume that, for a given $x_0\in\partial B_1(0)$ and $0<r<1$,
$B_1(0)\cap B_r(x_0)$ contains none of the $a_i$'s and there holds
\begin{align*}
f&=\langle\nabla^\perp a,\nabla b\rangle
\text{ in }B_1(0)\cap B_r(x_0),\\
g&=\partial_\tau a\,b\text{ on }\partial B_1(0)\cap B_r(x_0),
\end{align*}
for some $a,b\in W^{1,2}(B_1(0)\cap B_r(x_0))$ so that the traces on $\partial B_1(0)$
$a|_{\partial B_1(0)}, b|_{\partial B_1(0)}$
belong to $W^{1,p}(\partial B_1(0)\cap B_r(x_0))$ for some $p>1$.
Then $u\in C^0(\overline{B_1(0)\cap B_{r/2}(x_0)})
\cap W^{1,2}(B_1(0)\cap B_{r/2}(x_0))$ and there exists
 constants $C(r)>0$  $C(r,p)>0$ so that
\begin{align}
&\phantom{{}\leq{}}
\inf_{c\in\R}\|u-c\|_{L^{\infty}(B_1(0)\cap B_{r/2}(x_0))}\\
\notag&\leq C(r)
\bigg(\|f\|_{L^1(B_1(0))} + \|g\|_{L^1(\partial B_1(0))} + \sum_{i=1}^N |\alpha_i|\bigg)
+C(r)
\|\nabla a\|_{L^{2}(B_1(0)\cap B_r(x_0))}\|b\|_{W^{1,2}(B_1(0)\cap B_r(x_0))}\\
\notag&\phantom{{}\leq{}}+C(r,p)
\|\partial_\tau a|_{\partial B_1(0)}\|_{L^{p}(\partial B_1(0)\cap B_r(x_0))}
\|b|_{\partial B_1(0)}\|_{W^{1,p}(\partial B_1(0)\cap B_r(x_0))} .
\end{align}
\end{lemma}

\begin{bfproof}[Proof of Lemma \ref{lemma:estimate-localised-Wente-Neumann}]
We let $\chi$ be a function in $C^\infty(\overline{B_1(0)})$
so that $\chi =1$ in $B_1(0)\cap B_{3r/4}(x_0)$
and whose support is contained in $B_1(0)\cap B_{7r/8}(x_0)$,
and we let $\tilde{a}$ be an extension of $a$ to $B_1(0)$,
obtained through a suitable Moebius tranformation of $B_1(0)$
so that $\|\nabla\tilde{a}\|_{L^2(B_1(0))}\leq C
\|\nabla a\|_{L^2(B_1(0)\cap B_r(x_0))}$
and 
$\|\partial_\tau\tilde{a}\|_{L^p(\partial B_1(0))}\leq C(p)
\|\partial_\tau a\|_{L^p(\partial B_1(0)\cap B_r(x_0))}$
for constants $C, C(p)>0$.
Up to a constant, we may write $u=u_1+u_2$, with
\begin{align*}
&\left\{
\begin{aligned}
-\Delta u_1&=\langle\nabla^\perp \tilde{a},\nabla (\chi b)\rangle &&\text{in }B_1(0),\\
\partial_\nu u_1 &=[\partial_\tau \tilde{a}]\,(\chi b) &&\text{on }\partial B_1(0),
\end{aligned}
\right.
\end{align*}
and
\begin{align*}
&\left\{
\begin{aligned}
-\Delta u_2&=\langle\nabla^\perp a,\nabla ((1-\chi) b)\rangle
+\sum_{i}\alpha_i\delta_{a_i} &&\text{in }B_1(0),\\
\partial_\nu u_2 &=[\partial_\tau a]\,((1-\chi) b) &&\text{on }\partial B_1(0),
\end{aligned}
\right.
\end{align*}
 with the convention that
\begin{align*}
\langle\nabla^\perp a,\nabla ((1-\chi) b)\rangle = f
\text{ in }B_1(0)\setminus B_{7r/8}(x_0)
\quad\text{and}\quad[\partial_\tau a]\,((1-\chi) b) = g
\text{ on }\partial B_1(0)\setminus B_{7r/8}(x_0).
\end{align*}
From lemma \ref{lemma:estimate-Wente-Neumann},
we deduce that for some constant $c_1\in\R$ there holds
\begin{align}
\label{eq:estimate-localisation-u1}
&\phantom{{}={}}
\|u_1-c_1\|_{L^\infty(B_1(0))}
{+\|\nabla u_1\|_{L^2(B_1(0))}}\\
&\leq\notag
C\|\nabla \tilde{a}\|_{W^{1,2}(B_1(0))}\|\chi b\|_{W^{1,2}(B_1(0))}\\
\notag&\phantom{{}\leq{}}
+C(p)\left(\left\|\partial_\tau \tilde{a}|_{\partial B_1(0)}\right\|_{W^{1,p}(\partial B_1(0))}
\|(\chi b)|_{\partial B_1(0)}\|_{W^{1,p}(\partial B_1(0))}\right).
\end{align}

To estimate $u_2$, we use the representation formula:
\begin{align*}
u_2(x) - \overline{ u}_2
&=\int_{B_1(0)} \mathcal{G}(x,y)
\langle\nabla^\perp a,
\nabla((1-\chi)b)\rangle(y)\,\dd y\\
&\phantom{{}+{}}+\int_{\partial B_1(0)} \mathcal{G}(x,y)
\langle\partial_\tau a,(1-\chi)b\rangle\,\dd\hau^1(y)
+\sum_{i=1}^{N}\alpha_i\mathcal{G}(x,a_i),
\end{align*}
since none of the $a_i$'s is in $B_1(0)\cap B_r(x_0)$
and $1-\chi$ vanishes $B_1(0)\cap B_{3r/4}(x_0)$,
we may estimate on $B_1(0)\cap B_{r/2}(x_0)$:
\begin{align}
\label{eq:estimate-localisation-u2}
&\phantom{{}\leq{}}
\|u_2 - \overline{u}_2\|_{L^\infty(B_1(0)\cap B_{r/2}(x_0))}\\
&\leq C(r)\bigg(\|\langle\nabla^\perp a,
\nabla((1-\chi)b)\rangle\|_{L^1(B_1(0))}
\notag+\|[\partial_\tau a](1-\chi)b\rangle\|_{L^1(\partial B_1(0))}
+\sum_{i=1}^N|\alpha_i|\bigg).
\end{align}
Since  $\chi$ can always be chosen so that 
$\|\nabla \chi\|_{L^\infty(B_1(0))}\leq C/r$,
by joining estimates \eqref{eq:estimate-localisation-u1} and
 \eqref{eq:estimate-localisation-u2} we reach the conclusion.
\end{bfproof}

\subsection{Facts About Moving Frames}
\label{app:frames}
A \emph{moving frame}, or simply a \emph{frame}, from a domain $\Omega\subseteq\R^2$ into $\R^m$
is a map $\vec{f}=(\vec{f}_1,\vec{f}_2):\Omega\to\R^m\times\R^m$ so that $\langle\vec{f}_{i},\vec{f}_j\rangle=\delta_{ij}$.
If $\vec{f}$ and $\vec{g}$ are two frames from $\Omega$ into $\R^m$, a function $\phi:\Omega\to\R$
defines  a way to pass from $\vec{f}$ to $\vec{g}$, i.e. a \emph{gauge transformation}:
\begin{align*}
\left\{
\begin{aligned}
\vec{g}_1(x) &=\vec{f}_1(x)\cos\phi(x) + \vec{f}_2(x)\sin\phi(x),\\
\vec{g}_2(x) &=-\vec{f}_1(x)\sin\phi(x) + \vec{f}_2(x)\cos\phi(x),\\
\end{aligned}
\right.
\end{align*}
which can be written using the complex notation as
\begin{align*}
\vec{g}_1 + i\vec{g}_2 =\e^{-i\phi}(\vec{f}_1 +i\vec{f}_2)
\quad\text{in }\Omega.
\end{align*}
Differentiating this relation, we deduce the following \emph{Change of gauge formula}:
\begin{align}
\label{eq:change-of-gauge}
\langle\nabla\vec{g}_1,\vec{g}_2\rangle
=\nabla\phi +\langle\nabla\vec{f}_1,\vec{f}_2\rangle
\quad\text{in }\Omega.
\end{align}
For a map $\vec{n}:B_1(0)\to\Gr_{m-2}(\R^m)$
expressed as $\vec{n}=\vec{n}_1\wedge\ldots\wedge\vec{n}_{m-2}$
for some $(m-2)$-tuple of ortho-normal sections: $\langle\vec{n}_i,\vec{n}_j\rangle=\delta_{ij}$,
we say that the frame $(\vec{f}_1,\vec{f}_2)$
\emph{is a lift for $\vec{n}$} if:
\begin{align*}
\vec{n} = \star(\vec{f}_1\wedge \vec{f}_2)
\quad\text{in }B_1(0),
\end{align*}
where
$\star$ denotes the Euclidean Hodge operator in $\R^m$
transforming $2$ vectors into $m-2$ vectors and vice-versa
(see e.g. \cite{MR615912}).
By orthonormality, there holds:
\begin{align*}
\nabla\vec{f}_1
&=\langle\nabla\vec{f}_1,\vec{f}_2\rangle\vec{f}_2
+\sum_{i=1}^{m-2}\langle\nabla\vec{f}_1,\vec{n}_i\rangle \vec{n}_i
=\langle\nabla\vec{f}_1,\vec{f}_2\rangle\vec{f}_2
-\sum_{i=1}^{m-2}\langle\vec{f}_1,\nabla \vec{n}_i\rangle \vec{n}_i,\\
\nabla\vec{f}_2
&=\langle\nabla\vec{f}_2,\vec{f}_1\rangle\vec{f}_1
+\sum_{i=1}^{m-2}\langle\nabla\vec{f}_2,\vec{n}_i\rangle \vec{n}_i
=\langle\nabla\vec{f}_2,\vec{f}_1\rangle\vec{f}_1
-\sum_{i=1}^{m-2}\langle\vec{f}_2,\nabla \vec{n}_i\rangle \vec{n}_i,\\
\nabla \vec{n}_i &= \langle\nabla \vec{n}_i,\vec{f}_1\rangle\vec{f}_1
+\langle\nabla \vec{n}_i,\vec{f}_2\rangle\vec{f}_2\quad (i=1,\ldots,m-2),
\end{align*}
in particular:
\begin{align}
\label{eq:decomposition-frame-energy}
|\nabla\vec{f}_1|^2+|\nabla\vec{f}_2|^2
&=2|\langle\nabla \vec{f}_1,\vec{f}_2\rangle|^2 + |\nabla \vec{n}|^2.
\end{align}
When $\vec{n}=\vec{n}_\Phi$ is the Gauss map of an immersion $\Phi:B_1(0)\to\R^m$
and $(\vec{f},\vec{f}_2)$ is a positively oriented ortho-normal basis of the tangent space.
When $\vec{f}$ and $\vec{n}$ correspond to a orthonormal frame and Gauss map of an immersion $\Phi$,
an elementary computation reveals that
\begin{align}
\label{eq:Gauss-curvature-with-frames}
K_\Phi\,\dd vol_\Phi
=\langle\nabla^\perp\vec{f}_1,\vec{n}\rangle
\langle\nabla\vec{f}_2,\vec{n}\rangle\,\dd x
=\langle\nabla^\perp \vec{f}_1,\nabla\vec{f}_2\rangle\,\dd x,
\end{align}
where $K_\Phi$ is the Gauss curvature of $\Phi$ and $\dd\mu_\Phi$ its area element.
This equation has two important consequences:
first, it reveals that the the left-hand-side is a sum of Jacobians.
Second, if $\vec{n}$ is sufficiently regular (namely, that it admits a lifting frame $\vec{f}$:
see lemmas \ref{lemma:Helein} and \ref{lemma:frame-prescribed-boundary} below)
then it has a meaning also when $\Phi$ is singular, for example at a branch point,
when a tangent plane is not defined.
Finally, it is useful to note that
\begin{align}
\label{eq:pointwise-estimate-Gauss-second-fund}
|\langle\nabla^\perp \vec{f}_1,\nabla \vec{f}_2\rangle|
&\leq\frac{1}{2}\left(
|\langle\nabla^\perp\vec{f}_1,\vec{n}\rangle|^2
+|\langle\nabla\vec{f}_2,\vec{n}\rangle|^2\right)
=\frac{|\nabla \vec{n}|^2}{2}.
\end{align}

Recall H\'elein's lifting lemma
(\cite[Lemma 5.1.4]{MR1913803}, see also \cite{MR3016499}).
\begin{lemma}
\label{lemma:Helein}
There is an $\varepsilon_0>0$ so that, for every $0<\varepsilon<\varepsilon_0$,
there exists a constant $C>0$ independent of $\varepsilon$ with the following property.
If a map
$\vec{n}\in W^{1,2}(B_1(0),\Gr_2(\R^m))$ satisfies
\begin{align*}
\|\nabla \vec{n}\|^2_{L^2(B_1(0))}\leq\varepsilon
\end{align*} 
for some $0<\varepsilon<\varepsilon_0$,
then there exist an orthonormal frame 
$\vec{f}=(\vec{f}_1,\vec{f}_2)\in W^{1,2}(B_1(0),\R^m\times\R^m)$
which is a positive orthonormal basis for $\vec{n}$, i.e.:
\begin{align*}
\vec{n} = \star(\vec{f}_1\wedge \vec{f}_2)\quad\text{in }B_1(0),
\end{align*}
that satisfies the following Coulomb condition:
\begin{align}
\label{eq:Helein-Coulomb}
\begin{cases}
\divop(\langle\nabla \vec{f}_1,\vec{f}_2\rangle)&=0 \text{ in }B_1(0),\\
\langle\partial_\nu \vec{f}_1,\vec{f}_2\rangle &= 0 \text{ on }\partial B_1(0),
\end{cases}
\end{align}
and whose energy is controlled as follows:
\begin{align}
\label{eq:Helein-Coulomb-energy}
\|\nabla \vec{f}\|^2_{L^2(B_1(0))}
\leq C\|\nabla \vec{n}\|^2_{L^2(B_1(0))}.
\end{align}
\end{lemma}
The proof consists on a minimizing procedure and 
the $L^2$-bound relies on Wente's lemma.
The following variant concerns the existence of
a energy-controlled lift with a prescribed boundary value.

\begin{lemma}
\label{lemma:frame-prescribed-boundary}
There is an $\varepsilon_0>0$ so that, for every $0<\varepsilon<\varepsilon_0$,
there exists a constant $C>0$ independent of $\varepsilon$ with the following property.
For any map $\vec{n}\in W^{1,2}(B_1(0),\Gr_{m-2}(\R^m))$
and any ortho-normal frame
$\vec{e}=(\vec{e}_1,\vec{e}_2)\in H^{1/2}(\partial B_1(0),\R^m\times\R^m)$
lifting $\vec{n}$:
\begin{align*}
\vec{n}=\star(\vec{e}_1\wedge\vec{e}_2)
\quad\text{on }\partial B_1(0),
\end{align*}
and satisfying the estimate
\begin{align}
\label{eq:boundary-frame-smallness-assumption}
\|\nabla\vec{n}\|_{L^2(B_1(0))}
+[\vec{e}\,]_{W^{1/2,2}(\partial B_1(0))}\leq \varepsilon,
\end{align}
there exists an ortho-normal frame 
$\vec{g}=(\vec{g}_1,\vec{g}_2)\in W^{1,2}(B_1(0),\R^m\times\R^m)$ lifting $\vec{n}$:
\begin{align}
\label{eq:lifting-g}
\vec{n}=\star(\vec{g}_1\wedge\vec{g}_2)
\quad\text{in }B_1(0),
\end{align}
whose trace on $\partial_{B_1(0)}$ coincides with $\vec{e}$, satisfying
 the Coulomb condition
\begin{align}
\label{eq:coulomb-g}
\divop\left(\langle\nabla\vec{g_1},g_2\rangle\right)=0
\quad\text{in }B_1(0),
\end{align}
and the estimate
\begin{align}
\|\nabla\vec{g}\|_{L^2(B_1(0))}\leq C\left(
\|\nabla\vec{n}\|_{L^2(B_1(0))}+[\vec{e}\,]_{W^{1/2,2}(\partial B_1(0))}\right).
\end{align}
\end{lemma}
\begin{bfproof}[Proof of Lemma \ref{lemma:frame-prescribed-boundary}]
Let us start by fixing $\varepsilon_0<2\pi$, so that 
by lemma \eqref{lemma:Helein} we deduce the existence of a 
Coulomb frame $\vec{f}$ on $B_1(0)$ satisfying
\begin{align}
\label{eq:estimate-small-Helein-frame}
\|\nabla\vec{f}\|_{L^2(B_1(0))}\leq \sqrt{2}\|\nabla\vec{n}\|_{L^2(B_1(0))}.
\end{align}
We now want to identify the angle $\alpha_0:\partial B_1(0)\to\R$ which rotates $f$ 
to $\vec{e}$, implicitly defined in complex notation by
$\vec{e}_1+i\vec{e}_2 = \e^{-i\alpha}(\vec{f}_1+i\vec{f}_2)$.
To this aim, let us define the $S^1$-valued function:
\begin{align*}
u=\langle\vec{e}_1,\vec{f}_1\rangle
-i\langle\vec{e}_2,\vec{f}_1\rangle
\quad\text{on }\partial B_1(0),
\end{align*}
and note that it belongs to $W^{1/2,2}(\partial B_1(0),S^1)$,
satisfying
\footnote{
Recall that if $\Omega$ is a domain of $\R$ or of $S^1$, the space 
$(W^{1/2,2}\cap L^\infty)(\Omega)$
is an algebra with:
\begin{align*}
[ab]_{W^{1/2,2}(\Omega)}^2\leq 2
(\|a\|_{L^\infty(\Omega)}^2[b]_{W^{1/2,2}(\Omega)}^2
+\|b\|_{L^\infty(\Omega)}^2[a]_{W^{1/2,2}(\Omega)}^2).
\end{align*}}
because of \eqref{eq:boundary-frame-smallness-assumption}
and \eqref{eq:estimate-small-Helein-frame} the estimate
\begin{align}
\label{eq:estimate-implicit-lift}
[u]_{W^{1/2,2}(\partial_{B_1(0)})}
&\leq 2\left(
[\vec{e}\,]_{W^{1/2,2}(\partial B_1(0))}+
[\vec{f}\,]_{W^{1/2,2}(\partial B_1(0))}
\right)
\leq C\varepsilon,
\end{align} 
for some constant $C>0$.
By choosing $\varepsilon_0>0$ sufficiently small,
we may then invoke theorem \cite[Theorem 1]{MR3693661}
and deduce the existence of a function $U\in W^{1,2}(B_1(0),S^1)$
whose trace on $\partial_{B_1(0)}$ coincides with
$u$
and satisfying the estimate
\begin{align}
\label{eq:estimate-controlled-extension}
\|\nabla U\|_{L^2(B_1(0))}\leq C[u]_{W^{1/2,2}(\partial_{B_1(0)})},
\end{align}
for some  $C>0$. 
For such $U$, we may now deduce the existence
of a lift $\alpha\in W^{1,2}(B_1(0))$ (see e.g. \cite[Theorem 3]{MR1771523}),
that is,
\begin{align*}
U(x)=\e^{i\alpha(x)}\quad x\in B_1(0),
\end{align*}
satisfying the point-wise almost everywhere estimate $|\nabla U|=|\nabla \alpha|$.
 if $\widetilde{\alpha}$ denotes the harmonic extension of
the trace of $\alpha$ on $\partial B_1(0)$, 
the Dirichlet principle together with the inequalities
\eqref{eq:boundary-frame-smallness-assumption} and \eqref{eq:estimate-implicit-lift}
imply
\begin{align}
\label{eq:estimate-harmonic-angle}
\|\nabla\widetilde{\alpha}\|_{L^2(B_1(0))}
\leq C\left(
\|\nabla\vec{n}\|_{L^2(B_1(0))}+[\vec{e}\,]_{W^{1/2,2}(\partial B_1(0))}\right).
\end{align}
Finally, we set:
\begin{align*}
\vec{g}_1 +i\vec{g}_2
=\e^{-i\widetilde{\alpha}}(\vec{f}_1+i\vec{f}_2)
\quad \text{in }B_1(0).
\end{align*}
By construction, the frame $\vec{g}$
has trace equal to $\vec{e}$ on $\partial B_1(0)$,
and satisfies conditions \eqref{eq:lifting-g}, \eqref{eq:coulomb-g}
and from formula \eqref{eq:decomposition-frame-energy}, we see that almost everywhere in 
$B_1(0)$
the relation
\begin{align*}
|\nabla\vec{g}|^2 &= 2|\langle\nabla\vec{g}_1,\vec{g}_2\rangle|^2+|\nabla\vec{n}|^2\\
&=2|\langle\nabla\vec{f}_1,\vec{f}_2\rangle+\nabla\widetilde{\alpha}|^2+|\nabla\vec{n}|^2\\
&\leq4 |\nabla\vec{f}|^2+|\nabla\widetilde\alpha|^2
\end{align*}
holds,
from which we deduce, thanks to inequalities \eqref{eq:estimate-small-Helein-frame}
and \eqref{eq:estimate-harmonic-angle},
the validity of \eqref{eq:estimate-controlled-extension}.
This concludes the proof of the lemma.
\end{bfproof}
A localised version of the above lemma is as follows.
\begin{lemma}
\label{lemma:local-frame-prescribed-boundary}
There is an $\varepsilon_0>0$ so that, for every $0<\varepsilon<\varepsilon_0$,
a constant $C>0$ independent of $\varepsilon$ with the following property exists.
Let $x_0\in\partial B_1(0)$ and $0<r<1$ be fixed. 
For any map $\vec{n}\in W^{1,2}(B_1(0),\Gr_{m-2}(\R^m))$
and any ortho-normal frame
$\vec{e}=(\vec{e}_1,\vec{e}_2)\in W^{1/2,2}(\partial B_1(0)\cap B_r(x_0),\R^m\times\R^m)$
lifting $\vec{n}$:
\begin{align*}
\vec{n}=\star(\vec{e}_1\wedge\vec{e}_2)
\quad\text{on }\partial B_1(0)\cap B_r(x_0),
\end{align*}
and satisfying the estimate
\begin{align*}
\|\nabla\vec{n}\|_{L^2(B_1(0)\cap B_r(x_0))}
+[\vec{e}\,]_{W^{1/2,2}(\partial B_1(0)\cap B_r(x_0))}\leq \varepsilon,
\end{align*}
there exists an ortho-normal frame 
$\vec{g}=(\vec{g}_1,\vec{g}_2)\in W^{1,2}(B_1(0)\cap B_r(x_0),\R^m\times\R^m)$ 
lifting $\vec{n}$:
\begin{align*}
\vec{n}=\star(\vec{g}_1\wedge\vec{g}_2)
\quad\text{in }B_1(0)\cap B_r(x_0),
\end{align*}
whose trace on $\partial B_1(0)\cap B_r(x_0)$ coincides with $\vec{e}$, satisfying
 the Coulomb condition
\begin{align*}
\divop\left(\langle\nabla\vec{g_1},g_2\rangle\right)=0
\quad\text{in }B_1(0)\cap B_r(x_0),
\end{align*}
and the estimate
\begin{align*}
\|\nabla\vec{g}\|_{L^2(B_1(0)\cap B_r(x_0))}\leq C\left(
\|\nabla\vec{n}\|_{L^2(B_1(0)\cap B_r(x_0))}
+[\vec{e}\,]_{W^{1/2,2}(\partial B_1(0)\cap B_r(x_0))}\right).
\end{align*}
\end{lemma}

The proof makes use of the following elementary result.
\begin{lemma}
\label{lemma:extension-in-S^1}
Let $0<\theta_0<\pi$ be fixed and $f:(-\theta_0,\theta_0)\to\C$
be a $W^{1/2,2}$-function.
Let $F$ be its extension to $S^1\simeq(-\pi,\pi]/\sim$
by even reflection defined by:
\begin{align}
\label{eq:estension-S^1-reflection}
F(x)=\begin{cases}
f(x) &\text{if }x\in(-\theta_0,\theta_0),\\
f(m(x-\sign(x)\pi)) &\text{if }x\in (-\pi,\pi]\setminus(-\theta_0,\theta_0),
\end{cases}
\end{align}
where $m=\frac{\theta_0}{\theta_0-\pi}$.
Then  $F\in W^{1/2,2}(S^1)$ and there holds:
\begin{align*}
[F]_{W^{1/2,2}(S^1)}
\leq 2[{f}]_{W^{1/2,2}((-\theta_0,\theta_0)}.
\end{align*}
\end{lemma}

\begin{bfproof}[Proof of Lemma \ref{lemma:extension-in-S^1}]
Note first of all that we may equivalently write $F=f\circ j$, where 
$j:[-\pi,\pi]\to[-\theta_0,\theta_0]$ is defined as:
\begin{align*}
j(x)=\begin{cases}
x &\text{if }x\in(-\theta_0,\theta_0),\\
m(x-\sign(x)\pi) &\text{if }x\in (-\pi,\pi]\setminus(-\theta_0,\theta_0).
\end{cases}
\end{align*}
Using the invariance by rescaling of the $W^{1/2,2}$-seminorm and Tonelli's theorem, we see that:
\begin{align*}
[F]_{W^{1/2,2}(S^1)}^2
=2[f]_{W^{1/2,2}((-\theta_0,\theta_0))}^2
+2\int_{S^1\setminus(-\theta_0,\theta_0)}\int_{(-\theta_0,\theta_0)}
\frac{|F(x)-F(y)|^2}{|\e^{ix}-\e^{iy}|^2}\,\dd x\,\dd y.
\end{align*}
Thinking of $j$ as a diffeomorphism from $S^1\setminus[-\theta_0,\theta_0]$
to $[-\theta_0,\theta_0]$ with $j'=m$, we perform a change variable in the above integral
as follows:
\begin{align*}
&\phantom{{}={}}
\int_{S^1\setminus(-\theta_0,\theta_0)}\int_{(-\theta_0,\theta_0)}
\frac{|F(x)-F(y)|^2}{|\e^{ix}-\e^{iy}|^2}\,\dd x\,\dd y\\
&=\int_{S^1\setminus(-\theta_0,\theta_0)}\int_{(-\theta_0,\theta_0)}
\frac{|{f}(x)-{f}(j(y))|^2}{|\e^{ix}-\e^{iy}|^2}\,\dd x\,\dd y\\
&=\frac{1}{|m|}\int_{(-\theta_0,\theta_0)}\int_{(-\theta_0,\theta_0)}
\frac{|{f}(x)-{f}(\eta)|^2}{|\e^{ix}-\e^{ i(j^{-1}(\eta))}|^2}\,\dd x\,\dd \eta.\\
\end{align*}
Suppose now that $\theta_0=\pi/2$:
we have $j^{-1}(\eta)=-\eta+\sign(\eta)\pi$ and
since $x,\eta\in(-\pi/2,\pi/2)$, there holds
$|\e^{ix}-\e^{i(j^{-1}(\eta))}|\geq |\e^{ix}-\e^{i\eta}|$,
hence:
\begin{align*}
&\int_{S^1\setminus(-\pi/2,\pi/2)}\int_{(-\pi/2,\pi/2)}
\frac{|F(x)-F(y)|^2}{|\e^{ix}-\e^{iy}|^2}\,\dd x\,\dd y\\
&\leq \int_{(-\pi/2,\pi/2)}\int_{(-\pi/2,\pi/2)}
\frac{|{f}(x)-{f}(\eta)|^2}{|\e^{ix}-\e^{i\eta}|^2}\,\dd x\,\dd \eta
=[f]_{W^{1/2,2}((-\pi/2,\pi/2))}^2,
\end{align*}
So we conclude that:
\begin{align*}
[F]_{W^{1/2,2}(S^1)}^2\leq 4
[f]_{W^{1/2,2}((-\pi/2,\pi/2))}^2.
\end{align*}
For a general $0<\theta_0<\pi$, we may reduce to the case $\theta_0=\pi/2$
by using the fact that the $H^{1/2}$-seminorm is invariant with respect to
the restriction to $S^1$ of Moebius transformation of $D$.
In our particular case, the transformation we need is:
\begin{align*}
M_a(z)=\frac{z+a}{\overline{a}z+1}
\quad\text{with }a = \frac{\pi/2-\theta_0}{1-(\pi/2)\theta_0},
\quad\text{for }z\in S^1.
\end{align*}
In other words:
\begin{align*}
[F]_{W^{1/2,2}(S^1)}
=[F\circ M_a]_{W^{1/2,2}(S^1)},
\end{align*}
so we may apply the previous inequality and reach the conclusion.
This concludes the proof of the lemma.
\end{bfproof}
\begin{remark}
\label{remark:extension-in-S^1}
Lemma \ref{lemma:extension-in-S^1} holds also for domains which are conformally equivalent to $B_1(0)$.
\end{remark}
\begin{bfproof}[Proof of Lemma \ref{lemma:local-frame-prescribed-boundary}]
Without loss of generality we may assume $x_0=1$,
so that we have the identification $\partial_{B_1(0)}\cap B_r(x_0)\simeq (-\theta_0,\theta_0)$
for some $0<\theta_0<\pi$. \\
We define a map $\vec{N}\in W^{1,2}(B_1(0),\Gr_{m-2}(\R^m))$
which coincides with the given $\vec{n}$ in $B_1(0)\cap B_r(x_0)$ and has
globally controlled energy, as follows.
First,
define $\vec{n}'\in \dot{W}^{1,2}(\C,\Gr_{m-2}(\R^m))$ as the extension of $\vec{n}$ to $\C$ through even reflection:
\begin{align*}
\vec{n}'(z)=
\begin{cases}
\vec{n}(z) &\text{if }z\in B_1(0),\\
\vec{n}\left(\frac{z}{|z|^2}\right) &\text{if }z\in \C\setminus B_1(0).
\end{cases}
\end{align*}
By the conformal invariance of the Dirichlet energy, there holds
$\|\nabla\vec{n}'\|^2_{L^2(\C)}=2\|\nabla\vec{n}\|^2_{L^2(B_1(0))}$
and
\footnote{
If $I(z)=z/|z|^2$ denotes the inversion with respect to the unit circle,
then:
\begin{align*}
(\C^2\setminus D)\cap B_r(x_0))\subset I(D\cap B_r(x_0)).
\end{align*}}
\begin{align}
\label{eq:estimate-n-dash}
\|\nabla\vec{n}'\|^2_{L^2(B_r(x_0))}\leq
2\|\nabla\vec{n}\|^2_{L^2(D\cap B_r(x_0))}.
\end{align}
Consider now $\vec{n}'$ as a map from $B_r(x_0)$ and
define $\vec{N}\in \dot{W}^{1,2}(\C,\Gr_{m-2}(\R^m))$
to be its extension through even reflection:
\begin{align*}
\vec{N}(z)=
\begin{cases}
\vec{n}'(z) &\text{if }z\in B_r(x_0),\\
\vec{n}'\left(\frac{r^2}{|z-x_0|^2}(z-x_0)\right)
&\text{if }z\in\C\setminus B_r(x_0).
\end{cases}
\end{align*}
By the conformal invariance of the Dirichlet energy
and \eqref{eq:estimate-n-dash}, there holds
$\|\nabla\vec{N}\|^2_{L^2(\C)}\leq 4\|\nabla\vec{n}\|^2_{L^2(B_1(0)\cap B_r(x_0))}$,
hence a fortiori:
\begin{align*}
\|\nabla\vec{N}\|^2_{L^2(B_1(0))}
\leq 4\|\nabla\vec{n}\|^2_{L^2(B_1(0)\cap B_r(x_0))}.
\end{align*}
Consequently, by assuming $4\varepsilon_0<2\pi$,
we may invoke lemma
\ref{lemma:Helein} and find a Coulomb orthonormal frame
$\vec{f}=(\vec{f}_1,\vec{f}_2)\in W^{1,2}(B_1(0),\R^m\times\R^m)$ lifting $\vec{N}$
in $B_1(0)$ and satisfying:
\begin{align*}
\|\nabla\vec{f}\|_{L^2(B_1(0))}
\leq 2\sqrt{2} \|\nabla\vec{n}\|^2_{L^2(B_1(0)\cap B_r(x_0))}.
\end{align*}
As in the proof of lemma \ref{lemma:frame-prescribed-boundary},
the angle $\alpha_0:\partial B_1(0)\cap B_{r}(x_0)\to\R$ which rotates $\vec{f}$ 
to $\vec{e}$ is implicitly defined trough the $S^1$-valued function
$u=\langle\vec{e}_1,\vec{f}_1\rangle
-i\langle\vec{e}_2,\vec{f}_1\rangle$,
which belongs to $W^{1/2,2}((-\theta_0,\theta_0),S^1)$
and satisfies the estimate
\begin{align*}
[u]_{W^{1/2,2}((-\theta_0,\theta_0))}
&\leq 2\left(
[\vec{e}\,]_{W^{1/2,2}((-\theta_0,\theta_0))}+
[\vec{f}\,]_{W^{1/2,2}((-\theta_0,\theta_0))}
\right)
\leq C\varepsilon
\end{align*} 
for some constant $C>0$.
By means of lemma \ref{lemma:extension-in-S^1},
we may extend $u$ to $S^1=\partial B_1(0)$, thus obtaining a function 
$v\in W^{1/2,2}(\partial B_1(0),S^1)$
satisfying $[v]_{W^{1/2,2}(\partial B_1(0))}\leq 2[u]_{W^{1/2,2}((-\theta_0,\theta_0))}$.
The rest of the argument is now similar to that of the proof of lemma
 \ref{lemma:frame-prescribed-boundary}, with $v$ in place of $u$.
This concludes the proof of the lemma.
\end{bfproof}

\subsection{Singular Points of Lipschitz Immersions}
\label{sec:branch-points}
\begin{lemma}[\cite{MR3276154}, Lemma A.5]
\label{lemma:interior-branch-points}
Let $\Phi:B_1(0)\to\R^m$ be a measurable map so that,
for every $\delta>0$  $\Phi:B_1(0)\setminus B_\delta(0)\to\R^m$ defines a conformal Lipschitz
immersion with $L^2$-bounded second fundamental form.
Assume that $\Phi$ extends to a map in $W^{1,2}(B_1(0),\R^m)$ and that the  Gauss map
$\vec{n}_\Phi$ also extends to a map in $W^{1,2}(B_1(0),\Gr_{m-2}(\R^m))$.
Then $\Phi$ realises a Lipschitz map of the whole disk $B_1(0)$ and there exists a non-negative
integer $n$ and a costant $C>0$ such that
\begin{align*}
(C-o(1))|z|^{n}\leq |\nabla\Phi(z)|\leq (C+o(1))|z|^{n}
\quad\text{as }z\to 0. 
\end{align*}
\end{lemma}
\begin{lemma}
\label{lemma:boundary-singular-points}
Let $\{b_1,\ldots,b_M\}$ be points on $\partial B_1(0)$
and let $\Phi:B_1(0)\to\R^m$ be a measurable map so that,
for every $\delta>0$, $\Phi:B_1(0)\setminus \cup_{i=1}^MB_\delta (b_i)\to\R^m$
defines a conformal Lipschitz immersion with $L^2$-bounded
second fundamental form, possibly branched at  finite number of points 
$\{a_1,\ldots,a_N\}\subset D$, with $N$ independent of $\delta$.
Assume that
\begin{enumerate}
\item $\Phi$ extends to a map in $W^{1,2}(B_1(0),\R^m)$ and $\vec{n}_\Phi$ extends to a map in $W^{1,2}(B_1(0),\Gr_{m-2}(\R^m))$
\item $\log|\nabla\Phi|$ extends to a map in $W^{1,p}(B_1(0))$ for some $p>1$,
\item $\Phi|_{\partial B_1(0)} = \gamma\circ\sigma$ 
and $\vec{n}_{\Phi}|_{\partial B_1(0)}=\vec{n}_0\circ\sigma$
for some homeomorphism $\sigma$,
\end{enumerate}
where $\gamma$ is an arc-length parametrization of a closed, simple
curve $\Gamma$ in $\R^m$  of class $C^{1,1}$ and $\vec{n}_0$ is a unit-normal $m-2$ vector field
along $\Gamma$ of class $C^{1,1}$.
Then $\Phi$ extends to a weak Lipschitz immersion at every point $b_i$, $i=1,\ldots,M$.
\end{lemma}
\begin{bfproof}[Proof of Lemma \ref{lemma:boundary-singular-points}]
We call $\lambda =\log(|\nabla\Phi|/\sqrt{2})$.
It is enough to prove that, for every $i=1,\ldots, M$,
there exists some $0<s<1$ so that
\begin{align}
\label{eq:objective-singularity-removability}
\|\lambda\|_{L^{\infty}(B_1(0)\cap B_s(b_j))}<+\infty.
\end{align}

\emph{Claim 1:} For every $\varepsilon>0$,
 the coordinate ortho-normal frame of $\Phi$ denoted by
$\vec{e}=(\vec{e}_1,\vec{e}_2)$ 
extends to a map in 
$W^{1,2}(B_1(0)\setminus \cup_{j=1}^NB_\varepsilon(a_j),\R^m\times\R^m)$.\\
\emph{Proof of claim 1.}
We need to prove that $\vec{e}$ extends to a a $W^{1,2}$-map 
in a  neighbourhood of each $b_i$.
From the  relation 
\eqref{eq:decomposition-frame-energy} we have
\begin{align*}
|\nabla\vec{e}|^2
=2|\langle\nabla\vec{e}_1,\vec{e}_2\rangle|^2+|\nabla\vec{n}|^2
=2|\nabla\lambda|^2+|\nabla\vec{n}|^2
\quad\text{in }\mathcal{D}'(B_1(0)\setminus \{a_1,\ldots,a_N\})
\end{align*}
consequently, since $\nabla\lambda$ belongs to $L^{p}(B_1(0))$
and $|\nabla\vec{n}|$ belongs to  $L^2(B_1(0))$,
 we deduce that 
$|\nabla\vec{e}|$ belongs to $L^{p}(B_1(0)\setminus \cup_{j=1}^MB_\varepsilon(a_j))$
for every $\varepsilon>0$.
Hence $\vec{e}$ belongs to 
$W^{1,p}(B_1(0)\setminus \cup_{j=1}^MB_\varepsilon(a_j),\R^m\times\R^m)$ and
the trace of $\vec{e}$ on $\partial_{B_1(0)}$ is well-defined and belongs
to $W^{1/2,1-1/p}(\partial B_1(0),\R^{m}\times\R^{m})$.
Moreover,
if $\mathbf{t}$ is the unit tangent vector of $\Gamma$,
from the boundary conditions this trace is given in
complex notation by
\begin{align*}
\vec{e}_1+i\vec{e}_2=\e^{-i\theta}(\star(\mathbf{t}\wedge\vec{n}_0)+i\mathbf{t})(\sigma)
\quad\text{on }\partial B_1(0),
\end{align*}
so we see that it actually lies in
$(W^{1/2,2}\cap C^0)(\partial B_1(0),\R^{m}\times\R^{m})$.
Fix now $i=1,\ldots, M$ and choose a sufficiently small $0<r<1$
so that no branch point $a_j$, $j=1,\ldots M$, lies in $B_r(b_i)\cap D$ and
so that, thanks to lemma \ref{lemma:local-frame-prescribed-boundary},
we find an ortho-normal Coulomb frame $\vec{g}=(\vec{g}_1,\vec{g}_2)$
belonging to $W^{1,2}(B_1(0)\cap B_{r}(b_i),\R^m\times\R^m)$,
which lifts $\vec{n}_\Phi$ and whose trace on $\partial B_1(0)\cap B_r(b_i)$
coincides with that of $\vec{e}$.
If $\varphi$ denotes the angle which rotates $\vec{g}$ to $\vec{e}$ ,
from the change of Gauge formula \eqref{eq:change-of-gauge} we 
deduce that $\varphi$ is harmonic in $B_1(0)\cap B_r(b_i)$, moreover
\begin{align*}
|\nabla\varphi|\leq |\langle\nabla\vec{g}_1,\vec{g}_2\rangle|
+ |\langle\nabla\vec{e}_1,\vec{e}_2\rangle|
\quad\text{in. }\mathcal{D}'(B_1(0)\cap B_r(b_i)).
\end{align*}
Hence $\varphi\in W^{1,p}(B_1(0)\cap B_r(b_j))$ and thus
has a well-defined trace in $\partial B_1(0) \cap B_r(b_j)$ which is
is zero  by construction.
Hence $\varphi$ is smooth on $B_1(0)\cap B_{r/2}(b_i)$ and so
we deduce that $\vec{e}\in W^{1,2}(B_1(0)\cap B_{r/2}(b_j),\R^m\times\R^m)$.
Since $i=1,\dots,N$ was arbitrary, claim 1 follows.

\emph{Claim 2.} $\sigma'$ extends to a map
in $L^1(\partial B_1(0))$.\\
\emph{Proof of Claim 2.}
From the boundary conditions on $\Phi$ we have that, uniformly
on $\delta>0$,
\begin{align*}
\int_{\partial B_1(0)\setminus \cup_{i}B_\delta(b_i)}\sigma'
=\int_{\partial B_1(0)\setminus \cup_{i}B_\delta(b_i)}\e^{\lambda}|_{\partial B_1(0)}
\leq \mathcal{H}^1(\Gamma),
\end{align*}
hence, since $\sigma$ is continuous,
the classical Schwartz lemma for distributions
implies $\sigma'=\e^{\lambda}|_{\partial B_1(0)}$ extends
to a map in $L^1(\partial B_1(0))$.
This proves claim 2.

Combining claims 1 and 2, we deduce that $\lambda$ 
is a weak solution of Liouville's equation \ref{eq:Liouville}.
From claim 1, $k_g(\sigma)\e^{\lambda}-1$ belongs to $L^1(\partial B_1(0))$
hence we may find a sufficiently small $0<r<1$ so that, from lemma \ref{lemma:local-exp-int}
there holds
$\|e^{\lambda-\overline\lambda}\|_{L^p(\partial B_1(0)\cap B_{r/2}(b_i)}$
for some $p>1$.
Thanks to claim 2 and possibly reducing $r$ so that no branch point $a_j$ lies
in $D\cap B_r(b_i)$, we can invoke theorem \ref{lemma:estimate-Wente-Neumann}
and conclude that 
$\|\lambda-\overline{\lambda}\|_{L^\infty(\partial B_1(0)\cap B_{r/4}(b_i))}$ is finite,
which gives the desired estimate \ref{eq:objective-singularity-removability}.
This concludes the proof of the lemma.
\end{bfproof}

\bibliographystyle{plain}
\bibliography{/home/francesco/Documents/LaTeX/Bibliography/BiblioMath.bib}


\end{document}